\documentclass{article}
\usepackage{graphicx}
\usepackage{amsmath}
\usepackage{float}
\usepackage{color}
\usepackage{amsfonts}
\usepackage{mathtools}
\usepackage{amssymb}
\usepackage{amsthm}
\usepackage{verbatim}
\usepackage{enumitem}
\usepackage{bbold}
\usepackage{url}
\usepackage{fullpage}
\usepackage{caption}
\usepackage{subcaption}

\newcommand{\abs}[1]{\left|#1\right|}
\newcommand{\R}{\mathbb{R}}
\newcommand{\N}{\mathbb{N}}
\newcommand{\Z}{\mathbb{Z}}
\newcommand{\norm}[1]{||#1||}

\newtheorem{proposition}{PROPOSITION}[section]
\newtheorem{corollary}[proposition]{COROLLARY}
\newtheorem{remark}[proposition]{REMARK}
\newtheorem{lemma}[proposition]{LEMMA}
\newtheorem{theorem}[proposition]{THEOREM}
\newtheorem{definition}[proposition]{DEFINITION}

\usepackage{xcolor}
\definecolor{darkgreen}{RGB}{0,100,0}

\title{Null controllability of the 1D heat equation with interior inverse square potential}
\date{\today}
\author{
Pierre Lissy\thanks{CERMICS, ENPC, Institut Polytechnique de Paris,  Marne-la-Vall\'ee, France. (\texttt{pierre.lissy@enpc.fr})}
\and Tanguy Lourme \thanks{CERMICS, ENPC, Institut Polytechnique de Paris, Marne-la-Vall\'ee, France. (\texttt{tanguy.lourme@gmail.com)}}}

\numberwithin{equation}{section}

\begin{document}

\maketitle
\begin{abstract}
This paper aims to answer an open problem posed by Morancey in 2015 concerning the null controllability of the heat equation on $(-1,1)$ with an internal inverse square potential located at $x=0$. For the range of singularity under study, after having introduced a suitable self-adjoint extension that enables to transmit information from one side of the singularity to another, we prove null-controllability in arbitrary small time, firstly with an internal control supported in an arbitrary measurable set of positive measure, secondly with a boundary control acting on one side of the boundary. Our proof is mainly based on a precise spectral study of the singular operator together with some recent  refinements of the Fattorini-Russell moment method. This in particular requires to use some fine (and sometimes new) properties of Bessel functions and their zeros.
\end{abstract}



\tableofcontents

\section{Introduction}
\subsection{Statement of the problem and main result}
 This article is mainly dedicated to investigating the null-controllability properties of the following singular parabolic equation with inverse square potential
\begin{equation} \left \{
\begin{aligned}\partial_tf(t,x)-\partial_{xx}f(t,x)+\dfrac{c}{x^2}f(t,x)&=u(t, x)\mathbb{1}_\omega(x),& (t,x)\in(0,T)\times(-1, 1),\\f(t,-1)=f(t, 1)&=0,&t\in(0,T),\\
f(0,x)&=f^0(x),&x\in(-1, 1),
\end{aligned} \right .
\label{1}
\end{equation}
where $c\in \mathbb R$, $T>0$ is a time horizon, $\omega$ is a measurable subset of $(-1,1)$ assumed to have positive Lebesgue measure, $\mathbb{1}_\omega$ denotes the characteristic function of $\omega$, $f^0\in L^2(-1,1)$ and $u \in L^2\left((0,T)\times (-1,1)\right)$. 

Such potentials naturally appear in the context of combustion theory (\cite{BE}) or quantum mechanics (\cite{BG2}).

It is proved in \cite{BG} that \eqref{1} is ill-posed in $L^2(0,1)$ (with Dirichlet boundary conditions in $0$ and $1$) if $c<-\frac{1}{4}$ (see also \cite{VZ2}). On  $L^2(-1,1)$, we can show that \eqref{1} is ill-posed, in the sense that there does not exist any self-adjoint extension  of $\partial_{xx}-\dfrac{c}{x^2}\mathrm{Id}$, posed on $C^\infty_0((-1,1)\setminus \{0\})$, generating a $C^0$-semigroup, if $c<-\frac{1}{4}$. The proof of this fact is postponed to Appendix \ref{appen}.
 
Furthermore, from \cite[Chapter X]{RS}, we know that for $c\geqslant \frac{3}{4}$, the operator $-\partial_{xx}+\dfrac{c}{x^2}\mathrm{Id}$ is essentially self-adjoint on $C^\infty_0(\mathbb R\setminus\{0\})$, meaning that no information passes through the singularity $x=0$ (see \textit{e.g.} \cite[pp. 1744-1745]{BL} for an interesting discussion on this subject). In particular, in this context, it is not possible to control with a measurable set $\omega$ that lies only on one part of the singularity $x=0$. Indeed,  from \cite[Corollary 1.3]{BL}, in the absence of a control, an initial condition supported in $(-1,0)$, for instance, leads to a solution of \eqref{1} that is also supported  in $(-1,0)$ for any $t\geqslant 0$ (the same reasoning applies on $(0,1)$). By the Duhamel formula, in absence of an initial condition, a control acting on a measurable set of $(0,1)$ also leads to a solution of \eqref{1} that is also supported  in $(0,1)$ for any $t\geqslant 0$. It is then clear that an initial condition supported in $(-1,0)$ cannot be controlled from a measurable set included in $(0,1)$ (see also Remark \ref{rem:decoupled-not-interesting} for a similar reasoning).

For these reasons, from now on, we will restrict to the coefficient $c$ lying in  $\left(-\frac14,+\frac34\right)$ (the case $c=\frac{1}{4}$ cannot be  handled directly with the techniques developed in this article, see Section \ref{sec:c} for more explanations).

\noindent Throughout this paper, for the sake of consistency with the existing literature, the coefficient of the singular potential is parameterized as \(c = c_\nu\), where
\begin{equation}
c_{\nu}:=\nu^{2}-\frac{1}{4}, \quad  \nu \in \mathbb R.
    \label{1.2}
\end{equation}
Then, $c \in  \left(-\frac14,+\frac34\right)$  if and only if $|\nu| \in (0,1)$, and in this case, $c_{\nu}=c_{-\nu}$, which justifies to consider only positive $\nu$.
Our main result is the following null-controllability result in arbitrary small time.
\begin{theorem} 
      Assume that $\nu \in (0,1)$. For any $f^0 \in L^2(-1,1)$ and $T>0$,  there exists a control  $u \in L^2\left((0,T)\times (-1,1)\right)$ { such that }the corresponding solution $f$ to \eqref{1} verifies $f(T) = 0$ in $L^2(-1,1)$.
    \label{theoremNullControlIntern}
\end{theorem}
We also have a version of this Theorem with a boundary control on one side of the boundary, which requires some nontrivial intermediate results and more notations. For the sake of simplicity, we chose to postpone it at the end of the paper (see Theorem \ref{mainb}).
\begin{remark}
In the previous Theorem, the solutions of \eqref{1} have to be understood in a sense explained in Section \ref{sec:WP}. In particular, in the case where $c_\nu=0$, we obtain a particular self-adjoint extension of the Laplace operator on $C^\infty_0((-1,0)\cup(0,1))$ that does not correspond to the usual Dirichlet-Laplace operator on $(-1,1)$. In  Remark \ref{posrem} of Appendix \ref{appex}, in the case $c_\nu=0$, we explain how to create appropriate transmission conditions that enable to recover the usual Dirichlet-Laplace operator on $(-1,1)$. However, for $c_\nu\not =0$, it is not clear for us whether this extension is non-negative or not, so we chose not to treat this extension in more details.
\end{remark}
\begin{remark}
Here, $\omega$ may be any measurable set of positive measure, and not only a nonempty open set. As we will see during the proof, extending to this more general case is not completely straightforward (see Remark \ref{dif:NL}). For seminal works on controllability by measurable sets for the heat equation, one can read \cite{AE,PW,AEWZ}.

\end{remark}
 The controllability properties of heat equations with inverse-square potentials have attracted some attention over the last fifteen years. In dimension larger than three, the case of an internal inverse square singularity with inverse square potential and internal control  was first treated in  \cite{VZ} in a particular geometry and then extended in \cite{Erv} (see also \cite{Z} for the study of an optimal control problem in this setting). The case of a variable diffusion has been treated in \cite{QL}. The case of approximate boundary control has been studied in \cite{sh2}, again in dimension greater than three. The case of an inverse square singularity located at the boundary and internal control was first treated in any space dimension in \cite{caz}. A related result is also \cite{BZ}, involving the distance to the boundary in dimension larger than three. 

 Concerning specifically the case of the one-dimensional heat equation with inverse square potential at one point of the boundary, the question of the cost of controllability with a boundary control located at the other point of the boundary has been studied in \cite{MV}, whereas \cite{B} proved a result where the boundary control acts directly at the singularity location. The case of mixed degenerate diffusion and singular potentials is treated in \cite{FS,V}, and a model with memory is studied in \cite{AFS}.

 Remark that in all the literature cited above, in the 1D case, the case where the singularity lies \emph{inside} the domain has not been addressed, except in \cite[Theorem 1.5]{Morancey2015}, where the author proved an approximate controllability result in this context. 
In \cite{Morancey2015}, the author raised the open question of obtaining a null-controllability result in the context of equation \eqref{1}. The goal of the present article is to give a positive answer to this question.

That this question has not been investigated earlier comes from the fact that this situation is much more difficult from a theoretical point of view. Indeed, the domain $(-1,1)$ is split into two subintervals by the singularity of the inverse square potential. 
A natural (but naive) idea is to treat separately the Laplace operator with inverse square potential on $(-1,0)$ and on $(0,1)$, and to use on each of these intervals the domain coming from \cite{V}. In other words, to consider the operator $\Delta_0$ with domain 
\[
\{\, f \in L^{2}(-1,1) ;
f_{|[0,1]} \in H^2_{\mathrm{loc}}((0,1]) \cap H^1_0(0,1),\;
f_{|[-1,0]} \in H^2_{\mathrm{loc}}([-1,0)) \cap H^1_0(-1,0),\ 
-\partial_{xx}^2 f + \frac{c_\nu}{x^2} f \in L^2(-1,1)\}.
\]

The operator $\Delta_{0}$ is non-negative and self-adjoint, so the associated heat equation is well-posed by standard semigroup theory. 
However, in this framework, it is impossible to control equation~\eqref{1} by acting on only one side of the origin, since the corresponding heat equation has completely decoupled dynamics on $(-1,0)$ and $(0,1)$. This will be explained in more details in Appendix \ref{appex}, in particular in Remark \ref{rem:decoupled-not-interesting}.

 In fact, we need to consider another self-adjoint extension of $A_\nu$  on $C^\infty_0((-1,1) \setminus{\{0\}})$. Of course, what is important is to specify what to prescribe at the boundary $x=0$. As we will see in Section \ref{sec:WP}, we are able to introduce an appropriate self-adjoint extension, with appropriate transmission conditions, that enables to pass information through the singularity. The main difficulty is then to give an appropriate spectral decomposition of this self-adjoint extension. In particular, we will give some very precise estimates on the eigenvalues of $A_\nu$ and also on the eigenfunctions of $A_\nu$.

The plan of the paper is as follows. In Subsection \ref{sec:WP}, we investigate the well-posedness of \eqref{1}. In Subsection \ref{ss:b}, we give some already-known results on Bessel functions, together with new technical results that are needed for our study. In Section \ref{s:d}, we study the spectral decomposition of the singular elliptic operator under consideration. In Section \ref{s:as}, we give some precise results on the asymptotic behavior of the eigenvalues of our singular operator. In Section \ref{s:nc}, we prove our main Theorem \ref{theoremNullControlIntern} and give an extension to the case of a boundary control. To conclude, in Section \ref{sec:c},  we present several open problems related to the present study.

\subsection{Well-posedness}
\label{sec:WP}
\noindent
We consider the following homogeneous problem: 
\begin{equation}\left \{
 \begin{aligned}
 \partial_t f - \partial_{xx} f + \frac{c_{\nu}}{x^2} f &= 0, &  (t, x) \in(0, T) \times (-1, 1),\\
f(t, -1) = f(t, 1) &= 0, & t \in (0, T),\\f(0,x)&=f^0(x),&x\in(-1, 1).
 \end{aligned} \right .
     \label{1.3}
\end{equation}
We begin by addressing the well-posedness of (\ref{1.3}). The elliptic differential operator under consideration is for a fixed $\nu\in(0,1)$
\[
A_\nu f:=-\partial_{xx}f+\frac{c_{\nu}}{x^{2}}f.
\]
It is clear that \(A_\nu\) is well-defined on \(C_0^\infty((-1,0)\cup(0,1))\). As already explained, for $\nu \in (0,1)$, $A_\nu$ is not essentially self-adjoint on \(C_0^\infty((-1,0)\cup(0,1))\). Further details on the self-adjoint extensions of $A_\nu$ are available in Appendix \ref{appex}. We now specify the self-adjoint extension to be employed here, whose relevance is discusses in Remark \ref{rem:pos}.

\noindent As in \cite[Section 2.1]{Morancey2015}, we introduce:
\begin{equation}
    \tilde{H}^2_0(-1,1):=\left\{f\in H^2(-1,1), \ f(0)=f'(0)=0\right\},
    \label{H20}
\end{equation}
\begin{equation}
    \mathcal{F}_s^\nu:= \Bigg\{ f\in L^2(-1,1), \ \exists c_1^+, c_1^-, c_2^+, c_2^- \in \R, \quad f(x)=
    \begin{cases}
    c_1^-|x|^{\nu+\frac12}+c_2^-|x|^{-\nu+\frac12} \ & \mathrm{on} \ (-1,0),\\
    c_{1}^{+}|x|^{\nu+\frac12}+c_{2}^{+}|x|^{-\nu+\frac12}\ & \mathrm{on} \ (0,1).
    \end{cases}\Bigg\}.
    \label{Fs}
\end{equation}
\noindent Notice that for any \(f_s \in \mathcal{F}_s^\nu\), $f_s \in C^\infty((-1,0)\cup(0,1))$ and
\begin{equation}  
(-\partial_{xx}+\frac{c_\nu}{x^2})f_s(x)=0, \quad \forall x\in(-1,0)\cup(0,1).
    \label{1.4}
\end{equation}
\noindent The domain of $A_\nu$  is defined as
\begin{equation}
\begin{split}
    D(A_\nu)&:=\Big\{f=f_{r}+f_{s}; \ f_{r}\in\tilde{H}_{0}^{2}(-1,1),\ f_{s}\in\mathcal{F}_s^\nu \ \mathrm{such\ that\ } f(-1)=f(1)=0,\\ &c_{1}^{-}+c_{2}^{-}+c_{1}^{+}+c_{2}^{+}=0, \ \mathrm{and} \ \left(\nu+\frac{1}{2}\right)c_{1}^{-}+\left(-\nu+\frac{1}{2}\right)c_{2}^{-}=\left(\nu+\frac{1}{2}\right)c_{1}^{+}+\left(-\nu+\frac{1}{2}\right)c_{2}^{+}\Big\},
\end{split}
    \label{1.5}
\end{equation}
which is a linear subspace of \(L^2(-1, 1)\). By \cite[Remark 2.3]{Morancey2015}, we have $A_\nu(D(A_\nu)) \subset L^2(-1,1)$, confirming that $(A_\nu, D(A_\nu))$ is indeed an unbounded operator on $L^2(-1,1)$, which is moreover densely defined since 
\[
\Big\{\phi \in C_0^\infty((-1,0) \cup (0,1)), \text{ extended at } 0 \text{ by } \phi(0) = 0\Big\} \subset \tilde{H}_0^2(-1,1) \subset D(A_\nu).
\]
To conclude, using \cite[Proposition 2.2]{Morancey2015},  $(A_\nu, D(A_\nu))$ is self-adjoint on $L^2(-1,1)$, and for any $f \in D(A_\nu)$, the following inequality holds:

\begin{equation}
    \langle A_\nu f, f \rangle \geqslant \min\{1, 4\nu^2\} \int_{-1}^1 \left(\partial_x f_r(x)\right)^2 \, \mathrm{d}x.
    \label{1.6}
\end{equation}
\begin{remark}
    For \(\nu \in (0,1)\), it is easy to verify that we have a direct sum \(\tilde{H}_{0}^{2}(-1,1)\oplus\mathcal{F}_s^\nu
\). Since \(D(A_\nu)\subset\tilde{H}_{0}^{2}(-1,1)\oplus\mathcal{F}_s^\nu\), the unique decomposition of functions in \(D(A_\nu)\)  given by $f=f_r+f_s$, where $f\in D(A_\nu),\,f_r\in \tilde{H}_{0}^{2}(-1,1),\,f_s\in\mathcal{F}_s^\nu$ will be referred to as the decomposition into the regular part \(f_r\) and the singular part \(f_s\). The conditions imposed on the coefficients  of the singular part in (\ref{1.5}) will be referred to as transmission conditions.
    \label{rqUniqueDecomposition}
\end{remark}
\begin{remark}  On \(\mathcal{F}_s^\nu\), the operator \(A_\nu\) acts independently on \((0, 1)\) and \((-1, 0)\). Thus, the symbol \(\partial_{xx}\) should not be interpreted as a derivative in the distributional sense over \(\Omega = (-1, 1)\). Instead, it operates separately on \((-1, 0)\) and \((0, 1)\). In other words, \(\partial_{xx}\) has to be understood as a distributional derivative on \(\Omega = (-1, 0) \cup (0, 1)\). Therefore, from (\ref{1.4}), we have \(A_\nu f_s = 0\).
\end{remark}

\noindent We have already remarked that $(A_\nu,D(A_\nu))$ is a densely defined self-adjoint and monotone operator (by \eqref{1.6}). From \cite[Corollary 2.4.8]{CazenaveHaraux1990}, we deduce that $-A_\nu$ is a $m$-dissipative operator, so that applying the Hille-Yosida theorem for self-adjoint monotone operators (see, for instance, \cite[Theorem 3.2.1]{CazenaveHaraux1990}), we obtain that 
$-A_\nu$ generates a strongly continuous semigroup, denoted by $(e^{-A_\nu t})_{t \geqslant 0}$. In other words, the homogeneous problem \eqref{1.3} is well-posed: for any $f^0 \in L^2(-1,1)$, there exists a unique  
$$f \in C^{0}([0,+\infty), L^{2}(-1,1)) \cap C^{0}((0,+\infty), D(A_\nu)) \cap C^{1}((0,+\infty), L^{2}(-1,1))$$
such that $f$ verifies \eqref{1.3}.

\noindent Our next goal is to define the notion of a solution to the non-homogeneous problem (\ref{1}).

\noindent Let us introduce the operator
\[
B\ :\ L^2(-1,1) \to L^2(-1,1), \quad g \mapsto \mathbb{1}_\omega g.
\]
Since $B \in \mathcal{L}_c(L^2(-1,1))$, the non-homogeneous problem (\ref{1}) is well-posed. In other words, there exists a unique weak solution (defined for in the sense of transpositions as in \cite[Section 2.3]{Coron2007} for instance) to (\ref{1}) in $C^0([0,T], L^2(-1,1))$, that is ``explicitly'' given by

\[
f = e^{-A_\nu t}f^0 + \int_0^t e^{-A_\nu (t-s)} \mathbb{1}_\omega u(s) \, \mathrm{d}s, \quad \forall t \in [0,T].
    \]

\subsection{Generalities on Bessel functions}
\label{ss:b}
\noindent
In all what follows, we consider some coefficient  ${\tilde{\nu}}\in \mathbb R$ (sometimes, $\tilde \nu$ will play the role of $\nu$ in \eqref{1.2}). The Bessel functions of the first kind, denoted by $J_{\tilde{\nu}}$, are solutions to the differential equation given by
\begin{equation}
    x^2 \frac{d^2 y}{\mathrm{d}x^2} + x \frac{dy}{\mathrm{d}x} + (x^2 - {\tilde{\nu}}^2) y = 0, \quad \forall x>0.
    \label{besseleq}
\end{equation}

\noindent The Bessel function of the first kind $J_{\tilde{\nu}}(x)$ can be defined for all $x>0$ by its power series expansion
\begin{equation}
    J_{\tilde{\nu}}(x) = \sum_{k=0}^{+\infty} \frac{(-1)^k}{k! \, \Gamma(k + {\tilde{\nu}} + 1)} \left(\frac{x}{2}\right)^{2k + {\tilde{\nu}}},\quad\forall x>0,
    \label{besselseries}
\end{equation}
where $\Gamma(z)$ is the Euler Gamma function (by convention, $\frac{1}{\Gamma(x)}=0$ if $x$ is a non-positive integer). For all ${\tilde{\nu}}\in\mathbb R$ and all $n\in\N\setminus\{0\}$, we denote by $j_{{\tilde{\nu}}, n}$ the $n$-th positive real zero of $J_{\tilde{\nu}}$. We refer the reader to  \cite{Watson1944,Luke1962} for more details on Bessel functions.

\vspace{0.2cm}

In Appendix \ref{AB},  we gather several results on Bessel functions that will be useful in the proofs of this section.

\vspace{0.2cm}

The first result that we need is the following inequality, which we have not found in the literature and may be of independent interest.

\begin{lemma}Assume that ${\tilde{\nu}}\in(0,1)$. Then, 
\[    J_{\tilde{\nu}}(x)J_{-{\tilde{\nu}}}(x)<\frac{\sin({\tilde{\nu}}\pi)}{{\tilde{\nu}}\pi}, \quad \forall x>0.
\]
\label{lemmeInequality}
\end{lemma}
\noindent \textbf{\textit{Proof of Lemma \ref{lemmeInequality}}}. Let $x>0$. We introduce
\[
\varphi_x\ :\ (0,1] \longrightarrow \R, \quad {\tilde{\nu}} \mapsto \frac{\sin({\tilde{\nu}}\pi)}{{\tilde{\nu}}\pi}-J_{\tilde{\nu}}(x)J_{-{\tilde{\nu}}}(x).
\]
From \cite[p. 150, (1)]{Watson1944}, we know that
\[
            J_{\tilde{\nu}}(z)J_{-{\tilde{\nu}}}(z) = \frac{2}{\pi}\int_0^{\dfrac{\pi}{2}}J_0(2z\cos(\theta))\cos(2{\tilde{\nu}}\theta)\mathrm{d}\theta, \quad \forall {\tilde{\nu}}\in\R, \ \forall z > 0. 
\]
Therefore, $\varphi_x$ is differentiable on $(0,1]$ and
\[
\begin{split}
{\varphi_x}'({\tilde{\nu}})&=\frac{\cos({\tilde{\nu}}\pi)}{{\tilde{\nu}}}-\frac{\sin({\tilde{\nu}}\pi)}{\pi{\tilde{\nu}}^2}+\frac{2}{\pi}\int_0^{\dfrac{\pi}{2}}J_0(2x\cos(\theta))2\theta\sin(2{\tilde{\nu}}\theta)\mathrm{d}\theta,\quad \forall{\tilde{\nu}}\in(0,1].
\end{split}
\]
From the integral representation given in \cite[p. 176, (4)]{Watson1944}, we have 
$$J_0(x)=\frac{1}{\pi}\int_0^{\pi}\cos(x\sin(\theta))d\theta.$$
We deduce that  $|J_0(x)|<1$ for $x>0$.
Thus, we have
\[
\left |\frac{2}{\pi}\int_0^{\dfrac{\pi}{2}}J_0(2x\cos(\theta))2\theta\sin(2{\tilde{\nu}}\theta)\mathrm{d}\theta \right|\leqslant  \frac{2}{\pi}\int_0^{\dfrac{\pi}{2}}\abs{J_0(2x\cos(\theta))}2\theta\sin(2{\tilde{\nu}}\theta)\mathrm{d}\theta< \frac{2}{\pi}\int_0^{\dfrac{\pi}{2}}2\theta|\sin(2{\tilde{\nu}}\theta)|\mathrm{d}\theta,
\]
 and since ${\tilde{\nu}}\in(0,1)$, we have $\sin(2{\tilde{\nu}}\theta) > 0$ for all $\theta \in (0, \frac{\pi}{2})$. Furthermore,
\[
\frac{2}{\pi}\int_0^{\dfrac{\pi}{2}}2\theta\sin(2{\tilde{\nu}}\theta)\mathrm{d}\theta=\frac{2}{\pi}\left[-\frac{\theta\cos(2{\tilde{\nu}}\theta)}{{\tilde{\nu}}}+\frac{\sin(2{\tilde{\nu}}\theta)}{2{\tilde{\nu}}^2}\right]^{\dfrac{\pi}{2}}_0=-\frac{\cos({\tilde{\nu}}\pi)}{{\tilde{\nu}}}+\frac{\sin({\tilde{\nu}}\pi)}{\pi{\tilde{\nu}}^2}.
\]
Gathering the previous computations, we deduce that  ${\varphi_x}'({\tilde{\nu}})< 0$ for any ${\tilde{\nu}}\in(0,1]$, so $\varphi_x$ is  decreasing on $(0,1]$. Hence,
\[\forall {\tilde{\nu}}\in(0,1),\quad \varphi_x({\tilde{\nu}})>\varphi_x(1)=-J_1(x)J_{-1}(x).\]
From \cite[(2), p.15]{Watson1944}, we have $J_{-1}(x)=-J_1(x)$. We deduce that 
\[\forall {\tilde{\nu}}\in(0,1),\quad \varphi_x({\tilde{\nu}})>\varphi_x(1)=J_1(x)^2\geqslant0,
\]
which concludes the proof.
\qed

The next lemma gives an interlacing property on the zeros of $J_{\tilde{\nu}}$ and $J_{-{\tilde{\nu}}}$.

\begin{lemma} Assume that ${\tilde{\nu}}\in(0, 1)$. The real positive zeros of $J_{\tilde{\nu}}$ and $J_{-{\tilde{\nu}}}$ are strictly interlaced, and
\[
j_{-{\tilde{\nu}},1}<j_{{\tilde{\nu}},1}<j_{-{\tilde{\nu}},2}<j_{{\tilde{\nu}},2}<...<j_{-{\tilde{\nu}},n}<j_{{\tilde{\nu}},n}<...
\]
    \label{lemmaInterlaced}
\end{lemma}

\noindent \textbf{\textit{Proof of Lemma \ref{lemmaInterlaced}}}. Since ${\tilde{\nu}} \in (0,1)$, we have ${\tilde{\nu}}>-1$ and $-{\tilde{\nu}}>-1$, and we know from \cite[Theorem 3 and Remarks]{PA} that the positive real zeros of $J_{\tilde{\nu}}$ are interlaced. Thus, it suffices to show that $j_{-{\tilde{\nu}},1} < j_{{\tilde{\nu}},1}$.

\noindent To prove this, we will analyze the function $g : x \mapsto \frac{J_{\tilde{\nu}}(x)}{J_{-{\tilde{\nu}}}(x)}$ and show that $g$ is increasing on $(0, j_{-{\tilde{\nu}},1})$, and that $\lim_{x\to 0^+} g(x) = 0$. Consequently, we will have $g(x) > 0$ on $(0, j_{-{\tilde{\nu}},1})$, which implies that $j_{-{\tilde{\nu}},1} < j_{{\tilde{\nu}},1}$. By Lemma \ref{lemmaWronskien}, we know that 
\[
\begin{aligned}
    \quad g'(x)&=\frac{{J_{\tilde{\nu}}}'(x)J_{-{\tilde{\nu}}}(x)-{J_{\tilde{\nu}}}(x){J_{-{\tilde{\nu}}}}'(x)}{{J_{-{\tilde{\nu}}}}^2(x)}=\frac{-W(J_{\tilde{\nu}}, J_{-{\tilde{\nu}}})(x)}{{J_{-{\tilde{\nu}}}}^2(x)}=\frac{2\sin({\tilde{\nu}}\pi)}{\pi x {J_{-{\tilde{\nu}}}}^2(x)} > 0, \quad \forall x \in (0, j_{-{\tilde{\nu}},1}).
\end{aligned}
\]
We have shown that $g$ is increasing on $(0, j_{-{\tilde{\nu}},1})$. Finally, since $\lim_{x\to 0^+} J_{\tilde{\nu}}(x) = 0$ and $\lim_{x\to 0^+} J_{-{\tilde{\nu}}}(x) = +\infty$, it follows that $\lim_{x\to 0^+} g(x) = 0$. We conclude as stated above.
\qed

    \begin{lemma} Let ${\tilde{\nu}}\in\R$.
\[
    J_{\tilde{\nu}}(x) \underset{x\to+\infty}{=} \mathcal{O}\left(\frac{1}{x^{1/2}}\right), \quad
    {J_{\tilde{\nu}}}'(x) \underset{x\to+\infty}{=} \mathcal{O}\left(\frac{1}{x^{1/2}}\right).
    \]
        \label{corrAsymptBesselAndDerivatives}
    \end{lemma}

\noindent \textbf{\textit{Proof of Corollary \ref{corrAsymptBesselAndDerivatives}}}. The first asymptotic behavior is directly obtained from Lemma \ref{lemmaExpansion}. The second asymptotic behavior follows immediately from Lemma \ref{lemmaDerivateBessel} and the first asymptotic behavior.
\qed

\begin{lemma} 
Assume ${\tilde{\nu}}\in(0,1)$.

    \begin{subequations}
    \setlength{\jot}{10pt}
    \renewcommand{\theequation}{\theparentequation.\alph{equation}}
    \begin{align}
        &\int^1_0xJ_{\pm{\tilde{\nu}}}(ax)^2\mathrm{d}x = \frac12\left(\left(1-\frac{{\tilde{\nu}}^2}{a^2}\right)J_{\pm{\tilde{\nu}}}(a)^2+{J_{\pm{\tilde{\nu}}}}'(a)^2\right),&\quad \forall a>0,\label{intSame-NuSameA}\\
        &\int^1_0xJ_{{\tilde{\nu}}}(ax)J_{-{\tilde{\nu}}}(ax)\mathrm{d}x = \frac12\left(\left(1-\frac{{\tilde{\nu}}^2}{a^2}\right)J_{\tilde{\nu}}(a)J_{-{\tilde{\nu}}}(a)+{J_{\tilde{\nu}}}'(a){J_{-{\tilde{\nu}}}}'(a)\right)+\frac{{\tilde{\nu}} \sin (  {\tilde{\nu}}\pi)}{\pi  a^2}, &\quad \forall a>0, \label{intNotSameNuSameA}\\
        &\int^\beta_\alpha xJ_{\pm{\tilde{\nu}}}(ax)^2\mathrm{d}x = \frac12 \beta^2\left(\left(1-\frac{{\tilde{\nu}}^2}{a^2\beta^2}\right)J_{\pm{\tilde{\nu}}}(a\beta)^2+{J_{\pm{\tilde{\nu}}}}'(a\beta)^2\right)&\quad \forall \alpha, \beta > 0, \ \forall a>0, \nonumber\\
        & \qquad \qquad \qquad \quad \ \  -\frac12 \alpha^2\left(\left(1-\frac{{\tilde{\nu}}^2}{a^2\alpha^2}\right)J_{\pm{\tilde{\nu}}}(a\alpha)^2+{J_{\pm{\tilde{\nu}}}}'(a\alpha)^2\right), \label{intSame-NuSameA - alpha-beta}\\
        &\int^\beta_\alpha xJ_{{\tilde{\nu}}}(ax)J_{-{\tilde{\nu}}}(ax)\mathrm{d}x = \frac12\beta^2\left(\left(1-\frac{{\tilde{\nu}}^2}{a^2\beta^2}\right)J_{\tilde{\nu}}(a\beta)J_{-{\tilde{\nu}}}(a\beta)+{J_{\tilde{\nu}}}'(a\beta){J_{-{\tilde{\nu}}}}'(a\beta)\right)&\quad \forall \alpha, \beta > 0, \ \forall a>0, \nonumber\\
        &\qquad \qquad \qquad \qquad \qquad -\frac12\alpha^2\left(\left(1-\frac{{\tilde{\nu}}^2}{a^2\alpha^2}\right)J_{\tilde{\nu}}(a\alpha)J_{-{\tilde{\nu}}}(a\alpha)+{J_{\tilde{\nu}}}'(a\alpha){J_{-{\tilde{\nu}}}}'(a\alpha)\right). \label{intNotSameNuSameA - alpha-beta}
    \end{align}
    \end{subequations}

    \label{propintegrals}
\end{lemma}

\noindent \textbf{\textit{Proof of Lemma \ref{propintegrals}.}}
\noindent \textbf{Proof of (\ref{intSame-NuSameA})} \cite[(4), p. 255]{Luke1962} states that for any $a>0$,
\[
\begin{split}
    \int^1_0xJ_{\pm{\tilde{\nu}}}(ax)^2\mathrm{d}x=\frac12\left(\left(1-\frac{{\tilde{\nu}}^2}{a^2}\right)J_{\pm{\tilde{\nu}}}(a)^2+{J_{\pm{\tilde{\nu}}}}'(a)^2\right)-\lim_{x\to0^+}\frac12 x^2\left(\left(1-\frac{{\tilde{\nu}}^2}{a^2x^2}\right)J_{\pm{\tilde{\nu}}}(ax)^2+{J_{\pm{\tilde{\nu}}}}'(ax)^2\right).
\end{split}
\]
Moreover, by \eqref{besselseries},
\[
 J_{\pm{\tilde{\nu}}}(ax)\underset{x \to 0^+}{\sim}\frac{1}{\Gamma(\pm{\tilde{\nu}}+1)}\left(\frac{ax}{2}\right)^{\pm{\tilde{\nu}}},
\]
so
\[
xJ_{\pm{\tilde{\nu}}}(ax)\underset{x \to 0^+}{\sim}\frac{1}{\Gamma(\pm{\tilde{\nu}}+1)}\left(\frac{a}{2}\right)^{\pm{\tilde{\nu}}}x^{\pm{\tilde{\nu}}+1}\underset{x\to 0^+}\longrightarrow0 \quad \mathrm{because} \ {\tilde{\nu}}\in(0,1).
\]
Therefore, $\lim_{x\to0^+}x^2J_{\pm{\tilde{\nu}}}(ax)^2=0$, and we have
\[
\lim_{x\to0^+}\frac12 x^2\left(\left(1-\frac{{\tilde{\nu}}^2}{a^2x^2}\right)J_{\pm{\tilde{\nu}}}(ax)^2+{J_{\pm{\tilde{\nu}}}}'(ax)^2\right)=\lim_{x\to0^+}\frac{1}{2a^2}\left(-{\tilde{\nu}}^2J_{\pm{\tilde{\nu}}}(ax)^2+(ax)^2{J_{\pm{\tilde{\nu}}}}'(ax)^2\right).
\]
Therefore, to conclude the proof of (\ref{intSame-NuSameA}), we only need to prove that
\[
\lim_{x\to0^+}x^2{J_{\pm{\tilde{\nu}}}}'(x)^2 -{\tilde{\nu}}^2J_{\pm{\tilde{\nu}}}(x)^2=0.
\]
By examining the series expansion of ${J_{\pm{\tilde{\nu}}}}'$ obtained from \eqref{besselseries}, we have
\[
{J_{\pm{\tilde{\nu}}}}'(x) \underset{x\to0^+}{=}\frac{\pm{\tilde{\nu}}}{\Gamma(\pm{\tilde{\nu}}+1)}\frac{1}{2^{\pm{\tilde{\nu}}}}x^{\pm{\tilde{\nu}}-1} + \mathcal{O}(x^{1\pm{\tilde{\nu}}}).
\]
Therefore, we have
\begin{equation}
    x^2{J_{\pm{\tilde{\nu}}}}'(x)^2 \underset{x\to0^+}{=}\left[\frac{{\tilde{\nu}}}{\Gamma(\pm{\tilde{\nu}}+1)}\frac{1}{2^{\pm{\tilde{\nu}}}}\right]^2x^{\pm 2{\tilde{\nu}}} + \mathcal{O}(x^{2\pm2{\tilde{\nu}}}) \quad \mathrm{because} \ {\tilde{\nu}}\in(0,1),
    \label{term1}
\end{equation}
and by examining the series expansion of $J_{\pm{\tilde{\nu}}}$ given in \eqref{besselseries}, we have
\[
{\tilde{\nu}} J_{\pm{\tilde{\nu}}}(x)\underset{x\to0^+}{=}\frac{{\tilde{\nu}}}{\Gamma(\pm{\tilde{\nu}}+1)}\frac{1}{2^{\pm{\tilde{\nu}}}}x^{\pm{\tilde{\nu}}} + \mathcal{O}(x^{2\pm{\tilde{\nu}}}).
\]
so
\begin{equation}
        {\tilde{\nu}}^2 J_{\pm{\tilde{\nu}}}(x)^2 \underset{x\to0^+}{=}\left[\frac{{\tilde{\nu}}}{\Gamma(\pm{\tilde{\nu}}+1)}\frac{1}{2^{\pm{\tilde{\nu}}}}\right]^2x^{\pm2{\tilde{\nu}}} + \mathcal{O}(x^{2\pm2{\tilde{\nu}}}) \quad \mathrm{because} \ {\tilde{\nu}}\in(0,1).
        \label{term2}
\end{equation}

\noindent Finally, from (\ref{term1}) and (\ref{term2}), we obtain
\[
x^2{J_{\pm{\tilde{\nu}}}}'(x)^2-{\tilde{\nu}}^2 J_{\pm{\tilde{\nu}}}(x)^2\underset{x\to0^+}{=}\mathcal{O}(x^{2\pm2{\tilde{\nu}}})\underset{x\to0^+}{\longrightarrow}0 \quad \mathrm{because} \ {\tilde{\nu}}\in(0,1),
\]
which concludes the proof of (\ref{intSame-NuSameA}).
\vspace{0.2cm}

\noindent \textbf{Proof of (\ref{intNotSameNuSameA})} Let $A,B,C,D\in\R$. We define $C_{\tilde{\nu}}:=AJ_{\tilde{\nu}}+BY_{\tilde{\nu}}$ and $D_{\tilde{\nu}}:=CJ_{\tilde{\nu}}+DY_{\tilde{\nu}}$, with\\
$Y_{\tilde{\nu}}=\dfrac{J_{\tilde{\nu}}\cos({\tilde{\nu}}\pi)-J_{-{\tilde{\nu}}}}{\sin({\tilde{\nu}}\pi)}$ the second kind Bessel function. \cite[
(3), p. 254]{Luke1962} states that for any $a>0$,
\[
\begin{split}
\int^1_0xC_{{\tilde{\nu}}}(ax)D_{{\tilde{\nu}}}(ax)\mathrm{d}x=&\frac14 \left(2C_{{\tilde{\nu}}}(a)D_{{\tilde{\nu}}}(a)-C_{{\tilde{\nu}}-1}(a)D_{{\tilde{\nu}}+1}(a)-C_{{\tilde{\nu}}+1}(a)D_{{\tilde{\nu}}-1}(a)\right)\\
-\lim_{x\to0^+}&\frac14x^2 \left(2C_{{\tilde{\nu}}}(ax)D_{{\tilde{\nu}}}(ax)-C_{{\tilde{\nu}}-1}(ax)D_{{\tilde{\nu}}+1}(ax)-C_{{\tilde{\nu}}+1}(ax)D_{{\tilde{\nu}}-1}(ax)\right).
\end{split}
\]
We choose $A=1$, $B=0$, $C=\cos({\tilde{\nu}}\pi)$ and $D=-\sin({\tilde{\nu}}\pi)$. Thus, we have $C_{\tilde{\nu}}=J_{\tilde{\nu}}$ and $D_{\tilde{\nu}}=J_{-{\tilde{\nu}}}$. One needs to be careful when computing $D_{{\tilde{\nu}}+1}$ and $D_{{\tilde{\nu}}-1}$. We have
\[
D_{{\tilde{\nu}}+1}=C J_{{\tilde{\nu}}+1}+D\dfrac{J_{{\tilde{\nu}}+1}\cos(({\tilde{\nu}}+1)\pi)-J_{-{\tilde{\nu}}-1}}{\sin(({\tilde{\nu}}+1)\pi)}
=\cos({\tilde{\nu}}\pi) J_{{\tilde{\nu}}+1}-\sin({\tilde{\nu}}\pi)\dfrac{-J_{{\tilde{\nu}}+1}\cos({\tilde{\nu}}\pi)-J_{-{\tilde{\nu}}-1}}{-\sin({\tilde{\nu}}\pi)}
=-J_{-{\tilde{\nu}}-1}.
\]
Performing similar computations on $D_{{\tilde{\nu}}-1}$ gives $D_{{\tilde{\nu}}-1}=-J_{1-{\tilde{\nu}}}$. Therefore, we obtain
\[
\begin{split}
\int^1_0xJ_{{\tilde{\nu}}}(ax)J_{-{\tilde{\nu}}}(ax)\mathrm{d}x=&\frac14 \left(2J_{{\tilde{\nu}}}(a)J_{-{\tilde{\nu}}}(a)+J_{{\tilde{\nu}}-1}(a)J_{-{\tilde{\nu}}-1}(a)+J_{{\tilde{\nu}}+1}(a)J_{1-{\tilde{\nu}}}(a)\right)\\
-\lim_{x\to0^+}&\frac14x^2 \left(2J_{{\tilde{\nu}}}(ax)J_{-{\tilde{\nu}}}(ax)+J_{{\tilde{\nu}}-1}(ax)J_{-{\tilde{\nu}}-1}(ax)+J_{{\tilde{\nu}}+1}(ax)J_{1-{\tilde{\nu}}}(ax)\right).
\end{split}
\]
We first compute the limit term. By examining the series expansion of $J_{{\tilde{\nu}}}$, $J_{-{\tilde{\nu}}}$, $J_{{\tilde{\nu}}-1}$, $J_{1-{\tilde{\nu}}}$, $J_{{\tilde{\nu}}+1}$, and $J_{-{\tilde{\nu}}-1}$ given in \eqref{besselseries}, we find that
\[
\begin{split}
x^2J_{{\tilde{\nu}}}(ax)J_{-{\tilde{\nu}}}(ax)&\underset{x \to 0^+}{\sim}x^2\frac{1}{\Gamma({\tilde{\nu}}+1)}\left(\frac{ax}{2}\right)^{{\tilde{\nu}}}\frac{1}{\Gamma(1-{\tilde{\nu}})}\left(\frac{ax}{2}\right)^{-{\tilde{\nu}}}=\frac{1}{\Gamma({\tilde{\nu}}+1)}\frac{1}{\Gamma(1-{\tilde{\nu}})}x^2\underset{x\to 0^+}\longrightarrow0,\\
x^2J_{{\tilde{\nu}}+1}(ax)J_{1-{\tilde{\nu}}}(ax)&\underset{x \to 0^+}{\sim}x^2\frac{1}{\Gamma({\tilde{\nu}}+2)}\left(\frac{ax}{2}\right)^{{\tilde{\nu}}+1}\frac{1}{\Gamma(-{\tilde{\nu}}+2)}\left(\frac{ax}{2}\right)^{1-{\tilde{\nu}}}
=\frac{(ax/2)^2}{\Gamma({\tilde{\nu}}+2)\Gamma(-{\tilde{\nu}}+2)}x^2\underset{x\to 0^+}\longrightarrow0,\\
x^2J_{{\tilde{\nu}}-1}(ax)J_{-{\tilde{\nu}}-1}(ax)&\underset{x \to 0^+}{\sim}x^2\frac{1}{\Gamma({\tilde{\nu}})}\left(\frac{ax}{2}\right)^{{\tilde{\nu}}-1}\frac{1}{\Gamma(-{\tilde{\nu}})}\left(\frac{ax}{2}\right)^{-{\tilde{\nu}}-1}=\frac{1}{\Gamma({\tilde{\nu}})\Gamma(-{\tilde{\nu}})}\frac{4}{a^2}=\frac{-{\tilde{\nu}}}{\Gamma({\tilde{\nu}})\Gamma(1-{\tilde{\nu}})}\frac{4}{a^2},
\end{split}
\]
where the argument in the Euler Gamma function is always valid because ${\tilde{\nu}}\in(0,1)$. We recall the reflection formula for the Euler Gamma function, which is valid since ${\tilde{\nu}}\in(0,1)$, gives
\[\Gamma({\tilde{\nu}})\Gamma(1-{\tilde{\nu}})=\frac{\pi}{\sin({\tilde{\nu}}\pi)}.\]
Therefore, we conclude that
\[
\lim_{x\to0^+}\frac14x^2 \left(2J_{{\tilde{\nu}}}(ax)J_{-{\tilde{\nu}}}(ax)+J_{{\tilde{\nu}}-1}(ax)J_{-{\tilde{\nu}}-1}(ax)+J_{{\tilde{\nu}}+1}(ax)J_{1-{\tilde{\nu}}}(ax)\right)=-\frac{{\tilde{\nu}}\sin({\tilde{\nu}}\pi)}{\pi a^2}.
\]
Thus, to conclude the proof, we need to show that
\begin{equation}
    \frac14 \left(2J_{{\tilde{\nu}}}(a)J_{-{\tilde{\nu}}}(a)+J_{{\tilde{\nu}}-1}(a)J_{-{\tilde{\nu}}-1}(a)+J_{{\tilde{\nu}}+1}(a)J_{1-{\tilde{\nu}}}(a)\right)=\frac12\left(\left(1-\frac{{\tilde{\nu}}^2}{a^2}\right)J_{\tilde{\nu}}(a)J_{-{\tilde{\nu}}}(a)+{J_{\tilde{\nu}}}'(a){J_{-{\tilde{\nu}}}}'(a)\right).
    \label{wantedEquality}
\end{equation}
To prove this, we need to use the following identities
\[
\begin{split}
     \forall {\tilde{\nu}}\in(0,1), \forall x>0, \quad J_{{\tilde{\nu}}+1}(x)=\frac{{\tilde{\nu}}}{x}J_{\tilde{\nu}}(x)-{J_{\tilde{\nu}}}'(x), \quad J_{1-{\tilde{\nu}}}(x)=\frac{-{\tilde{\nu}}}{x}J_{-{\tilde{\nu}}}(x)-{J_{-{\tilde{\nu}}}}'(x), \\
    J_{{\tilde{\nu}}-1}(x)=\frac{{\tilde{\nu}}}{x}J_{\tilde{\nu}}(x)+{J_{\tilde{\nu}}}'(x), \quad J_{-{\tilde{\nu}}-1}(x)=\frac{-{\tilde{\nu}}}{x}J_{-{\tilde{\nu}}}(x)+{J_{-{\tilde{\nu}}}}'(x),
    \end{split}
\]
which are stated in \cite[Section 3.2, (3) and (4), p.45]{Watson1944}. Substituting $J_{{\tilde{\nu}}+1}(a)$, $J_{1-{\tilde{\nu}}}(a)$, $J_{{\tilde{\nu}}-1}(a)$, and $J_{-{\tilde{\nu}}-1}(a)$ into $\frac14 \left(2J_{{\tilde{\nu}}}(a)J_{-{\tilde{\nu}}}(a)-J_{{\tilde{\nu}}-1}(a)J_{1-{\tilde{\nu}}}(a)-J_{{\tilde{\nu}}+1}(a)J_{-{\tilde{\nu}}-1}(a)\right)$ yields directly (\ref{wantedEquality}) and thus concludes the proof of (\ref{intNotSameNuSameA}).
\vspace{0.2cm}

 \noindent \textbf{Proof of (\ref{intSame-NuSameA - alpha-beta}) and (\ref{intNotSameNuSameA - alpha-beta})}. The result in (\ref{intSame-NuSameA - alpha-beta}) is directly obtained from \cite[Section 11.2, (4)]{Luke1962}. The result in (\ref{intNotSameNuSameA - alpha-beta}) is obtained from \cite[Section 11.2, (3)]{Luke1962}, using a similar approach as in the proof of (\ref{intNotSameNuSameA}) to substitute \(J_{{\tilde{\nu}}+1}\), \(J_{1-{\tilde{\nu}}}\), \(J_{{\tilde{\nu}}-1}\), and \(J_{-{\tilde{\nu}}-1}\), which gives the desired result.
\qed

\section{Diagonalization}
\label{s:d}
\noindent From now on, we will always assume that  $\nu \in (0, 1)$. Recall that $c_{\nu}$ is defined by (\ref{1.2}). We know that $A_\nu$ is self-adjoint and non-negative by \eqref{1.6}. Therefore, we know that the spectrum of $A_\nu$ is a subset of $\R^+$, $\sigma (A_\nu) \subset [0, +\infty)$. Our goal here is to  prove that  $A_\nu$ has a Hilbert basis of eigenfunctions, and to give a rather precise description of the spectrum, together with some useful estimates on the eigenfunctions.

\subsection{Complete inventory of eigenvalues and eigenfunctions}
First of all, let us examine the kernel of $A_\nu$.
\begin{proposition}
\[
    \begin{aligned}
        \mathrm{Ker}(A_\nu)&=\mathcal{F}_s^\nu\cap D(A_\nu)\\
        &=\mathrm{Span}\left(x\mapsto|x|^{\nu+\frac12}-|x|^{-\nu+\frac12}\right).
    \end{aligned}
    \]
$\lambda_{0}=0$ is therefore the smallest eigenvalue of $A_\nu$, of multiplicity 1.
\label{PropKernel}
\end{proposition}

\noindent \textbf{\textit{Proof of Proposition \ref{PropKernel}}}. First, using (\ref{1.4}), we directly obtain $\mathcal{F}_s^\nu\cap D(A_\nu)\subset \mathrm{Ker}(A_\nu)$.

\noindent Now, let $f\in D(A_\nu)$ be such that $A_\nu f=0$. Using (\ref{1.6}), we obtain $\norm{\partial_xf_r}_{L^2(-1, 1)}=0$. 

\noindent Since $\partial_xf_r \in H^1(-1, 1) \hookrightarrow C^0(-1, 1)$ and $f_r'(0)=0$, we obtain that  $f_r=0$. Therefore, $f=f_s \in \mathcal{F}_s^\nu\cap D(A_\nu)$.

\noindent Thus, $\mathrm{Ker}(A_\nu) = \mathcal{F}_s^\nu\cap D(A_\nu)$.

\noindent Finally, it is easy to check that $\mathcal{F}_s^\nu\cap D(A_\nu)=\mathrm{Span}\left(x\mapsto |x|^{\nu+\frac12}-|x|^{-\nu+\frac12}\right)$, using the transmission conditions in the definition of $A_\nu$ together with the Dirichlet boundary conditions at $\pm 1$.
\qed

\vspace{0.3cm}

\noindent Let $E>0$. Now, using Bessel functions, we want to find $f\in D(A_\nu)$ such that $A_\nu f=Ef$. Let us first investigate the differential equation $A_\nu f=Ef$, without taking into account the transmission and boundary conditions.

\begin{proposition} Let $E>0$. Let $f\in C^\infty((-1,0)\cup(0,1))$. Then, the function $f$ satisfies
\begin{equation}
        \forall x\in(-1, 0)\cup(0, 1), \quad \left(-\frac{d^2}{\mathrm{d}x^2} + \frac{\nu^2-\frac14}{x^2}\right)f(x) = E f(x),
        \label{equivalentEigen}
\end{equation}
    if and only if $f$ is of the form
    \begin{equation}
        f(x) = \begin{cases}
            a_\nu^-\sqrt{-x}J_{\nu}(-\sqrt{E}x)+a_{-\nu}^-\sqrt{-x}J_{-\nu}(-\sqrt{E}x), &\quad \forall x \in (-1, 0),
                       \\ a_\nu^+\sqrt{x}J_{\nu}(\sqrt{E}x)+a_{-\nu}^+\sqrt{x}J_{-\nu}(\sqrt{E}x), &\quad \forall x \in (0, 1),
        \end{cases}
        \label{formieigenE}
    \end{equation}
with $a_\nu^-, a_{-\nu}^-,a_\nu^+,a_{-\nu}^+\in\R$.

\vspace{0.2cm}

\noindent Therefore, a function $f$ is in the eigenspace associated with the eigenvalue $E$ if and only if $f\in D(A_\nu)$ and $f$ is of the form (\ref{formieigenE}).
    
    \label{propFormEigen}
    \end{proposition}

\noindent \textbf{\textit{Proof of Proposition \ref{propFormEigen}}}.

\noindent \textbf{On (0, 1)}. We consider the function $\psi_{\nu,E}^+(x) = \sqrt{x} J_\nu(\sqrt{E}x)$, defined on $(0, 1)$. Taking the second derivative, we obtain
\[
\psi_{\nu,E}^{+''}(x) = -\frac{1}{4x^{3/2}} J_\nu(\sqrt{E}x) + \frac{1}{\sqrt{x}} \sqrt{E} J'_\nu(\sqrt{E}x) + E \sqrt{x} J''_\nu(\sqrt{E}x).
\]

\noindent Next, we use the Bessel equation (\ref{besseleq}) satisfied by $J_\nu$ at $\sqrt{E}x>0$, which gives
\[
Ex^2 J''_\nu(\sqrt{E}x) + x\sqrt{E} J'_\nu(\sqrt{E}x) + E x^2 J_\nu(\sqrt{E}x) = \nu^2 J_\nu(\sqrt{E}x).
\]
Therefore, we have, for all $x\in(0, 1)$,
\[
\begin{aligned}
\left(-\frac{d^2}{\mathrm{d}x^2} + \frac{\nu^2-\frac14}{x^2}\right)\psi_{\nu,E}^+(x) &= \frac{1}{4x^{3/2}} J_\nu(\sqrt{E}x) - \frac{1}{\sqrt{x}} \sqrt{E} J'_\nu(\sqrt{E}x) - E \sqrt{x} J''_\nu(\sqrt{E}x)\\
& - \frac{1}{4x^2}\sqrt{x}J_\nu(\sqrt{E}x)+ \frac{\sqrt{x}}{x^2} (Ex^2 J''_\nu(\sqrt{E}x) + x\sqrt{E} J'_\nu(\sqrt{E}x) + E x^2 J_\nu(\sqrt{E}x))\\
& = E \psi_{\nu,E}^+(x).
\end{aligned}
\]

\noindent Since $\nu\in\mathbb{R}\setminus\mathbb{Z}$, we have that $(J_\nu, J_{-\nu})$ is a basis of solutions of (\ref{besseleq}) by Lemma \ref{lemmaWronskien}. Therefore, by posing $\psi_{-\nu,E}^+(x) = \sqrt{x}J_{-\nu}(\sqrt{E}x)$ on $(0, 1)$, we also have that for all $x\in(0, 1)$,

\[
\left(-\frac{d^2}{\mathrm{d}x^2} + \frac{\nu^2-\frac14}{x^2}\right)\psi_{-\nu,E}^+(x) = E \psi_{-\nu,E}^+(x).
\]
\noindent \textbf{On (-1, 0)}. We follow the same computations as on $(0,1)$ and we observe that they are the same.
\vspace{0.2cm}

\noindent We proved that any function of the form (\ref{formieigenE}) satisfies (\ref{equivalentEigen}). We have found all the solutions to (\ref{equivalentEigen}) since the total solution space is of dimension 4, consisting of two independent solutions on $(0,1)$ and two independent solutions on $(-1,0)$. This matches the expected structure, as the equation is second-order and has a singularity at $x=0$, which separates the two intervals.
\qed

\vspace{0.3cm}

\noindent Now, we want to find conditions on $a_{-\nu}^-, \ a_{\nu}^-, \ a_{-\nu}^+$, and $a_{\nu}^+$ to ensure that the function $f$ defined in (\ref{formieigenE}) is in $D(A_\nu)$.

\begin{proposition} Let $\nu\in(0,1).$ Any function of the form (\ref{formieigenE}) is in $\tilde{H}_0^2(-1, 1)\oplus\mathcal{F}_s^\nu$.
    \label{propEigenInSum}
\end{proposition}

\noindent \textbf{\textit{Proof of Proposition \ref{propEigenInSum}}}. By examining the power series expansion (\ref{besselseries}) of the Bessel functions \(J_\nu\) and isolating the term for \(k=0\), we have
\[
\begin{split}
\sqrt{x}J_\nu(\sqrt{E}x)&=\frac{1}{\Gamma(\nu+1)}\left(\frac{\sqrt{E}}{2}\right)^\nu |x|^{\nu+\frac12}+
\displaystyle\sum_{k=1}^{\infty} \frac{(-1)^k}{k! \, \Gamma(k + \nu + 1)} \left(\frac{\sqrt{E}}{2}\right)^{2k + \nu}x^{2k + \nu + \frac12},\quad\forall x \in (0, 1),\\
\sqrt{x}J_{-\nu}(\sqrt{E}x)&=\frac{1}{\Gamma(1-\nu)}\left(\frac{\sqrt{E}}{2}\right)^{-\nu} |x|^{-\nu+\frac12}+
\displaystyle\sum_{k=1}^{\infty} \frac{(-1)^k}{k! \, \Gamma(k - \nu + 1)} \left(\frac{\sqrt{E}}{2}\right)^{2k - \nu}x^{2k - \nu + \frac12},\quad\forall x \in (0, 1).
\end{split}
\]

\noindent We introduce 
\[
f_r(x):= 
\begin{cases}
a_{\nu}^-\displaystyle\sum_{k=1}^{\infty} \frac{(-1)^k}{k! \, \Gamma(k + \nu + 1)} \left(\frac{\sqrt{E}}{2}\right)^{2k + \nu}(-x)^{2k + \nu + \frac12}
\\+a_{-\nu}^-
\displaystyle\sum_{k=1}^{\infty} \frac{(-1)^k}{k! \, \Gamma(k - \nu + 1)} \left(\frac{\sqrt{E}}{2}\right)^{2k - \nu}(-x)^{2k - \nu + \frac12}
, & \forall x \in (-1, 0),\\

    a_{\nu}^+\displaystyle\sum_{k=1}^{\infty} \frac{(-1)^k}{k! \, \Gamma(k + \nu + 1)} \left(\frac{\sqrt{E}}{2}\right)^{2k + \nu}x^{2k + \nu + \frac12}
    \\+ a_{-\nu}^+
\displaystyle\sum_{k=1}^{\infty} \frac{(-1)^k}{k! \, \Gamma(k - \nu + 1)} \left(\frac{\sqrt{E}}{2}\right)^{2k - \nu}x^{2k - \nu + \frac12}, &  \forall x \in (0, 1),
\end{cases}
\]

\begin{equation}
f_s(x):= 
\begin{cases}
     a_{\nu}^-\frac{1}{\Gamma(\nu+1)}\left(\frac{\sqrt{E}}{2}\right)^\nu |x|^{\nu+\frac12}
+a_{-\nu}^-
\frac{1}{\Gamma(1-\nu)}\left(\frac{\sqrt{E}}{2}\right)^{-\nu} |x|^{-\nu+\frac12}
, &\quad \forall x \in (-1, 0),\\
    a_{\nu}^+\frac{1}{\Gamma(\nu+1)}\left(\frac{\sqrt{E}}{2}\right)^\nu |x|^{\nu+\frac12}
    + a_{-\nu}^+
\frac{1}{\Gamma(1-\nu)}\left(\frac{\sqrt{E}}{2}\right)^{-\nu} |x|^{-\nu+\frac12}, &\quad \forall x \in (0, 1).
\end{cases}
\label{f_s}
\end{equation}

\noindent Thus, we have $f=f_s+f_r$. Since $\nu\in(0, 1)$, it is clear that $f_s\in L^2(-1, 1)$. By examining (\ref{f_s}), we see that $f_s$ has the right form with
\begin{equation}
    \begin{aligned}
        &c_1^+=a_{\nu}^+\frac{1}{\Gamma(\nu+1)}\left(\frac{\sqrt{E}}{2}\right)^\nu, \quad c_2^+=a_{-\nu}^+\frac{1}{\Gamma(1-\nu)}\left(\frac{\sqrt{E}}{2}\right)^{-\nu}, \\
        &c_1^-=a_{\nu}^-\frac{1}{\Gamma(\nu+1)}\left(\frac{\sqrt{E}}{2}\right)^\nu, \quad c_2^-=a_{-\nu}^-
\frac{1}{\Gamma(1-\nu)}\left(\frac{\sqrt{E}}{2}\right)^{-\nu} .
    \end{aligned}
    \label{c}
\end{equation}
Therefore, $f_s\in\mathcal{F}_s^\nu$. Finally, we need to verify that $f_r \in \tilde{H}_{0}^{2}(-1,1)$. We introduce
\[
g_r(x)=
\begin{cases}
a_{\nu}^-\sum_{k=1}^{\infty} \frac{(-1)^k}{k! \, \Gamma(k + \nu + 1)} \left(\frac{\sqrt{E}}{2}\right)^{2k + \nu}(-x)^{2k + \nu + \frac12}
\\+a_{-\nu}^-
\sum_{k=2}^{\infty} \frac{(-1)^k}{k! \, \Gamma(k - \nu + 1)} \left(\frac{\sqrt{E}}{2}\right)^{2k - \nu}(-x)^{2k - \nu + \frac12}
, &\ \forall x \in (-1, 0),\\
    a_{\nu}^+\sum_{k=1}^{\infty} \frac{(-1)^k}{k! \, \Gamma(k + \nu + 1)} \left(\frac{\sqrt{E}}{2}\right)^{2k + \nu}x^{2k + \nu + \frac12}
    \\+ a_{-\nu}^+
\sum_{k=2}^{\infty} \frac{(-1)^k}{k! \, \Gamma(k - \nu + 1)} \left(\frac{\sqrt{E}}{2}\right)^{2k - \nu}x^{2k - \nu + \frac12}, &\ \forall x \in (0, 1), 
\end{cases}
\]
\[
h_r(x)=
\begin{cases}
    a_{-\nu}^-\frac{-1}{\Gamma(1 - \nu + 1)} \left(\frac{\sqrt{E}}{2}\right)^{2 - \nu}(-x)^{2 - \nu + \frac12}
, &\ \forall x \in (-1, 0),\\
    a_{-\nu}^+\frac{-1}{\Gamma(1 - \nu + 1)} \left(\frac{\sqrt{E}}{2}\right)^{2 - \nu}x^{2 - \nu + \frac12}, &\ \forall x \in (0, 1).
\end{cases}
\]

\noindent Thus, $f_r=g_r+h_r$. Since $\nu\in(0, 1)$, we have $g_r\in C^2([-1, 1]) \subset H^2(-1, 1)$ and $g_r(0)=g_r'(0)=0$, therefore, $g_r\in \tilde{H}_0^2(-1, 1)$. We also have $h_r\in C^1([-1, 1])$ and $h_r(0)=h_r'(0)=0$. We only need to check that $h_r''\in L^2(-1, 1)$. Since for $\nu > \frac12, \ h_r'$ is only differentiable on $(-1, 0)\cup(0, 1)$ and not at $0$, we need to compute $h_r''$ using the derivative in the distributional sense. Let $\phi \in C_c^{\infty}(-1, 1)$,
\[
\begin{aligned}
    \langle h_r'', \phi \rangle &= - \int^1_{-1} h_r' \phi'= - \int^0_{-1} h_r' \phi' -\int^1_{0} h_r' \phi'= \int^0_{-1} h_r'' \phi - [h_r'\phi]_{-1}^0 + \int^1_{0} h_r'' \phi - [h_r'\phi]_{0}^1\\
    &= \int^0_{-1} h_r'' \phi + \int^1_{0} h_r'' \phi \quad \mathrm{since} \ h_r'(0)=0.
\end{aligned}
\]
Therefore, $h_r''$ is obtained by examining the derivative of $h_r'$ on $(-1, 0)$ and $(0, 1)$. We find
\[
h_r''(x)=
\begin{cases}
    a_{-\nu}^-\frac{-1}{\Gamma(1 - \nu + 1)} \left(\frac{\sqrt{E}}{2}\right)^{2 - \nu}(2-\nu+\frac12)(1-\nu+\frac12)(-x)^{ - \nu + \frac12}
, \quad \forall x \in (-1, 0),\\
    a_{-\nu}^+\frac{-1}{\Gamma(1 - \nu + 1)} \left(\frac{\sqrt{E}}{2}\right)^{2 - \nu}(2-\nu+\frac12)(1-\nu+\frac12)x^{- \nu + \frac12}, \quad \forall x \in (0, 1), \\
\end{cases}
\]
which is in $L^2(-1, 1)$ since $\nu\in(0, 1)$. Thus, $h_r\in\tilde{H}_0^2(-1, 1)$ and $g_r\in\tilde{H}_0^2(-1, 1)$. Therefore, $f_r=g_r+h_r\in\tilde{H}_0^2(-1, 1)$. We have shown that $f=f_s+f_r\in\tilde{H}_0^2(-1, 1)\oplus\mathcal{F}_s^\nu$.
\qed

\begin{proposition} Let $\nu\in(0, 1)$, let $E>0$. The value $E$ is a positive eigenvalue of $A_\nu$ if and only if 
\begin{equation}
J_\nu(\sqrt{E}) \left(\frac{\sqrt{E}}{2}\right)^{-2\nu}\frac{\Gamma(\nu+1)}{\Gamma(1-\nu)}=J_{-\nu}(\sqrt{E}),
    \label{maincond+}
\end{equation}
or
\begin{equation}
J_\nu(\sqrt{E}) \left(\frac{\sqrt{E}}{2}\right)^{-2\nu}\frac{\Gamma(\nu+1)}{\Gamma(1-\nu)}\frac{1-2\nu}{1+2\nu}=J_{-\nu}(\sqrt{E}).
    \label{maincond-}
\end{equation}

\noindent Furthermore, $J_\nu(\sqrt{E})\neq0$ and
\begin{itemize}
    \item if (\ref{maincond+}) is satisfied, the eigenspace associated with the eigenvalue $E$ is
    \begin{equation}
            \mathrm{Span}\left(x\mapsto 
    \begin{cases}
        -\dfrac{\Gamma(\nu+1)}{\Gamma(1-\nu)}\left(\dfrac{\sqrt{E}}{2}\right)^{-2\nu}\sqrt{-x}J_\nu(-\sqrt{E}x)+\sqrt{-x}J_{-\nu}(-\sqrt{E}x)
         \ &\text{on} \ (-1,0)
        \\
        -\dfrac{\Gamma(\nu+1)}{\Gamma(1-\nu)}\left(\dfrac{\sqrt{E}}{2}\right)^{-2\nu}\sqrt{x}J_\nu(\sqrt{E}x)+\sqrt{x}J_{-\nu}(\sqrt{E}x) \ & \text{on} \ (0,1)
    \end{cases}
    \right),
    \label{eigenSpace+}
        \end{equation}
    \item if (\ref{maincond-}) is satisfied, the eigenspace associated with the eigenvalue $E$ is
    \begin{equation}
            \mathrm{Span}\left(x\mapsto
    \begin{cases}
        \dfrac{1-2\nu}{1+2\nu}\dfrac{\Gamma(\nu+1)}{\Gamma(1-\nu)}\left(\dfrac{\sqrt{E}}{2}\right)^{-2\nu}\sqrt{-x}J_\nu(-\sqrt{E}x)-\sqrt{-x}J_{-\nu}(-\sqrt{E}x)
         \ &\text{on} \ (-1,0)
        \\
        -\dfrac{1-2\nu}{1+2\nu}\dfrac{\Gamma(\nu+1)}{\Gamma(1-\nu)}\left(\dfrac{\sqrt{E}}{2}\right)^{-2\nu}\sqrt{x}J_\nu(\sqrt{E}x)+\sqrt{x}J_{-\nu}(\sqrt{E}x) \ &\text{on} \ (0,1)
    \end{cases}\right),
    \label{eigenSpace-}
    \end{equation}
\end{itemize}
which are both of dimension 1.
\label{PropCondSimpler}
\end{proposition}

\noindent \textbf{\textit{Proof of Proposition \ref{PropCondSimpler}}}. We want to find conditions that allow candidate eigenfunctions written as (\ref{formieigenE}) to be in $D(A_\nu)$. Referring to the definition of $D(A_\nu)$ in (\ref{1.5}), Proposition \ref{propEigenInSum}, and (\ref{c}), the function $f$ defined in (\ref{formieigenE}) is in $D(A_\nu)$ if and only if 

\begin{equation}
\begin{aligned}
    (a_{\nu}^++a_{\nu}^-)\frac{1}{\Gamma(\nu+1)}\left(\frac{\sqrt{E}}{2}\right)^\nu 
    + (a_{-\nu}^++a_{-\nu}^-)
\frac{1}{\Gamma(1-\nu)}\left(\frac{\sqrt{E}}{2}\right)^{-\nu}=0,
\end{aligned}
\label{cond1}
\end{equation}

\begin{equation}
(a_{\nu}^+-a_{\nu}^-)(\nu+\frac{1}{2})\frac{1}{\Gamma(\nu+1)}\left(\frac{\sqrt{E}}{2}\right)^\nu +(a_{-\nu}^+-a_{-\nu}^-)(-\nu+\frac{1}{2})\frac{1}{\Gamma(1-\nu)}\left(\frac{\sqrt{E}}{2}\right)^{-\nu} =0,
\label{cond2}
\end{equation}

\begin{equation}
     f(1)=a_\nu^+J_\nu(\sqrt{E})+a_{-\nu}^+J_{-\nu}(\sqrt{E})=0,
    \label{cond3}
\end{equation}

\begin{equation}
     f(-1)=a_\nu^-J_\nu(\sqrt{E})+a_{-\nu}^-J_{-\nu}(\sqrt{E})=0.
    \label{cond4}
\end{equation}
We introduce the useful quantity
\begin{equation}
        D_\nu:=\frac{1}{1+2\nu}\frac{\Gamma(\nu+1)}{\Gamma(1-\nu)}\left(\frac{\sqrt{E}}{2}\right)^{-2\nu}.
\label{Dnu}
\end{equation}

\noindent $\frac12\left[\Gamma(\nu+1)\left(\frac{\sqrt{E}}{2}\right)^{-\nu}\times(\ref{cond1})+\frac{\Gamma(\nu+1)}{\nu+\frac12}\left(\frac{\sqrt{E}}{2}\right)^{-\nu}\times(\ref{cond2})\right]$ gives

\begin{equation}
    a^+_\nu=-D_\nu(a^+_{-\nu}+2\nu a^-_{-\nu}),
    \label{newcond1}
\end{equation}
\noindent and $\frac12\left[\Gamma(\nu+1)\left(\frac{\sqrt{E}}{2}\right)^{-\nu}\times(\ref{cond1})-\frac{\Gamma(\nu+1)}{\nu+\frac12}\left(\frac{\sqrt{E}}{2}\right)^{-\nu}\times(\ref{cond2})\right]$ gives

\begin{equation}
    a^-_\nu=-D_\nu(a^-_{-\nu}+2\nu a^+_{-\nu}).
    \label{newcond2}
\end{equation}
\noindent We reason by contradiction and assume that $J_\nu(\sqrt{E})=0$. Then, (\ref{cond3}) and (\ref{cond4}) become $a^+_{-\nu}J_{-\nu}(\sqrt{E})=0$ and $a^-_{-\nu}J_{-\nu}(\sqrt{E})=0$. Assume once more by contradiction that $J_{-\nu}(\sqrt{E})\neq0$. Then, we must have $a^+_{-\nu}=a^-_{-\nu}=0$. Using (\ref{newcond1}) and (\ref{newcond2}), we also must have $a^+_{\nu}=a^-_{\nu}=0$. Thus, $f=0$. This is not possible since we are looking for eigenfunctions. Hence, we must have $J_{-\nu}(\sqrt{E})=0$, which is not possible since Lemma \ref{lemmaInterlaced} shows that the positive zeros of $J_\nu$ and $J_{-\nu}$ are strictly interlaced.

\noindent Therefore, $J_\nu(\sqrt{E})\neq0$. We substitute $a^-_\nu$ in (\ref{cond4}) using (\ref{newcond2}). This gives
\[
a^+_{-\nu}J_\nu(\sqrt{E})2\nu D_\nu=a^-_{-\nu}\left(J_{-\nu}(\sqrt{E})-J_\nu(\sqrt{E})D_\nu\right).
\]
\noindent Since $J_\nu(\sqrt{E})\neq0$, we have

\begin{equation}
a^+_{-\nu}=a^-_{-\nu}\frac{J_{-\nu}(\sqrt{E})-J_\nu(\sqrt{E})D_\nu}{J_\nu(\sqrt{E})2\nu D_\nu}.
\label{newcond4}
\end{equation}

\noindent Substituting $a^+_\nu$ in (\ref{cond3}) and using (\ref{newcond1}) gives
\[
a^-_{-\nu}J_\nu(\sqrt{E})2\nu D_\nu=a^+_{-\nu}\left(J_{-\nu}(\sqrt{E})-J_\nu(\sqrt{E})D_\nu\right).
\]
We now substitute $a^+_{-\nu}$ in the previous equality, using (\ref{newcond4}). This gives
\[
a^-_{-\nu}J_\nu(\sqrt{E})2\nu D_\nu=a^-_{-\nu}\frac{J_{-\nu}(\sqrt{E})-J_\nu(\sqrt{E})D_\nu}{J_\nu(\sqrt{E})2\nu D_\nu}\left(J_{-\nu}(\sqrt{E})-J_\nu(\sqrt{E})D_\nu\right),
\]
which can be rewritten as
\[
a^-_{-\nu}\left(J_\nu(\sqrt{E})2\nu D_\nu\right)^2=a^-_{-\nu}\left(J_{-\nu}(\sqrt{E})-J_\nu(\sqrt{E})D_\nu\right)^2.
\]
We cannot have $a^-_{-\nu}=0$ because this would lead to $a^+_{-\nu}=0$ using (\ref{newcond4}), and then $a^-_{\nu}=a^+_{\nu}=0$ using (\ref{newcond1}) and (\ref{newcond2}). Therefore, we must have
\begin{equation}
    \left(J_\nu(\sqrt{E})2\nu D_\nu\right)^2=\left(J_{-\nu}(\sqrt{E})-J_\nu(\sqrt{E})D_\nu\right)^2.
    \label{maincondproof}
\end{equation}
\noindent To summarize, we proved that if $E$ is an eigenvalue of $A_\nu$ for a function $f$ of the form (\ref{formieigenE}), then we must have (\ref{maincondproof}), and we must have $a^+_{-\nu}$ satisfying (\ref{newcond4}), while $a^+_\nu$ and $a^-_\nu$ must satisfy (\ref{newcond1}) and (\ref{newcond2}) respectively, with $a^-_{-\nu}$ chosen freely in $\R$.

\noindent On the other hand, if (\ref{maincondproof}) is satisfied, we can choose any function from the set
\begin{equation}
\left.
\begin{cases}
f \mathrm{\ of \ the\ form\ (\ref{formieigenE})\ such\ that:} \quad a^-_{-\nu} \in \mathbb{R}, \quad a^+_{-\nu}=a^-_{-\nu}\dfrac{J_{-\nu}(\sqrt{E})-J_\nu(\sqrt{E})D_\nu}{J_\nu(\sqrt{E})2\nu D_\nu},\\
    a^+_\nu = -D_\nu(a^+_{-\nu} + 2\nu a^-_{-\nu}), \quad a^-_\nu = -D_\nu(a^-_{-\nu} + 2\nu a^+_{-\nu})
\end{cases}\right\}
\label{proofSpaceEigen}
\end{equation}
and we indeed have $A_\nu f=Ef$ and $f\in D(A_\nu)$. Therefore, $E$ is an eigenvalue and the eigenspace associated with the eigenvalue $E$ is given by (\ref{proofSpaceEigen}).

\noindent Now, we will distinguish cases on the condition (\ref{maincondproof}). If (\ref{maincondproof}) holds, then we have either 
\begin{equation}
J_\nu(\sqrt{E})2\nu D_\nu=J_{-\nu}(\sqrt{E})-J_\nu(\sqrt{E})D_\nu
\label{cond1Eigen}
\end{equation}
or
\begin{equation}
    -J_\nu(\sqrt{E})2\nu D_\nu=J_{-\nu}(\sqrt{E})-J_\nu(\sqrt{E})D_\nu.
    \label{cond2Eigen}
\end{equation}
We notice that condition (\ref{cond1Eigen}) is equivalent to condition (\ref{maincond+}), and condition (\ref{cond2Eigen}) is equivalent to condition (\ref{maincond-}). Therefore, what remains to be proved is that the eigenspace (\ref{proofSpaceEigen}) takes the form stated in the proposition, whether (\ref{cond1Eigen}) or (\ref{cond2Eigen}) holds.

\noindent Suppose that (\ref{cond1Eigen}) holds. It is straightforward that the eigenspace (\ref{proofSpaceEigen}) can be rewritten as
\[
\left.
\begin{cases}
f \mathrm{\ of \ the\ form\ (\ref{formieigenE})\ such\ that:} \quad a^-_{-\nu} \in \mathbb{R}, \quad a^+_{-\nu}=a^-_{-\nu},\\
    a^+_\nu = -D_\nu(1+2\nu)a^-_{-\nu}, \quad a^-_\nu = -D_\nu(1+2\nu)a^-_{-\nu}
\end{cases}\right\}.
\]
Looking at the expression of $D_\nu$ in (\ref{Dnu}) and at functions of the form (\ref{formieigenE}), it is clear that this space is precisely (\ref{eigenSpace+}).

\noindent Now, suppose that (\ref{cond2Eigen}) holds. It is straightforward that the eigenspace (\ref{proofSpaceEigen}) can be rewritten as
\[
\left.
\begin{cases}
f \mathrm{\ of \ the\ form\ (\ref{formieigenE})\ such\ that:} \quad a^-_{-\nu} \in \mathbb{R}, \quad a^+_{-\nu}=-a^-_{-\nu},\\
    a^+_\nu = D_\nu(1-2\nu)a^-_{-\nu}, \quad a^-_\nu = -D_\nu(1-2\nu)a^-_{-\nu}
\end{cases}\right\}.
\]
By looking at the expression of $D_\nu$ in (\ref{Dnu}) and at the functions of the form (\ref{formieigenE}), it is clear that this space is precisely (\ref{eigenSpace-}). This concludes the proof of Proposition \ref{PropCondSimpler}.
\qed

\begin{remark} It is crucial to emphasize that we have found all the eigenvalues and eigenfunctions of our operator $A_\nu$, for any $\nu\in(0,1)$. Indeed, we know that $A_\nu$ is self-adjoint and non-negative by \eqref{1.6}. Therefore, $\sigma_{point}(A_\nu)\subset\mathbb{R}_+$. Proposition \ref{PropKernel} shows that 0 is an eigenvalue of multiplicity 1 and provides the kernel. Proposition \ref{PropCondSimpler} indicates that the other eigenvalues are given by the implicit equations (\ref{maincond+}) and (\ref{maincond-}), and provides the associated eigenspaces, which are also of dimension one.
    \label{rqAllEigenFound}
    \end{remark}
\subsection{Distribution structure of the eigenvalues}

Now that we have found the conditions (\ref{maincond+}) and (\ref{maincond-}) for our eigenvalues, we would like to determine which values of $E$ satisfy one of these conditions, and how these values are distributed along the positive real axis. This is the subject of this subsection, which concludes with a proposition characterizing the structure of the distribution of our eigenvalues

\begin{proposition} Let $\nu\in(0, 1)$, let $n\in\N\setminus\{0\}$. The function $f: x\mapsto \frac{J_\nu(x)}{J_{-\nu}(x)}x^{-2\nu}$ is a well-defined, increasing, and continuous bijection from $(j_{-\nu,n},j_{-\nu, n+1})$ to $\R$. Furthermore, $f$ is also a well-defined, increasing, and continuous bijection from $(0, j_{-\nu, 1})$ to $(\frac{\Gamma(1-\nu)}{\Gamma(\nu+1)}2^{-2\nu}, +\infty)$.\\
\label{coroMonotone}
\end{proposition}

\noindent \textbf{\textit{Proof of Proposition \ref{coroMonotone}.}} Let $n\in\N\setminus\{0\}$.
First, by Lemma \ref{lemmaInterlaced}, we know that $J_{-\nu}$ does not vanish on $(j_{-\nu,n},j_{-\nu, n+1})$ or on $(0, j_{-\nu, 1})$. Therefore, $f$ is well-defined on $(j_{-\nu,n},j_{-\nu, n+1})$ or $(0, j_{-\nu, 1})$. The function $f$ is also differentiable on these intervals, so to show that $f$ is increasing, it suffices to examine its derivative. Let $ x\in(j_{-\nu,n},j_{-\nu, n+1})$ or $x\in(0, j_{-\nu, 1})$. Then,
\[
\begin{aligned}
   f'(x)&=\frac{\left({J_\nu}'(x)x^{-2\nu}-2\nu J_\nu(x)x^{-2\nu-1}\right)J_{-\nu}(x)-J_\nu(x)x^{-2\nu}{J_{-\nu}}'(x)}{{J_{-\nu}}^2(x)}\\
   &=\frac{-x^{-2\nu}W\left(J_\nu,J_{-\nu}\right)(x)-2\nu x^{-2\nu-1} J_\nu(x) J_{-\nu}(x)}{{J_{-\nu}}^2(x)} \mathrm{ \ with\ } W\left(J_\nu,J_{-\nu}\right) \mathrm{ \ defined\ in\ Lemma\ \ref{lemmaWronskien},} \\
   &=\frac{2\nu x^{-2\nu-1}\left(\frac{\sin(\nu\pi)}{\nu\pi} - J_\nu(x)J_{-\nu}(x)\right)}{{J_{-\nu}}^2(x)} > 0 \quad \text{by Lemma \ref{lemmeInequality}.}
\end{aligned}
\]
Therefore, $f$ is increasing on $(j_{-\nu,n},j_{-\nu, n+1})$, and on $(0, j_{-\nu, 1})$. Furthermore, $J_{-\nu}(j_{-\nu,n}) = J_{-\nu}(j_{-\nu, n+1})=0$, and $J_{\nu}(j_{-\nu,n}) \neq 0$, $J_{\nu}(j_{-\nu, n+1})\neq0$ because the positive real zeros of $J_{\nu}$ and $J_{-\nu}$ are strictly interlaced as shown in Lemma \ref{lemmaInterlaced}. Therefore, we have
\[
\lim_{x \to j_{-\nu,n}^+} f(x) = -\infty, \quad
\lim_{x \to j_{-\nu,n+1}^-} f(x) = +\infty,\ \mathrm{and} \ 
\lim_{x \to j_{-\nu,1}^-} f(x) = +\infty.
\]
Finally, from the asymptotic expansion \eqref{besselseries}, we have  $$J_{\nu}(x) \sim_{0^+} \frac{1}{\Gamma(\nu+1)} \left(\frac{x}{2}\right)^{\nu} \mbox{ and }J_{-\nu}(x) \sim_{0^+} \frac{1}{\Gamma(1-\nu)} \left(\frac{x}{2}\right)^{-\nu}.$$
We deduce that  $\lim_{x \to 0^+}f(x)=\frac{\Gamma(1-\nu)}{\Gamma(\nu+1)}2^{-2\nu}$, which concludes the proof.
\qed

\noindent Using the above property, we can state a result on the distribution of the eigenvalues of $A_\nu$.

\begin{proposition}[Distribution of the eigenvalues]\label{thDistribution}
Fix $\nu \in (0,1)$.

\medskip
\noindent\textbf{Even eigenvalues.}
The even eigenvalues $(\lambda_{2n})_{n\geq 0}$ are characterized as follows.
\begin{itemize}
    \item $\lambda_{0} = 0$.
    \item $\nexists\, E\in(0, j^2_{-\nu, 1})$ such that $E$ satisfies \eqref{maincond+}.
    \item $\forall n\in\N\setminus\{0\},\ \nexists\, E\in[j^2_{-\nu,n}, j^2_{\nu,n}]$ such that $E$ satisfies \eqref{maincond+}.
    \item $\forall n\in\N\setminus\{0\},\ \exists!\, E\in(j^2_{\nu,n}, j^2_{-\nu,n+1})$ such that $E$ satisfies \eqref{maincond+}. We denote it by $\lambda_{2n}$.
\end{itemize}

\medskip
\noindent\textbf{Odd eigenvalues.}
The odd eigenvalues $(\lambda_{2n+1})_{n\geq 0}$ depend on the value of $\nu$ as follows.

\smallskip
\noindent\emph{Case $\nu\in(0,\frac12)$:}
\begin{itemize}
    \item $\exists!\, E\in (0, j^2_{-\nu, 1})$ such that $E$ satisfies \eqref{maincond-}. We denote it by $\lambda_{1}$.
    \item $\forall n\in\N\setminus\{0\},\ \nexists\, E\in[j^2_{-\nu,n}, j^2_{\nu,n}]$ such that $E$ satisfies \eqref{maincond-}.
    \item $\forall n\in \N\setminus\{0\},\ \exists!\, E\in(j^2_{\nu,n}, j^2_{-\nu,n+1})$ such that $E$ satisfies \eqref{maincond-}. We denote it by $\lambda_{2n+1}$.
\end{itemize}

\smallskip
\noindent\emph{Case $\nu\in(\frac12,1)$:}
\begin{itemize}
    \item $\nexists\, E\in (0, j^2_{-\nu, 1}]$ such that $E$ satisfies \eqref{maincond-}.
    \item $\forall n\in\N\setminus\{0\},\ \exists!\, E\in(j^2_{-\nu,n}, j^2_{\nu,n})$ such that $E$ satisfies \eqref{maincond-}. We denote it by $\lambda_{2n-1}$.
    \item $\forall n\in \N\setminus\{0\},\ \nexists\, E\in[j^2_{\nu,n}, j^2_{-\nu,n+1}]$ such that $E$ satisfies \eqref{maincond-}.
\end{itemize}

\smallskip
\noindent\emph{Case $\nu=\frac12$:}
\begin{itemize}
    \item $\forall n\in\N\setminus\{0\},\ E=j^2_{-\nu,n}$ satisfies \eqref{maincond-}. We denote it by $\lambda_{2n-1}$.
    \item $\nexists\, E \in \R^+ \setminus\Bigl(\{0\} \cup_{n \geq 1} \{ j^2_{-\nu,n} \}\Bigr)$ such that $E$ satisfies \eqref{maincond-}.
\end{itemize}

\medskip
In particular, for every $\nu\in(0,1)$ we obtain the increasing sequence
\begin{equation}
    0=\lambda_{0}<\lambda_{1}<\lambda_{2}<\dots<\lambda_{2n}<\lambda_{2n+1}<\dots
    \label{order}
\end{equation}
\end{proposition}

\begin{figure}[H]
  \centering

  \begin{subfigure}[b]{0.495\textwidth}
    \centering
    \includegraphics[width=\textwidth]{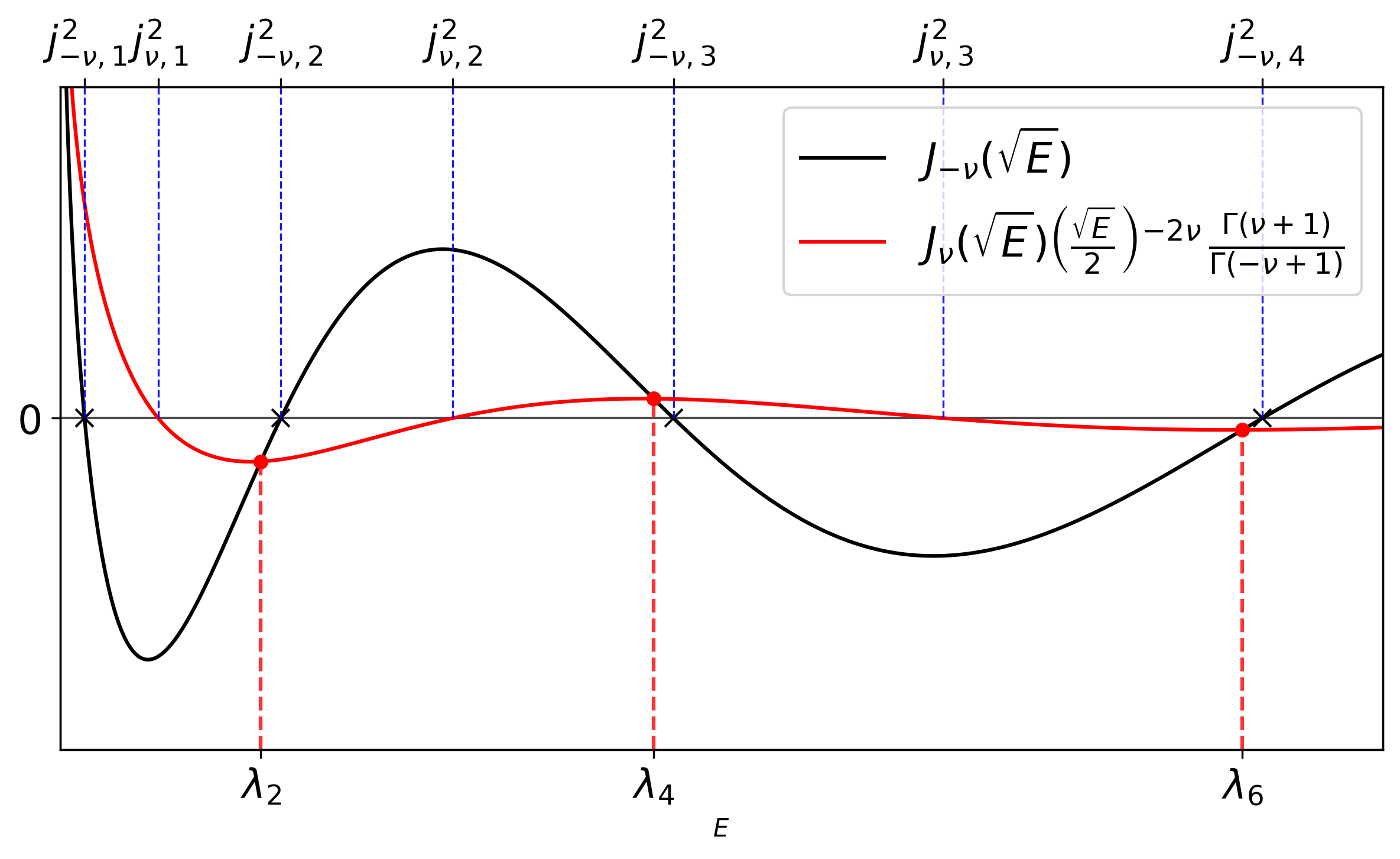}
    \caption{Even eigenvalues $\lambda_{2n}$ for $\nu=1/2$}
    \label{fig:subfig_a}
  \end{subfigure}
  \hfill
  \begin{subfigure}[b]{0.495\textwidth}
    \centering
    \includegraphics[width=\textwidth]{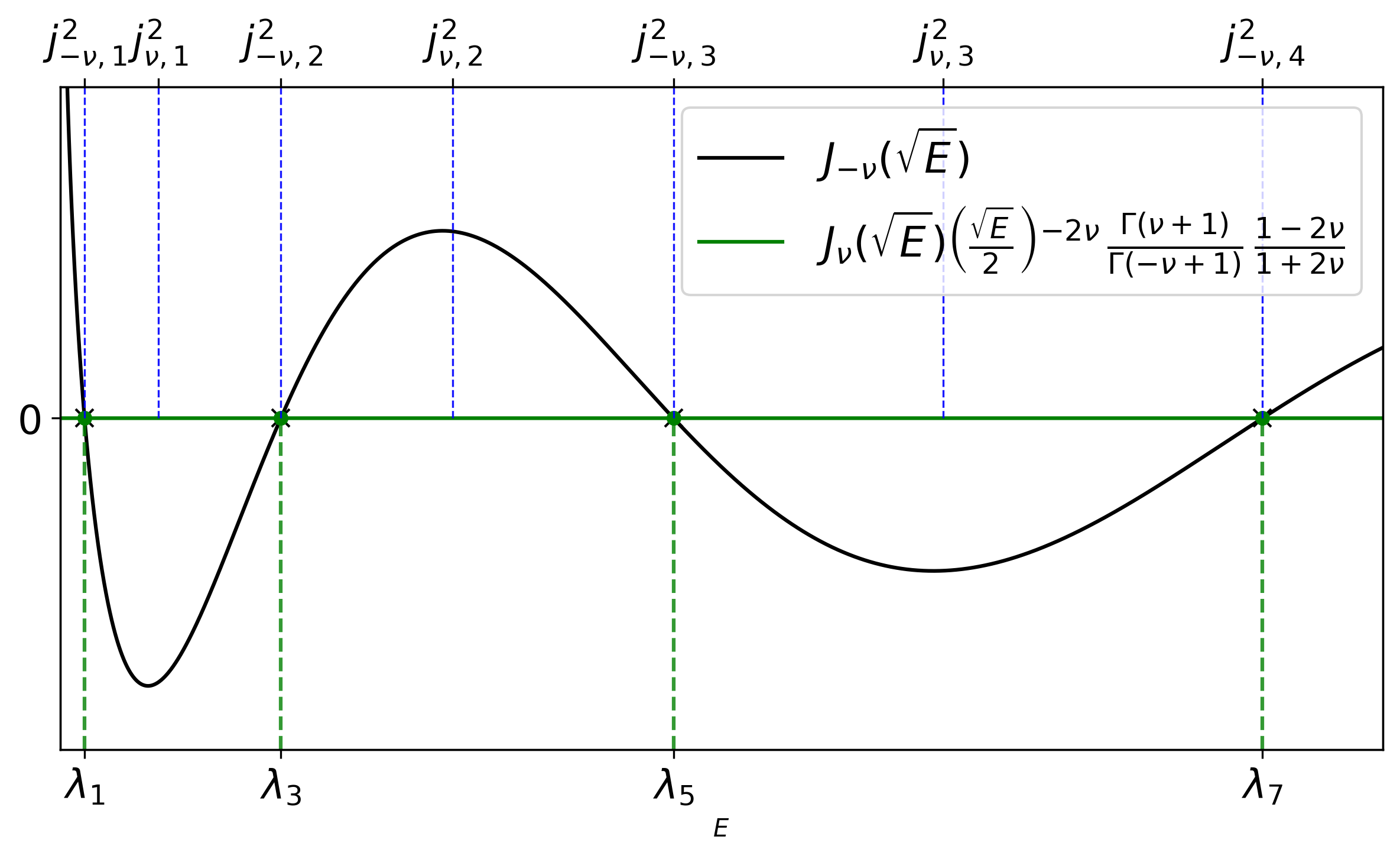}
    \caption{Odd eigenvalues $\lambda_{2n+1}$ for $\nu=1/2$}
    \label{fig:subfig_b}
  \end{subfigure}

  \begin{subfigure}[b]{0.495\textwidth}
    \centering
    \includegraphics[width=\textwidth]{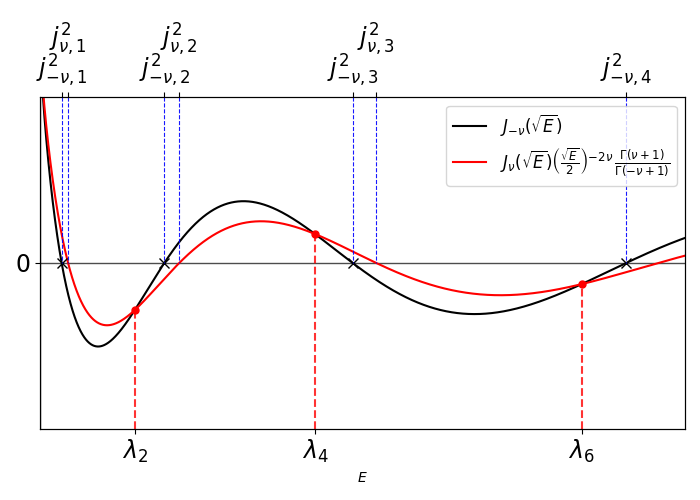}
    \caption{Even eigenvalues $\lambda_{2n}$ for $\nu=0.1\in(0,1/2)$}
    \label{fig:subfig_d}
  \end{subfigure}
  \hfill
  \begin{subfigure}[b]{0.495\textwidth}
    \centering
    \includegraphics[width=\textwidth]{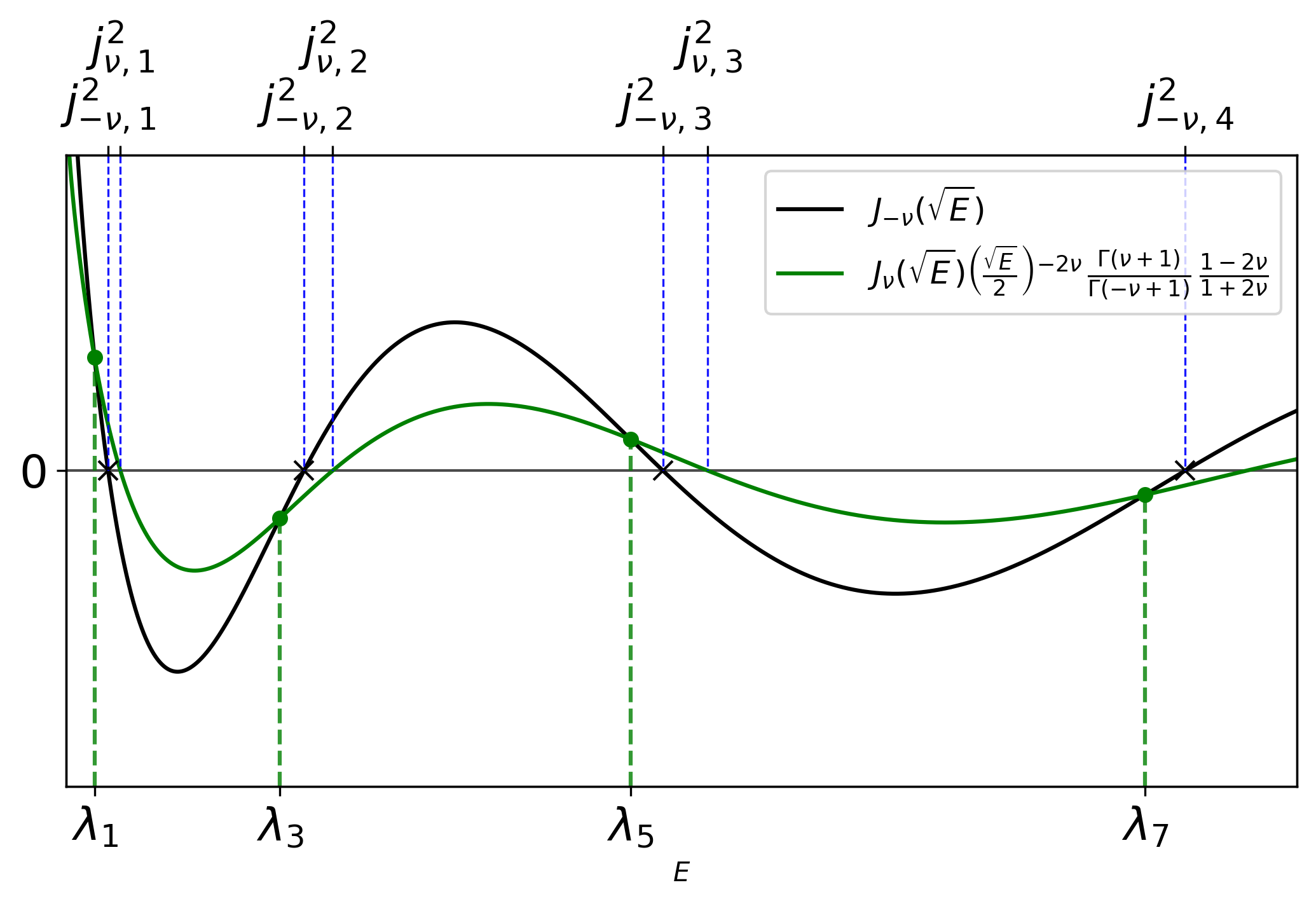}
    \caption{Odd eigenvalues $\lambda_{2n+1}$ for $\nu=0.1\in(0,1/2)$}
    \label{fig:subfig_c}
  \end{subfigure}
  
    \begin{subfigure}[b]{0.495\textwidth}
    \centering
    \includegraphics[width=\textwidth]{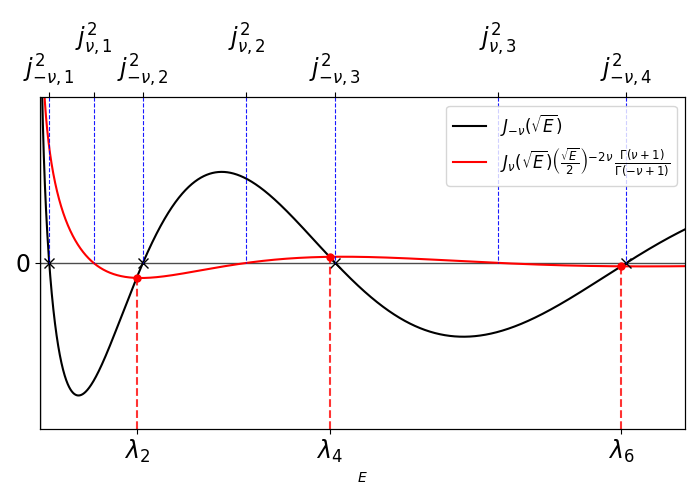}
    \caption{Even eigenvalues $\lambda_{2n}$ for $\nu=0.6\in(1/2,1)$}
    \label{fig:subfig_e}
  \end{subfigure}
  \hfill
  \begin{subfigure}[b]{0.495\textwidth}
    \centering
    \includegraphics[width=\textwidth]{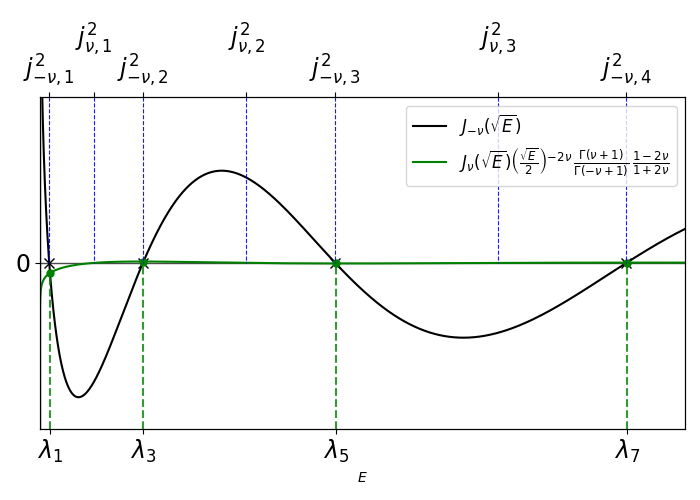}
    \caption{Odd eigenvalues $\lambda_{2n+1}$ for $\nu=0.6\in(1/2,1)$}
    \label{fig:subfig_f}
  \end{subfigure}

  \caption{Illustration of the evolution of eigenvalues as a function of $\nu$.}
  \label{fig:4plots}
\end{figure}

\begin{remark}.
\begin{itemize}
    \item The reader should note that the proposition is stated in such a way that for any fixed $\nu \in (0, 1)$, all values of $E$ that are solutions of (\ref{maincond+}) or (\ref{maincond-}) are listed. We ensure that we provide a complete inventory. Therefore, the point spectrum of $A_\nu$ is $\sigma_{point}=\{\lambda_n,n\in\N\}.$
    \item Proposition \ref{thDistribution} implies that $\lim_{n \to +\infty} \lambda_n = +\infty$, since $\lim_{n \to +\infty} j_{-\nu,n} = +\infty$.
\end{itemize}
    \label{remarkTheorem}
\end{remark}

\begin{remark} Proposition \ref{thDistribution} (and notably \eqref{order}) implies that
\begin{itemize}
    \item if $\nu\in(0,\frac12)$, $\lambda_{2n}<\lambda_{2n+1}<j^2_{-\nu,n+1}, \quad \forall n\in\N\setminus\{0\}$.
    \item if $\nu\in(\frac12,1)$, $\lambda_{2n}<j^2_{-\nu,n+1}<\lambda_{2n+1}, \quad \forall n\in\N\setminus\{0\}$.
    \item if $\nu=\frac12$, $\lambda_{2n}<\lambda_{2n+1}=j^2_{-\nu,n+1}, \quad \forall n\in\N\setminus\{0\}$.
\end{itemize}
    \label{remarkTheorem2}
\end{remark}

Figure~\ref{fig:4plots} illustrates the behavior of the eigenvalues $\lambda_n$ for different values of the parameter $\nu \in (0,1)$, distinguishing between even and odd eigenvalues. These plots visually confirm Proposition~\ref{thDistribution}, by explicitly showing where the solutions to (\ref{maincond+}) and (\ref{maincond-}) lie relative to the squared Bessel zeros $j_{\pm \nu,n}^2$.

The graphical layout shows that the eigenvalues are interlaced with the Bessel zeros. For instance, one observes that between each pair of consecutive Bessel zeros, exactly one eigenvalue may exist, depending on the parity and the value of $\nu$. The figures also show the transitions described in Remark~\ref{remarkTheorem2}: when $\nu$ increases from below $1/2$ to above it, the location of the odd eigenvalues $\lambda_{2n+1}$ shift from being to the left of $j^2_{-\nu,n+1}$ to being to its right.

This transition is particularly noticeable when comparing Figures~\ref{fig:subfig_b}, \ref{fig:subfig_c}, and \ref{fig:subfig_f}. However, in the case $\nu = 0.6$, corresponding to Figure~\ref{fig:subfig_f}, the values of $\lambda_{2n+1}$ and $j_{-\nu,n+1}^2$ are extremely close, and the main plot does not allow a precise visual comparison.

\begin{figure}[ht]
    \centering
    \includegraphics[width=0.4\textwidth]{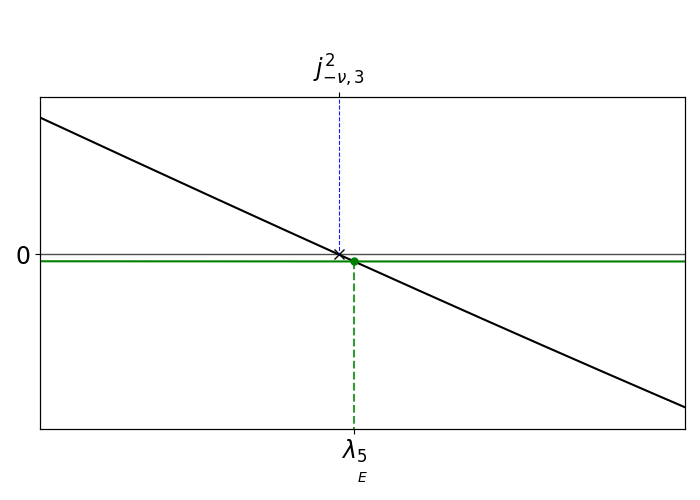}
    \caption{Zoom on the eigenvalue $\lambda_5$ and the squared Bessel zero $j_{-\nu,3}^2$ for $\nu = 0.6\in(1/2,1)$}
    \label{fig:label_de_ma_figure}
\end{figure}

To address this, Figure~\ref{fig:label_de_ma_figure} gives a magnified view centered around $\lambda_5$ and $j_{-\nu,3}^2$. This zoomed-in inset clearly confirms that $\lambda_5 > j_{-\nu,3}^2$, as given by Proposition~\ref{thDistribution} for $\nu \in (1/2,1)$.

\noindent \textbf{\textit{Proof of Proposition \ref{thDistribution}}}. There are four clear blocks in the proposition. We will go through them one by one in Parts 1, 2, 3, and 4 of the proof. First, we should underline that, by Lemma \ref{lemmaInterlaced}, we know that the family
\[
(0, j^2_{-\nu, 1}) \bigcup_{n \geq 1} \{ j^2_{-\nu,n} \} \bigcup_{n \geq 1} \{ j^2_{\nu,n} \} \bigcup_{n \geq 1} (j^2_{-\nu,n}, j^2_{\nu,n}) \bigcup_{n \geq 1} (j^2_{\nu,n}, j^2_{-\nu,n+1})
\]
is a partition of $\mathbb{R}_+\setminus\{0\}$. Therefore, the proposition provides a complete inventory of the solutions $E > 0$ of (\ref{maincond+}) or (\ref{maincond-}). We will use the function $f$ defined in Proposition \ref{coroMonotone}.
\vspace{0.2cm}

\noindent \textbf{Part 1.}  Let us introduce $h_1$, which is well-defined on $(\R_+\setminus\{0\})\bigcup_{n \geq 1} \{ j^2_{-\nu,n} \}$:
\[
h_1:x\mapsto f(\sqrt{x})\frac{1}{2^{-2\nu}}\frac{\Gamma(\nu+1)}{\Gamma(1-\nu)}.
\]
Since $x\mapsto\sqrt{x}$ is increasing on $\R_+\setminus\{0\}$, and $\frac{1}{2^{-2\nu}}\frac{\Gamma(\nu+1)}{\Gamma(1-\nu)} >0$, we have, using Proposition \ref{coroMonotone}, that for any $n\in\N\setminus\{0\}$, $h_1$ is a continuous increasing bijection from $(j^2_{-\nu,n}, j^2_{-\nu,n+1})$ to $\R$, and from $(0, j^2_{-\nu,1})$ to $(\frac{1}{2^{-2\nu}}\frac{\Gamma(\nu+1)}{\Gamma(1-\nu)}\times\frac{\Gamma(1-\nu)}{\Gamma(\nu+1)}2^{-2\nu} , +\infty)=(1, +\infty)$.

\noindent Let $E>0$ satisfying (\ref{maincond+}). Then $E\not\in \bigcup_{n \geq 1} \{ j^2_{-\nu,n} \}$ because otherwise it would imply that $J_\nu(\sqrt{E})=J_{-\nu}(\sqrt{E})=0$, which is not possible according to Lemma \ref{lemmaInterlaced}.

\noindent Therefore, $E>0$ satisfies (\ref{maincond+}) if and only if
\begin{equation}
    E\not\in \bigcup_{n \geq 1} \{ j^2_{-\nu,n} \} \quad 
    \text{and} \quad 
    h_1(E)=1.
\label{cond+withf}
\end{equation}

\noindent Using the intermediate value theorem on $h_1$, we directly prove the first three items of the proposition:
\begin{itemize}
    \item \textit{$\nexists E\in(0, j^2_{-\nu, 1}), \ E$ satisfies (\ref{maincond+})}. Indeed, we cannot have $h_1(E)=1$ on $(0, j^2_{-\nu, 1})$ since $h_1(0)=1$.
    \item \textit{$\forall n\in\N\setminus\{0\}, \ \nexists E\in[j^2_{-\nu,n}, j^2_{\nu,n}], \ E$ satisfies (\ref{maincond+})}. Indeed, for any $n\in\N\setminus\{0\}$, $h_1(j^2_{\nu,n})=0$, so ${h_{1|}}_{(j^2_{-\nu,n}, j^2_{\nu,n}]}$ is a bijection from $(j^2_{-\nu,n}, j^2_{\nu,n}]$ to $(-\infty,0]$. Therefore, if $E$ satisfies (\ref{maincond+}), then $E\not\in(j^2_{-\nu,n}, j^2_{\nu,n}]$. Additionally, $E\neq j^2_{-\nu,n}$ according to (\ref{cond+withf}), so $E\not\in[j^2_{-\nu,n}, j^2_{\nu,n}]$.
    \item \textit{$\forall n\in\N\setminus\{0\}, \ \exists ! E\in(j^2_{\nu,n}, j^2_{-\nu,n+1})$ such that $E$ satisfies (\ref{maincond+})}. Indeed, ${h_{1|}}_{(j^2_{\nu,n}, j^2_{-\nu,n+1})}$ is a bijection from $(j^2_{\nu,n}, j^2_{-\nu,n+1})$ to $(0,+\infty)$.
\end{itemize}
\noindent \textbf{Part 2.} Let $\nu\in(0, \frac12)$. We introduce $h_2$, which is well-defined on $\R_+\setminus\{0\} \bigcup_{n \geq 1} \{ j^2_{-\nu,n} \}$:
\[
h_2:x\mapsto f(\sqrt{x})\frac{1}{2^{-2\nu}}\frac{\Gamma(\nu+1)}{\Gamma(1-\nu)}\frac{1-2\nu}{1+2\nu}.
\]
Since $x\mapsto\sqrt{x}$ is increasing on $\R_+\setminus\{0\}$, and $\frac{1}{2^{-2\nu}}\frac{\Gamma(\nu+1)}{\Gamma(1-\nu)}\frac{1-2\nu}{1+2\nu} >0$, using Proposition \ref{coroMonotone},  we have that for any $n\in\N\setminus\{0\}$, $h_2$ is a continuous increasing bijection from $(j^2_{-\nu,n}, j^2_{-\nu,n+1})$ to $\R$  and from $(0, j^2_{-\nu,1})$ to $(\frac{1}{2^{-2\nu}}\frac{\Gamma(\nu+1)}{\Gamma(1-\nu)}\frac{1-2\nu}{1+2\nu}\times\frac{\Gamma(1-\nu)}{\Gamma(\nu+1)}2^{-2\nu} , +\infty)=(\frac{1-2\nu}{1+2\nu}, +\infty)$.

\noindent Using exactly the same arguments as in Part 1 of the proof, we find that $E>0$ satisfies (\ref{maincond-}) if and only if
\[
    E\not\in \bigcup_{n \geq 1} \{ j^2_{-\nu,n} \}\quad
    \text{and}\quad
    h_2(E)=1.
\]
\noindent Using the intermediate value theorem on $h_2$, we directly obtain the three required items, as done in Part 1. Finally, by Lemma \ref{lemmaInterlaced}, the only thing left to prove to have (\ref{order}) for $\nu\in(0,\frac12)$ is that for all $n\in\N\setminus\{0\}$, $\lambda_{2n}<\lambda_{2n+1}$. This is straightforward since $\lambda_{2n}$ and $\lambda_{2n+1}$ are both in $(j^2_{\nu,n}, j^2_{-\nu,n+1})$, $h_2=\frac{1-2\nu}{1+2\nu}h_1$ so $h_2<h_1$ on $(j^2_{\nu,n}, j^2_{-\nu,n+1})$, and $h_1(\lambda_{2n})=h_2(\lambda_{2n+1})=1$.
\vspace{0.2cm}

\noindent \textbf{Part 3.} Let $\nu\in\left(\frac{1}{2}, 1\right)$. We follow the exact same strategy as in Part 2, except that here, $\frac{1-2\nu}{1+2\nu} <0$, so the growth behavior of $h_2$ is inverted. The increasing parts are now turned into decreasing ones, and vice versa. Using the intermediate value theorem on $h_2$, we directly prove the three needed items. Finally, by Lemma \ref{lemmaInterlaced}, we have nothing left to prove to establish (\ref{order}) for $\nu\in\left(\frac{1}{2},1\right)$.
\vspace{0.2cm}

\noindent \textbf{Part 4.} Let $\nu=\frac{1}{2}$. For $\nu=\frac{1}{2}$, we directly find that the fact that $E>0$ satisfies (\ref{maincond-}) is equivalent to $J_{-\nu}(\sqrt{E})=0$. Therefore, the two items of this part are directly obtained. Finally, by Lemma \ref{lemmaInterlaced}, we have nothing left to prove in order to establish (\ref{order}) for $\nu=\frac{1}{2}$.
\qed

\subsection{Hilbert Basis of eigenfunctions}

In order to prove that the spectrum of $A_\nu$ is restricted to the point spectrum, we prove the following result.

\begin{proposition} $(A_\nu, D(A_\nu))$ has compact resolvent.
    \label{CorrCompactResolvant}
\end{proposition}

\noindent \textbf{\textit{Proof of proposition \ref{CorrCompactResolvant}}}. Since $(A_\nu, D(A_\nu))$ is self-adjoint, it is automatically closed. Hence, for any $\lambda \in \mathbb R$ such that $A_\nu-\lambda \mathrm{Id}$ is invertible, the closed graph theorem ensures that $(A_\nu-\lambda \mathrm{Id})^{-1} \in \mathcal L_c(L^2(-1,1),D(A_\nu))$, where $D(A_\nu)$ is equipped with the usual graph norm. So, by composition, $(A_\nu, D(A_\nu))$ has a compact resolvent as soon as the injection \((D(A_\nu), \norm{\cdot}_{D(A_\nu)}) \hookrightarrow (L^2(-1,1), \norm{\cdot}_{L^2(-1,1)})\) is compact.

\noindent Let \((f_n=f_r^n+f_s^n)_{n\in\N} \in D(A_\nu)^\N\) be a bounded sequence for the graph norm \(\norm{\cdot}_{D(A_\nu)}\)\footnote{In all what follows, for any set $E$, $E^{\N}$ stands for sequences of elements of $E$ indexed by $\N$. The notation $E^{\N\setminus \{0\}}$ is employed in a similar way.}, where for all \(n\in\N\), \(f_r^n\in\tilde{H}_0^2(-1,1)\) and \(f_s^n\in\mathcal{F}_s^\nu\). We have that
\[
\exists C_1>0, \ \forall n\in\N, \ \norm{f_n}_{L^2(-1,1)}\leq C_1, \quad \mathrm{and} \quad 
\exists C_2>0, \ \forall n\in\N, \ \norm{Af_n}_{L^2(-1,1)}\leq C_2.
\]
Therefore, we obtain
\[
\forall n\in\N, \quad \langle Af_n,f_n\rangle_{L^2(-1,1)} \leq \norm{f_n}_{L^2(-1,1)}\norm{Af_n}_{L^2(-1,1)}\leq C_1 C_2.
\]
Using (\ref{1.6}), we have
\[
\forall n\in\N, \quad \min\{1,4\nu^2\}\int_{-1}^1(\partial_x f^n_r)^2 \leq C_1 C_2.
\]
Thus, we conclude that \((\partial_x f^n_r)_{n\in\N}\) is a bounded sequence in \(L^2(-1,1)\). Combining this with Poincar\'e's inequality on (-1,1) (which is valid since $f^n_r\in \tilde{H}_0^2(-1,1)$ for any $n\in\N$), we deduce that \((f^n_r)_{n\in\N}\) is also a bounded sequence in \(L^2(-1,1)\). This implies that \((f^n_r)_{n\in\N}\) is a bounded sequence in \(H^1(-1,1)\), so that Rellich's theorem states that the injection \(H^1(-1,1) \hookrightarrow L^2(-1,1)\) is compact. Therefore, there exists a subsequence \((f^{\varphi_1(n)}_r)_{n\in\N}\) converging strongly in \(L^2(-1,1)\) to some \(f_r\in L^2(-1,1)\). Furthermore,
\[
\forall n\in\N, \quad \norm{f^n_s}_{L^2(-1,1)} = \norm{f^n_s + f^n_r - f^n_r}_{L^2(-1,1)} \leq \norm{f^n}_{L^2(-1,1)} + \norm{f^n_r}_{L^2(-1,1)}.
\]
Since \((f^n)_{n\in\N}\) is bounded in \(L^2(-1,1)\), and we have shown that \((f_r^n)_{n\in\N}\) is bounded in \(L^2(-1,1)\), we conclude that \((f^n_s)_{n\in\N}\) is also bounded in \(L^2(-1,1)\), so \((f^{\varphi_1(n)}_s)_{n\in\N}\) is also bounded. Moreover, \((f_s^{\varphi_1(n)})_{n\in\N} \in (\mathcal{F}_s^\nu)^\N\), and since \(\mathcal{F}_s^\nu\) is a finite-dimensional vector subspace of \(L^2(-1,1)\), there exists a subsequence \((f_s^{\varphi_2\circ\varphi_1(n)})_{n\in\N}\) that converges strongly in \(L^2(-1,1)\) to some \(f_s\in\mathcal{F}_s^\nu \subset L^2(-1,1)\). Thus, we have shown that
\[
\forall n\in\N, \quad f_{\varphi_2\circ\varphi_1(n)} = f_r^{\varphi_2\circ\varphi_1(n)} + f_s^{\varphi_2\circ\varphi_1(n)} \underset{n\to+\infty}{\longrightarrow} f_r + f_s \in L^2(-1,1).
\]
Therefore, \((D(A_\nu), \norm{\cdot}_{D(A_\nu)}) \hookrightarrow (L^2(-1,1), \norm{\cdot}_{L^2(-1,1)})\) is compact, and we conclude as explained at the beginning of the proof.
\qed

\vspace{0.3cm}

\noindent We are now ready to give a complete description of the spectral theory of $A_\nu$.

\begin{theorem} We introduce
\[
    \psi_{0}:=|x|^{\nu+\frac12}-|x|^{-\nu+\frac12} \ \mathrm{on} \ (-1,1), \quad \text{and} \quad
a_0:=\norm{\psi_0}_{L^2(-1,1)}.
\]
Let $n\in\N\setminus\{0\}$. We introduce
\[
\psi_{2n}:=
    \begin{cases}
    -\dfrac{\Gamma(\nu+1)}{\Gamma(1-\nu)}\left(\dfrac{\sqrt{\lambda_{2n}}}{2}\right)^{-2\nu}\sqrt{-x}J_\nu\left(-\sqrt{\lambda_{2n}}x \right)+\sqrt{-x}J_{-\nu}\left(-\sqrt{\lambda_{2n}}x\right),
     \ &\mathrm{on}\ (-1,0),
    \\
    -\dfrac{\Gamma(\nu+1)}{\Gamma(1-\nu)}\left(\dfrac{\sqrt{\lambda_{2n}}}{2}\right)^{-2\nu}\sqrt{x}J_\nu\left(\sqrt{\lambda_{2n}}x\right)+\sqrt{x}J_{-\nu}\left(\sqrt{\lambda_{2n}}x\right),
    \ &\mathrm{on}\ (0,1),
\end{cases}
    \]
\begin{equation}
\begin{split}
     &a_{2n}^2:=\norm{\psi_{2n}}^2_{L^2(-1,1)}=\left(1-\dfrac{\nu^2}{\lambda_{2n}}\right)J_{-\nu}\left(\sqrt{\lambda_{2n}}\right)^2+{J_{-\nu}}'\left(\sqrt{\lambda_{2n}}\right)^2\\
        &\begin{split}
        -2\dfrac{\Gamma(\nu+1)}{\Gamma(1-\nu)}\left(\dfrac{\sqrt{\lambda_{2n}}}{2}\right)^{-2\nu}\left[\left(1-\dfrac{\nu^2}{\lambda_{2n}}\right)\right.&J_{\nu}\left(\sqrt{\lambda_{2n}}\right)J_{-\nu}\left(\sqrt{\lambda_{2n}}\right)\\
        &\left.+ {J_{\nu}}'\left(\sqrt{\lambda_{2n}}\right){J_{-\nu}}'\left(\sqrt{\lambda_{2n}}\right) + \frac{2\nu\sin(\nu\pi)}{\pi \lambda_{2n}}\right]    
        \end{split}\\
        &+\left(\dfrac{\Gamma(\nu+1)}{\Gamma(1-\nu)}\left(\dfrac{\sqrt{\lambda_{2n}}}{2}\right)^{-2\nu}\right)^2\left[\left(1-\dfrac{\nu^2}{\lambda_{2n}}\right)J_{\nu}\left(\sqrt{\lambda_{2n}}\right)^2+{J_{\nu}}'\left(\sqrt{\lambda_{2n}}\right)^2\right] >0
        \end{split}
    \label{a2n}
\end{equation}
Let $n\in\N$, we introduce
\[\psi_{2n+1}:=
    \begin{cases}
    \dfrac{1-2\nu}{1+2\nu}\dfrac{\Gamma(\nu+1)}{\Gamma(1-\nu)}\left(\dfrac{\sqrt{\lambda_{2n+1}}}{2}\right)^{-2\nu}\sqrt{-x}J_\nu\left(-\sqrt{\lambda_{2n+1}}x\right)-\sqrt{-x}J_{-\nu}\left(-\sqrt{\lambda_{2n+1}}x\right),  &\mathrm{on}(-1,0),
    \\
    -\dfrac{1-2\nu}{1+2\nu}\dfrac{\Gamma(\nu+1)}{\Gamma(1-\nu)}\left(\dfrac{\sqrt{\lambda_{2n+1}}}{2}\right)^{-2\nu}\sqrt{x}J_\nu\left(\sqrt{\lambda_{2n+1}}x\right)+\sqrt{x}J_{-\nu}\left(\sqrt{\lambda_{2n+1}}x\right),
     \ &\mathrm{on}(0,1),
\end{cases}
    \]
\begin{equation}
\begin{split}
&a_{2n+1}^2:=\norm{\psi_{2n+1}}^2_{L^2(-1,1)}=\left(1-\frac{\nu^2}{\lambda_{2n+1}}\right)J_{-\nu}\left(\sqrt{\lambda_{2n+1}}\right)^2+{J_{-\nu}}'\left(\sqrt{\lambda_{2n+1}}\right)^2\\
        &\begin{split}
        -2\frac{1-2\nu}{1+2\nu}\frac{\Gamma(\nu+1)}{\Gamma(1-\nu)}
        \left(\frac{\sqrt{\lambda_{2n+1}}}{2}\right)^{-2\nu}\left[\left(1-\frac{\nu^2}{\lambda_{2n+1}}\right)\right.&J_{\nu}\left(\sqrt{\lambda_{2n+1}}\right)J_{-\nu}\left(\sqrt{\lambda_{2n+1}}\right).\\
        &\left.+ {J_{\nu}}'\left(\sqrt{\lambda_{2n+1}}\right){J_{-\nu}}'\left(\sqrt{\lambda_{2n+1}}\right)+ \frac{2\nu\sin(\nu\pi)}{\pi \lambda_{2n+1}}\right]
        \end{split}\\
        &+\left(\frac{1-2\nu}{1+2\nu}\frac{\Gamma(\nu+1)}{\Gamma(1-\nu)}\left(\frac{\sqrt{\lambda_{2n+1}}}{2}\right)^{-2\nu}\right)^2\left[\left(1-\frac{\nu^2}{\lambda_{2n+1}}\right)J_{\nu}\left(\sqrt{\lambda_{2n+1}}\right)^2+{J_{\nu}}'\left(\sqrt{\lambda_{2n+1}}\right)^2\right]>0.
        \end{split}    
        \label{a2n+1}
\end{equation}
    \noindent Let $n\in\N$, define
\[
\phi_n:=\frac{\psi_{n}}{a_n}.
\]
Then, $(\phi_n)_{n\in\N}$ forms a Hilbert basis of $L^2(-1,1)$ consisting of eigenfunctions of $A_\nu$.
    \label{propEigenBasis}
\end{theorem}

\noindent \textbf{\textit{Proof of Theorem \ref{propEigenBasis}}}. Using Proposition \ref{PropKernel} and Corollary \ref{PropCondSimpler}, we know that $\sigma_{point}(A_\nu)=\{\lambda_n,\ n\in\N\}$, and that all these eigenvalues are distinct 
and have multiplicity 1. Moreover, for any $n\in\N$, $\mathrm{Ker}(A_\nu-\lambda_n \ \mathrm{Id})=\mathrm{Span}(\psi_n)=\mathrm{Span}(\phi_n)$. Furthermore, using  (\ref{intSame-NuSameA}) and (\ref{intNotSameNuSameA}), we directly obtain 
\eqref{a2n} and \eqref{a2n+1}. Since $A_\nu$ is self-adjoint with compact resolvent, we deduce that $\sigma(A_\nu)=\{\lambda_n,\ n\in\N\}$ and that the sequence $(\phi_n)_{n\in\N}$ forms a Hilbert basis of $L^2(-1,1)$.
\qed

\section{Asymptotic behavior of the eigenvalues}
\label{s:as}
\subsection{Case $\nu\in\left(0, \frac12\right)$}

\noindent Throughout this subsection, we focus on the case $\nu\in(0,\frac12)$.
\begin{lemma} Let $\nu\in\left(0,\frac12\right)$,
    \begin{equation}
    \begin{split}
            \sqrt{\frac{\pi x}{2}}\left(J_\nu(x) \left(\frac{x}{2}\right)^{-2\nu}\frac{\Gamma(\nu+1)}{\Gamma(1-\nu)}-J_{-\nu}(x)\right)\underset{x\to+\infty}{=}&-\cos\left(x+\nu\frac{\pi}{2}-\frac{\pi}{4}\right)\\
            &+\frac{\Gamma(\nu+1)}{\Gamma(1-\nu)}\cos\left(x-\nu\frac{\pi}{2}-\frac{\pi}{4}\right)\left(\frac{2}{x}\right)^{2\nu}\\
            &+\mathcal{O}\left(\frac{1}{x}\right),
            \end{split}
            \label{troncated(0,1/2)3.20}
    \end{equation}
    \begin{equation}
    \begin{split}
            \sqrt{\frac{\pi x}{2}}\left(J_\nu(x) \left(\frac{x}{2}\right)^{-2\nu}\frac{1-2\nu}{1+2\nu}\frac{\Gamma(\nu+1)}{\Gamma(1-\nu)}-J_{-\nu}(x)\right)\underset{x\to+\infty}{=}&-\cos\left(x+\nu\frac{\pi}{2}-\frac{\pi}{4}\right)\\
            &+\frac{1-2\nu}{1+2\nu}\frac{\Gamma(\nu+1)}{\Gamma(1-\nu)}\cos\left(x-\nu\frac{\pi}{2}-\frac{\pi}{4}\right)\left(\frac{2}{x}\right)^{2\nu}\\
            &+\mathcal{O}\left(\frac{1}{x}\right).
            \end{split}
            \label{troncated(0,1/2)3.21}
    \end{equation}
        \label{lemmaTroncatedExpansion(0,1/2)}
    \end{lemma}
    
    \noindent \textbf{\textit{Proof of Lemma \ref{lemmaTroncatedExpansion(0,1/2)}}}. Let $\nu\in\left(0,\frac12\right)$, by Lemma \ref{lemmaExpansion}, we have
    
    \begin{equation}
        \sqrt{\frac{\pi x}{2}}J_{-\nu}(x)\underset{x\to+\infty}{=}\cos\left(x+\nu\frac{\pi}{2}-\frac{\pi}{4}\right)+\mathcal{O}\left(\frac{1}{x}\right),
        \label{prooftronc1-0,1/2}
    \end{equation}
    \begin{equation}
        \sqrt{\frac{\pi x}{2}}J_{\nu}(x)\left(\frac{x}{2}\right)^{-2\nu}\underset{x\to+\infty}{=}
        \cos\left(x-\nu\frac{\pi}{2}-\frac{\pi}{4}\right)\left(\frac{2}{x}\right)^{2\nu}+\mathcal{O}\left(\frac{1}{x^{1+2\nu}}\right).
        \label{prooftronc2-0,1/2}
    \end{equation}
    
    \noindent $\dfrac{\Gamma(\nu+1)}{\Gamma(1-\nu)}\times$(\ref{prooftronc2-0,1/2})$-$(\ref{prooftronc1-0,1/2})$=$(\ref{troncated(0,1/2)3.20}), \quad \text{and} \quad  $\dfrac{1-2\nu}{1+2\nu}\dfrac{\Gamma(\nu+1)}{\Gamma(1-\nu)}\times$(\ref{prooftronc2-0,1/2})$-$(\ref{prooftronc1-0,1/2})$=$(\ref{troncated(0,1/2)3.21}).
    \qed
    
    \begin{proposition} Let $\nu\in\left(0,\frac12\right)$,
    \begin{flalign}
        \label{sqrtEven(0,1/2)}
    \sqrt{\lambda_{2(n-1)}}\underset{n\to+\infty}{=}&\pi\left(n-\frac{\nu}{2}-\frac14\right)-\frac{\Gamma(\nu+1)}{\Gamma(1-\nu)}\sin(\nu\pi)\left(\frac{2}{\pi n}\right)^{2\nu}
    +\mathcal{O}\left(\frac{1}{n^{\min(1,4\nu)}}\right),&&
    \end{flalign}
    \begin{flalign}
        \label{sqrtOdd(0,1/2)}
    \sqrt{\lambda_{2(n-1)+1}}\underset{n\to+\infty}{=}&\pi\left(n-\frac{\nu}{2}-\frac14\right)-\frac{1-2\nu}{1+2\nu}\frac{\Gamma(\nu+1)}{\Gamma(1-\nu)}\sin(\nu\pi)\left(\frac{2}{\pi n}\right)^{2\nu}
    + \mathcal{O}\left(\frac{1}{n^{\min(1,4\nu)}}\right).&&
    \end{flalign}
    \label{propAsymptoticbehaviors(0,1/2)}
    \end{proposition}
    
\noindent \textbf{\textit{Proof of Proposition \ref{propAsymptoticbehaviors(0,1/2)}}}. Let $\nu\in\left(0,\frac12\right)$.
    
\noindent \textbf{Proof of (\ref{sqrtEven(0,1/2)})}. 

\noindent \textbf{Step 1.} Using (\ref{troncated(0,1/2)3.20}), we know that
\[
\sqrt{\frac{\pi x}{2}}\left(J_\nu(x) \left(\frac{x}{2}\right)^{-2\nu}\frac{\Gamma(\nu+1)}{\Gamma(1-\nu)}-J_{-\nu}(x)\right)\underset{x\to+\infty}{=}-\cos\left(x+\nu\frac{\pi}{2}-\frac{\pi}{4}\right)+\mathcal{O}\left(\frac{1}{x^{2\nu}}\right).
\]
We also need to have for $n>1$
\begin{equation}
    \sqrt{\frac{\pi \sqrt{\lambda_{2(n-1)}}}{2}}\left(J_\nu\left( \sqrt{\lambda_{2(n-1)}}\right) \left(\frac{ \sqrt{\lambda_{2(n-1)}}}{2}\right)^{-2\nu}\frac{\Gamma(\nu+1)}{\Gamma(1-\nu)}-J_{-\nu}\left( \sqrt{\lambda_{2(n-1)}}\right)\right)=0.
    \label{proof3.20-0,1/2}
\end{equation}
Thus we must have
\begin{equation}
    0\underset{n\to+\infty}{=}-\cos\left(\sqrt{\lambda_{2(n-1)}}+\nu\frac{\pi}{2}-\frac{\pi}{4}\right)+\mathcal{O}\left(\frac{1}{\sqrt{\lambda_{2(n-1)}}}\right),
    \label{proofverified-0,1/2}
\end{equation}
because $\sqrt{\lambda_{2(n-1)}}\underset{n\to+\infty}{\longrightarrow}+\infty$ as explained in Remark \ref{remarkTheorem}. Thus, we have
\[
\cos\left(\sqrt{\lambda_{2(n-1)}}+\nu\frac{\pi}{2}-\frac{\pi}{4}\right)\underset{n\to+\infty}{=}o(1).
\]
That implies that there must exist a sequence of integers $(k_n)_{n\in\N\setminus\{0\}}\in\Z^{\N\setminus\{0\}}$ such that
$$\sqrt{\lambda_{2(n-1)}}\underset{n\to+\infty}{=}\pi\left(k_n-\frac{\nu}{2}-\frac14\right)+o(1),$$
and using Proposition \ref{thDistribution} we know that for $n>1$, $\sqrt{\lambda_{2(n-1)}}\in(j_{\nu,n-1},j_{-\nu,n})$. Using Lemma \ref{lemmaAsymptoticZero}, we obtain, as $n\rightarrow +\infty$,
\[
\sqrt{\lambda_{2(n-1)}}{\in}\left(\pi\left(n-1+\frac{\nu}{2}-\frac14\right)+\mathcal{O}\left(\frac1n\right),\quad \pi\left(n-\frac{\nu}{2}-\frac14\right)+\mathcal{O}\left(\frac1n\right)\right).
\]
Therefore, since $\nu\in(0,1)$, for $n$ large enough, we must have $k_n=n$. This proves that
\begin{equation}
    \sqrt{\lambda_{2(n-1)}}\underset{n\to+\infty}{=}\pi\left(n-\frac{\nu}{2}-\frac14\right)+\varepsilon_{2\nu}(n) \quad \text{with} \ \varepsilon_{2\nu}(n)= o(1).
    \label{prooffirststep-0,1/2}
\end{equation}
\noindent \textbf{Step 2.} Now, we prove that $\varepsilon_{2\nu}(n)\underset{n\to+\infty}{=} \mathcal{O}\left(\dfrac1{n^{2\nu}}\right)$. We repeat the same argument. Using (\ref{proofverified-0,1/2}) and (\ref{prooffirststep-0,1/2}), we obtain
\[
\quad \cos\left(\pi n +\varepsilon_{2\nu}(n) -\dfrac{\pi}{2}\right)\underset{n\to+\infty}{=}\mathcal{O}\left(\dfrac1{n^{2\nu}}\right).
\]
\noindent Therefore, 
\[\quad \sin\left(\varepsilon_{2\nu}(n)\right)\underset{n\to+\infty}{=}\mathcal{O}\left(\dfrac1{n^{2\nu}}\right),\quad \mathrm{which\ gives\quad} \quad\varepsilon_{2\nu}(n)+\mathcal{O}\left(\varepsilon_{2\nu}(n)^3\right)\underset{n\to+\infty}{=}\mathcal{O}\left(\dfrac1{n^{2\nu}}\right).
\]
\noindent Thus,
\[\quad \varepsilon_{2\nu}(n)\left(1+\mathcal{O}\left(\varepsilon_{2\nu}(n)^2\right)\right)\underset{n\to+\infty}{=}\mathcal{O}\left(\dfrac1{n^{2\nu}}\right),\quad
\]
which implies that 
\noindent Thus,
\[\quad \varepsilon_{2\nu}(n)\left(1+o(1)\right)\underset{n\to+\infty}{=}\mathcal{O}\left(\dfrac1{n^{2\nu}}\right).\quad
\]
This implies that $\varepsilon_{2\nu}(n)\underset{n\to+\infty}{=}\mathcal{O}\left(\dfrac1{n^{2\nu}}\right)$.

\noindent \textbf{Step 3.} Now we  obtain an even more precise expression of $\varepsilon_{2\nu}$. We will use the same strategy. Using (\ref{troncated(0,1/2)3.20}), we know that
\[
\begin{split}
        \sqrt{\frac{\pi x}{2}}\left(J_\nu(x) \left(\frac{x}{2}\right)^{-2\nu}\frac{\Gamma(\nu+1)}{\Gamma(1-\nu)}-J_{-\nu}(x)\right)\underset{x\to+\infty}{=}&-\cos\left(x+\nu\frac{\pi}{2}-\frac{\pi}{4}\right)\\
        &+\frac{\Gamma(\nu+1)}{\Gamma(1-\nu)}\cos\left(x-\nu\frac{\pi}{2}-\frac{\pi}{4}\right)\left(\frac{2}{x}\right)^{2\nu}+\mathcal{O}\left(\frac{1}{x}\right).
        \end{split}
\]
\noindent (\ref{proof3.20-0,1/2}) is still true. Therefore, we must have
\[
\begin{split}
    &-\cos\left(\sqrt{\lambda_{2(n-1)}}+\nu\frac{\pi}{2}-\frac{\pi}{4}\right)
     \frac{\Gamma(\nu+1)}{\Gamma(1-\nu)}\cos\left(\sqrt{\lambda_{2(n-1)}}-\nu\frac{\pi}{2}-\frac{\pi}{4}\right)\left(\frac{2}{\sqrt{\lambda_{2(n-1)}}}\right)^{2\nu} \\&\underset{n\to+\infty}{=}
    \mathcal{O}\left(\frac{1}{\sqrt{\lambda_{2(n-1)}}}\right)
    .
    \end{split}
\]
Using (\ref{prooffirststep-0,1/2}), we obtain
\[
   (-1)^{n+1}\sin\left(\varepsilon_{2\nu}(n)\right) \ +\ \frac{\Gamma(\nu+1)}{\Gamma(1-\nu)}(-1)^n\sin\left(-\nu\pi+\varepsilon_{2\nu}(n)\right)\left(\frac{2}{\pi n}\right)^{2\nu}\left(1+\mathcal{O}\left(\frac{1}{n}\right)\right)^{-2\nu} \underset{n\to+\infty}{=}\mathcal{O}\left(\frac{1}{n}\right).
\]
\noindent We now use the asymptotic expansions of $\sin(x)$, $(1+x)^{-2\nu}$, and $\sin(-\nu\pi + x)$ around $0$. Recall that we already know that $\varepsilon_{2\nu}(n)\underset{n\to+\infty}{=}\mathcal{O}\left(\dfrac{1}{n^{2\nu}}\right)$ by the previous step, and that since $\nu\in\left(0,\dfrac12\right)$, $2\nu\in(0,1)$. We obtain
\[
\begin{split}  &\frac{\Gamma(\nu+1)}{\Gamma(1-\nu)}(-1)^n\left(\sin(-\nu\pi)+\mathcal{O}\left(\dfrac{1}{n^{2\nu}}\right)\right)\left(\frac{2}{\pi n}\right)^{2\nu}\left(1+\mathcal{O}\left(\dfrac{1}{n}\right)\right)
    \\&+   (-1)^{n+1}\left(\varepsilon_{2\nu}(n)+\mathcal{O}\left(\dfrac{1}{n^{6\nu}}\right)\right)\underset{n\to+\infty}{=}\mathcal{O}\left(\frac{1}{n}\right),
    \end{split}
\]
which yields
\[
\varepsilon_{2\nu}(n)\underset{x\to+\infty}{=}-\frac{\Gamma(\nu+1)}{\Gamma(1-\nu)}\sin(\nu\pi)\left(\frac{2}{\pi n}\right)^{2\nu}
+\mathcal{O}\left(\frac{1}{n^{\min(1,4\nu)}}\right),
\]
and this concludes the proof of (\ref{sqrtEven(0,1/2)}).
\vspace{0.2cm}

\noindent \textbf{Proof of (\ref{sqrtOdd(0,1/2)})}. It is crucial to note that the asymptotic expansion (\ref{troncated(0,1/2)3.21}) is simply the asymptotic expansion of (\ref{troncated(0,1/2)3.20}) but with $\dfrac{\Gamma(\nu+1)}{\Gamma(\nu-1)}$ replaced by $\dfrac{1-2\nu}{1+2\nu}\dfrac{\Gamma(\nu+1)}{\Gamma(\nu-1)}$ . Therefore, if we could prove that
    \begin{equation}
            \sqrt{\lambda_{2(n-1)+1}}\underset{n\to+\infty}{=}\pi\left(n-\frac{\nu}{2}-\frac14\right)+o(1),
        \label{step1goal}
        \end{equation}
    then we would follow the exact same proof as for (\ref{sqrtEven(0,1/2)}) except that $\dfrac{\Gamma(\nu+1)}{\Gamma(\nu-1)}$ would always be replaced by $\dfrac{1-2\nu}{1+2\nu}\dfrac{\Gamma(\nu+1)}{\Gamma(\nu-1)}$ . Therefore, we would have proven (\ref{sqrtOdd(0,1/2)}).\\
    \noindent Using (\ref{troncated(0,1/2)3.21}) instead of (\ref{troncated(0,1/2)3.20}) and for $n\in\N\setminus\{0\}$, the condition that $\sqrt{\lambda_{2(n-1)+1}}$ satisfies is now (\ref{maincond-}) instead of (\ref{maincond+}), we can follow the same computations as in Step 1 to obtain that there must exist a sequence of integers $(k_n)_{n\in\N\setminus\{0\}}\in\Z^{\N\setminus\{0\}}$ such that 
    $$\sqrt{\lambda_{2(n-1)+1}}\underset{n\to+\infty}{=}\pi\left(k_n-\frac{\nu}{2}-\frac14\right)+o(1).$$
    \noindent Since $\nu\in\left(0,\frac12\right)$, Proposition \ref{thDistribution} gives that for $n>1$, $\sqrt{\lambda_{2(n-1)+1}}\in(j_{\nu,n-1},j_{-\nu,n})$. Thus, Performing similar computations as in Step 1, we obtain (\ref{step1goal}). As explained before, this is enough to conclude the proof of (\ref{sqrtOdd(0,1/2)}).
\qed
\vspace{0.3cm}

\noindent
Taking the square of (\ref{sqrtEven(0,1/2)}) and (\ref{sqrtOdd(0,1/2)}) and keeping only the terms of the right order of precision leads to the following result.

    \begin{corollary} \label{Corasvp} Let $\nu\in\left(0,\frac12\right)$,
        \begin{flalign}
            \label{even(0,1/2)}
            \lambda_{2(n-1)}\underset{n\to+\infty}{=}&\pi^2\left(n-\frac{\nu}{2}-\frac14\right)^2-4\frac{\Gamma(\nu+1)}{\Gamma(1-\nu)}\sin(\nu\pi)\left(\frac{2}{\pi n}\right)^{2\nu-1}
        + \mathcal{O}\left(\frac{1}{n^{\min(0,4\nu-1)}}\right),&&
    \end{flalign}
    \begin{flalign}
                \label{odd(0,1/2)}
    \lambda_{2(n-1)+1}
    \underset{n\to+\infty}{=}&\pi^2\left(n-\frac{\nu}{2}-\frac14\right)^2-4\frac{1-2\nu}{1+2\nu}\frac{\Gamma(\nu+1)}{\Gamma(1-\nu)}\sin(\nu\pi)\left(\frac{2}{\pi n}\right)^{2\nu-1}
        + \mathcal{O}\left(\frac{1}{n^{\min(0,4\nu-1)}}\right),&&
    \end{flalign}
    \begin{flalign}
            \label{OddMinusEven(0,1/2)}
    \lambda_{2(n-1)+1}-\lambda_{2(n-1)}\underset{n\to+\infty}{=}\frac{16\nu}{1+2\nu}\frac{\Gamma(\nu+1)}{\Gamma(1-\nu)}\sin(\nu\pi)\left(\frac{2}{\pi n}\right)^{2\nu-1} +  \mathcal{O}\left(\frac{1}{n^{\min(0,4\nu-1)}}\right).&&
    \end{flalign}
    
    \label{corollary1(0,1/2)}
    \end{corollary}

\subsection{Case $\nu\in\left[\frac12, 1\right)$}

\noindent During this whole subsection, we will focus only on the case $\nu\in[\frac12,1)$.

\begin{lemma} Let $\nu\in\left[\frac12,1\right)$,
    \begin{equation}
    \begin{split}
            \sqrt{\frac{\pi x}{2}}\left(J_\nu(x) \left(\frac{x}{2}\right)^{-2\nu}\frac{\Gamma(\nu+1)}{\Gamma(1-\nu)}-J_{-\nu}(x)\right)\underset{x\to+\infty}{=}&-\cos\left(x+\nu\frac{\pi}{2}-\frac{\pi}{4}\right)\\
            &+\sin\left(x+\nu\frac{\pi}{2}-\frac{\pi}{4}\right)\frac{4\nu^2-1}{8}\times\frac1x\\
            &+\frac{\Gamma(\nu+1)}{\Gamma(1-\nu)}\cos\left(x-\nu\frac{\pi}{2}-\frac{\pi}{4}\right)\left(\frac{2}{x}\right)^{2\nu}\\
            &+\mathcal{O}\left(\frac{1}{x^2}\right).
            \end{split}
            \label{troncated(1/2,1)3.20}
    \end{equation}
        \label{lemmaTroncatedExpansion(1/2,1)}
    \end{lemma}
    
    \noindent \textbf{\textit{Proof of Lemma \ref{lemmaTroncatedExpansion(1/2,1)}}}. Let $\nu\in\left[\frac12,1\right)$, by Lemma \ref{lemmaExpansion}, we have
    \begin{equation}
        \sqrt{\frac{\pi x}{2}}J_{-\nu}(x)\underset{x\to+\infty}{=}\cos\left(x+\nu\frac{\pi}{2}-\frac{\pi}{4}\right)
        -\sin\left(x+\nu\frac{\pi}{2}-\frac{\pi}{4}\right)\frac{4\nu^2-1}{8}\times\frac1x
        +\mathcal{O}\left(\frac{1}{x^2}\right)
        \label{prooftronc1}
    \end{equation}
    and
    \begin{equation}
        \sqrt{\frac{\pi x}{2}}J_{\nu}(x)\left(\frac{x}{2}\right)^{-2\nu}\underset{x\to+\infty}{=}
        \cos\left(x-\nu\frac{\pi}{2}-\frac{\pi}{4}\right)\left(\frac{2}{x}\right)^{2\nu}
            +\mathcal{O}\left(\frac{1}{x^{1+2\nu}}\right).
        \label{prooftronc2}
    \end{equation}
    
    \noindent $\dfrac{\Gamma(\nu+1)}{\Gamma(1-\nu)}$(\ref{prooftronc2})$-$(\ref{prooftronc1})$=$(\ref{troncated(1/2,1)3.20}), and we keep $\mathcal{O}\left(\dfrac{1}{x^2}\right)$ instead of $\mathcal{O}\left(\dfrac{1}{x^{1+2\nu}}\right)$ because $\nu\in\left[\frac12,1\right)$. 
    \qed
    \vspace{0.3cm}
    
    \noindent As soon as we have proved this asymptotic expansion, it is not difficult to figure out that one can follow exactly the proof of  Proposition \ref{propAsymptoticbehaviors(0,1/2)} and obtain the following result.
    
    \begin{proposition} Let $\nu\in\left[\frac12,1\right)$,
    \begin{flalign}
        \label{sqrtEven(1/2,1)}
    \sqrt{\lambda_{2(n-1)}}\underset{n\to+\infty}{=}&+\pi\left(n-\frac{\nu}{2}-\frac14\right)-\dfrac{4\nu^2-1}{8\pi n}-\frac{\Gamma(\nu+1)}{\Gamma(1-\nu)}\sin(\nu\pi)\left(\frac{2}{\pi n}\right)^{2\nu}
    +\mathcal{O}\left(\frac{1}{n^2}\right),&&
    \end{flalign}
    \begin{flalign}
        \label{sqrtOdd(1/2,1)}
    \sqrt{\lambda_{2(n-1)+1}}\underset{n\to+\infty}{=}&+\pi\left(n-\frac{\nu}{2}-\frac14\right)-\dfrac{4\nu^2-1}{8\pi n}-\frac{1-2\nu}{1+2\nu}\frac{\Gamma(\nu+1)}{\Gamma(1-\nu)}\sin(\nu\pi)\left(\frac{2}{\pi n}\right)^{2\nu}
    + \mathcal{O}\left(\frac{1}{n^2}\right).
    &&
    \end{flalign}
    \label{propAsymptoticbehaviors(1/2,1)}
    \end{proposition}
    
\noindent Taking the square of (\ref{sqrtEven(1/2,1)}) and (\ref{sqrtOdd(1/2,1)}) leads to the following result.

    \begin{corollary} Let $\nu\in\left[\frac12,1\right)$,
        \begin{flalign}
            \label{even(1/2,1)}
            \lambda_{2(n-1)}\underset{n\to+\infty}{=}&\pi^2\left(n-\frac{\nu}{2}-\frac14\right)^2-2\dfrac{4\nu^2-1}{8}-4\frac{\Gamma(\nu+1)}{\Gamma(1-\nu)}\sin(\nu\pi)\left(\frac{2}{\pi n}\right)^{2\nu-1}+ \mathcal{O}\left(\frac{1}{n}\right),&&
    \end{flalign}
    \begin{flalign}
                \label{odd(1/2,1)}
    \lambda_{2(n-1)+1}\underset{n\to+\infty}{=}&\pi^2\left(n-\frac{\nu}{2}-\frac14\right)^2-2\dfrac{4\nu^2-1}{8}-4\frac{1-2\nu}{1+2\nu}\frac{\Gamma(\nu+1)}{\Gamma(1-\nu)}\sin(\nu\pi)\left(\frac{2}{\pi n}\right)^{2\nu-1}
        + \mathcal{O}\left(\frac{1}{n}\right),&&
    \end{flalign}
    \begin{flalign}
            \label{OddMinusEven(1/2,1)}
    \lambda_{2(n-1)+1}-\lambda_{2(n-1)}\underset{n\to+\infty}{=}\frac{16\nu}{1+2\nu}\frac{\Gamma(\nu+1)}{\Gamma(1-\nu)}\sin(\nu\pi)\left(\frac{2}{\pi n}\right)^{2\nu-1} + \mathcal{O}\left(\frac{1}{n}\right).&&
    \end{flalign}
    \label{corollary1(1/2,1)}
    \end{corollary}

We can also deduce the following less precise result, true for any $\nu\in (0,1)$.

\begin{corollary} 
    \[
        \sqrt{\lambda_{2(n-1)}}\underset{n\to+\infty}{=}\pi \left(n-\frac\nu2-\frac14\right)+\mathcal{O}\left(\frac{1}{n^{\min(1,2\nu)}}\right)
    \]
    and 
        \[
        \sqrt{\lambda_{2(n-1)+1}}\underset{n\to+\infty}{=}\pi \left(n-\frac\nu2-\frac14\right)+\mathcal{O}\left(\frac{1}{n^{\min(1,2\nu)}}\right).
    \]
    \label{corrAsymptSqrtEven(0,1)}
\end{corollary}
\noindent \textbf{\textit{Proof of Corollary \ref{corrAsymptSqrtEven(0,1)}}}. This is directly obtained from (\ref{sqrtEven(0,1/2)}), \eqref{sqrtOdd(0,1/2)},  (\ref{sqrtEven(1/2,1)}) and \eqref{sqrtOdd(1/2,1)}, with a case distinction on $\nu$.
\qed

\section{Null Controllability}
\label{s:nc}
\subsection{Internal Control}

To prove the null controllability property stated in Theorem \ref{theoremNullControlIntern}, we will follow the moment method, that was introduced in \cite{FR71} for the study of the boundary null controllability of the 1D heat equation,  in the spirit of the work \cite{Lag}, which explains how to treat the case of an internal control. The moment method is presented in \cite[Section 5.3.3]{Trelat2024} for instance.
To be able to find a control $u \in L^2\left((0,T)\times\omega\right)$ that is well-defined and that brings the final state to 0, we will have to
\begin{itemize}
    \item find a biorthogonal family $(q_n)_{n\in\N}$ to the family of exponentials $(t\mapsto e^{-\lambda_nt})_{n\in\N}$ in $L^2(0,T)$, \textit{i.e.} verifying 
\[
\int_0^T q_k(t) e^{-\lambda_jt}dt=\delta_{k,j}, \quad \forall k,j\in\N,
\]
with $\delta_{k,j}$ the Kronecker symbol; 
    \item obtain an estimate on the $L^2$-norm of the  $(q_n)_{n\in\N}$;
    \item prove that $\inf_{n\in\N}\int_\omega \phi_n^2 > 0,
        $ with $(\phi_n)_{n\in\N}$ the Hilbert basis of eigenfunctions defined in Theorem \ref{propEigenBasis}.
\end{itemize}
The first point is easily obtained. Moreover, in order to find the estimate on the $L^2$-norm, we will need to introduce the notion of condensation index, following the work of \cite[Section 3]{AmmarKhodja2013}. Then, using results from \cite[Section 4]{AmmarKhodja2013}, we will be able to obtain our second point. These two steps are stated in Proposition \ref{propBiorthogonal}. Finally, to obtain the third point, we will need several asymptotic behaviors that are stated in Lemma \ref{lemmaAsymptotics -nu(0,1)}. Then, we prove it on a control domain of the form $\omega=(\alpha, \beta)\subset(-1,1)$ in Proposition \ref{propInt_a^b}. Afterwards, we generalize the result to the measurable case in Corollary \ref{corIntBorel}. With all those intermediate results, we will finally prove Theorem \ref{theoremNullControlIntern}.

\begin{definition}{\cite[Definition 3.1 and Remark 3.10]{AmmarKhodja2013}} 

\noindent Let $\Lambda = (\mu_n)_{n\in\N}\in\R^\N$ be a positive and increasing sequence such that $\sum_{n\in\N}\dfrac{1}{\mu_n}<+\infty$. The condensation index of the sequence $\Lambda$ is defined by
\begin{equation}
c(\Lambda):= \underset{n\to+\infty}{\limsup}-\frac{1}{\mu_n}\ln{\abs{C'(\mu_n)}}\in\left[0,+\infty\right],
    \label{condeIndex}
\end{equation}
with
\begin{equation}
C(\mu)=\prod_{n=0}^{+\infty}\left(1-\frac{\mu^2}{\mu_n^2}\right), \quad \forall \mu>0.
    \label{C(lambda)}
\end{equation}
\label{defCondensationIndex}
\end{definition}

\begin{remark}\label{rem:cond}In fact, to define the condensation index, we only need to assume that $(\lambda_n)_{n\in\N}$ is positive for large enough $n$. Indeed, one easily proves that removing a finite number of terms in the sequence $\Lambda$ does not change the value of $c(\Lambda)$, in the sense that for any $N\in\mathbb N$, we have 
\begin{equation*}
c(\Lambda)= \underset{n\to+\infty}{\limsup}-\frac{1}{\mu_n}\ln{\abs{C_N'(\mu_n)}}\in\left[0,+\infty\right],
\end{equation*}
with
\begin{equation*}
C_N(\mu)=\prod_{n=N}^{+\infty}\left(1-\frac{\mu^2}{\mu_n^2}\right), \quad \forall \mu>0.
\end{equation*}
\end{remark}
\begin{remark}
We study the condensation index since it provides a measure of the separation of the elements $\lambda_n$ of the sequence $\Lambda$. It will be useful to obtain an essential estimate on the $L^2$-norm of a biorthogonal family in Proposition \ref{propBiorthogonal}. 
\label{remarkDefCondIndex}
\end{remark}

\begin{proposition}
We introduce $\Lambda=(\lambda_n)_{n\in\N}$, with $\{\lambda_n,n\in\N\}$ the eigenvalues of $A_\nu$ defined in Propositions \ref{PropKernel} and \ref{thDistribution}. Then, the condensation index $c(\Lambda)$ is $0$.
\label{propNullCondensIndex}
\end{proposition}

\noindent \textbf{\textit{Proof of Proposition \ref{propNullCondensIndex}}}. This proof takes some inspiration from the proof of \cite[Lemma 4.1]{OU}. First, from (\ref{order}), Corollary \ref{corollary1(0,1/2)}, and Corollary \ref{corollary1(1/2,1)}, we obtain that $(\lambda_n)_{n\in\N\setminus{\{0\}}}$ is a positive increasing sequence and  $\sum_{n\in\N\setminus{\{0\}}}\frac{1}{\lambda_n}<+\infty$. Therefore, the notion of condensation index is well-defined. By examining (\ref{C(lambda)}) (where we start the summation at $n=1$, see Remark \ref{rem:cond}), we notice that
\[
\abs{C'(\lambda_n)}=\frac{2}{\lambda_n}\prod_{\substack{j=1 \\ j \neq n}}^{+\infty}\abs{1-\frac{\lambda_n^2}{\lambda_j^2}},\quad\forall n\in\N\setminus{\{0\}}.
\]
Therefore, we get, for any $n\in\N\setminus\{0\}$,
\[\begin{split}
\frac{1}{\lambda_n}\ln{\abs{C'(\lambda_n})}&=\frac{1}{\lambda_n}\left(\ln{\frac{2}{\lambda_n}}+\sum_{j=1}^{n-2}\ln\abs{1-\frac{\lambda_n^2}{\lambda_j^2}}+\ln{\abs{1-\frac{\lambda_n^2}{\lambda_{n-1}^2}}}+\ln{\abs{1-\frac{\lambda_n^2}{\lambda_{n+1}^2}}}+\sum_{j=n+2}^{+\infty}\ln\abs{1-\frac{\lambda_n^2}{\lambda_j^2}}\right)\\
&=\frac{1}{\lambda_n}\ln{\frac{2}{\lambda_n}}+
\frac{1}{\lambda_n}\ln\left(\frac{\lambda_n^2}{\lambda_{n-1}^2}-1\right)+
\frac{1}{\lambda_n}\ln\left(1-\frac{\lambda_n^2}{\lambda_{n+1}^2}\right)+F_n+G_n,
\end{split}
\]
with $\displaystyle F_n=\frac{1}{\lambda_n}\sum_{j=1}^{n-2}\ln\left(\frac{\lambda_n^2}{\lambda_j^2}-1\right)$, and $\displaystyle G_n=\frac{1}{\lambda_n}\sum_{j=n+2}^{+\infty}\ln\left(1-\frac{\lambda_n^2}{\lambda_j^2}\right)$.

\noindent Moreover, using Corollary \ref{corollary1(0,1/2)} and Corollary \ref{corollary1(1/2,1)}, we know that the asymptotic growth speed of $(\lambda_n)_{n\in\N}$ is of order $n^2$, so it is also the asymptotic growth speed of $(\lambda_{n+1}+\lambda_n)_{n\in\N}$, and that the asymptotic growth speed of $(\lambda_{2n+2}-\lambda_{2n+1})_{n\in\N}$ is of order $n$. Moreover, (\ref{OddMinusEven(0,1/2)}) and (\ref{OddMinusEven(1/2,1)}) give that the asymptotic decay speed of $(\lambda_{2n+1}-\lambda_{2n})_{n\in\N}$ is of order $\displaystyle\frac{1}{n^{2\nu-1}}$. Therefore, we have
\[\begin{split}
    &\frac{1}{\lambda_n}\ln{\frac{2}{\lambda_n}}
    \underset{n\to+\infty}{\longrightarrow}0,\\
    &\frac{1}{\lambda_n}\ln\left(\frac{\lambda_n^2}{\lambda_{n-1}^2}-1\right)=-\frac{1}{\lambda_n}2\ln(\lambda_{n-1})+\frac{1}{\lambda_n}\ln(\lambda_n+\lambda_{n-1})+
    \frac{1}{\lambda_n}\ln(\lambda_n-\lambda_{n-1})
\underset{n\to+\infty}{\longrightarrow}0,\\
&\frac{1}{\lambda_n}\ln\left(1-\frac{\lambda_n^2}{\lambda_{n+1}^2}\right)=-\frac{1}{\lambda_n}2\ln(\lambda_{n+1})+\frac{1}{\lambda_n}\ln(\lambda_{n+1}+\lambda_n)+
    \frac{1}{\lambda_n}\ln(\lambda_{n+1}-\lambda_n)
\underset{n\to+\infty}{\longrightarrow}0.
\end{split}
\]
Thus, if we prove that $\displaystyle F_n\underset{n\to+\infty}{\longrightarrow}0$ and $\displaystyle G_n\underset{n\to+\infty}{\longrightarrow}0$, then, according to Definition \ref{defCondensationIndex}, we will have proved that $c(\Lambda)=0$.

\vspace{0.2cm}
\noindent \textbf{Proof of $\displaystyle F_n\underset{n\to+\infty}{\longrightarrow}0$}.

First, we notice that
\[
\abs{F_n}\leqslant      
\frac{1}{\lambda_n}\sum_{\substack{j=1\\\sqrt{2}\lambda_j<\lambda_n}}^{n-2}\ln\left(\frac{\lambda_n^2}{\lambda_j^2}-1\right)
+
\frac{1}{\lambda_n}\sum_{\substack{j=1\\\sqrt{2}\lambda_j>\lambda_n}}^{n-2}\ln\left(\frac{\lambda_j^2}{\lambda_n^2-\lambda_j^2}\right)
,\quad \forall n\in\N\setminus\{0\}.
\]
Moreover, for $n\geqslant 2$,
\[
0\leqslant\frac{1}{\lambda_n}\sum_{\substack{j=1\\\sqrt{2}\lambda_j<\lambda_n}}^{n-2}\ln\left(\frac{\lambda_n^2}{\lambda_j^2}-1\right)\leqslant
\frac{1}{\lambda_n}\sum_{\substack{j=1\\\sqrt{2}\lambda_j<\lambda_n}}^{n-2}\ln\left(\frac{\lambda_n^2}{\lambda_1^2}-1\right)\leqslant
\frac{n-2}{\lambda_n}\ln\left(\frac{\lambda_n^2}{\lambda_1^2}-1\right)\underset{n\to+\infty}{\longrightarrow}0,
\]
since the asymptotic growth speed of $(\lambda_n)_{n\in\N}$ is of order  $n^2$. Moreover, for $n\geqslant 2$,
\[0\leqslant
\frac{1}{\lambda_n}\sum_{\substack{j=1\\\sqrt{2}\lambda_j>\lambda_n}}^{n-2}\ln\left(\frac{\lambda_j^2}{\lambda_n^2-\lambda_j^2}\right)\leqslant
\frac{1}{\lambda_n}\sum_{\substack{j=1\\\sqrt{2}\lambda_j>\lambda_n}}^{n-2}\ln\left(\frac{\lambda_n^2}{\lambda_n^2-\lambda_{n-2}^2}\right)\leqslant
\frac{n-2}{\lambda_n}\ln\left(\frac{\lambda_n^2}{\lambda_n^2-\lambda_{n-2}^2}\right)
\underset{n\to+\infty}{\longrightarrow}0,
\]
using the already mentioned information about the asymptotic expansion of the eigenvalues and their differences.

Therefore, we proved that $\displaystyle F_n\underset{n\to+\infty}{\longrightarrow}0$.

\vspace{0.2cm}
\noindent \textbf{Proof of  $\displaystyle G_n\underset{n\to+\infty}{\longrightarrow}0$}. First, we notice that
\[
\abs{G_n}\leqslant
\frac{1}{\lambda_n}\sum_{j=n+2}^{+\infty}\ln\left(\frac{\lambda_j^2}{\lambda_j^2-\lambda_n^2}\right)
=\frac{1}{\lambda_n}\sum_{j=n+2}^{+\infty}\ln\left(1+\frac{\lambda_n^2}{\lambda_j^2-\lambda_n^2}\right),
\quad\forall n\in\N\setminus\{0\}.
\]
Moreover, we recall that for any $x>0$, we have $\ln(1+x)\leqslant x$. Therefore,
\[
\abs{G_n}\leqslant
\frac{1}{\lambda_n}\sum_{j=n+2}^{+\infty}\frac{\lambda_n^2}{\lambda_j^2-\lambda_n^2}
=\sum_{j=n+2}^{+\infty}\frac{\lambda_n}{(\lambda_j+\lambda_n)(\lambda_j-\lambda_n)}
\leqslant
\sum_{j=n+2}^{+\infty}\frac{1}{\lambda_j-\lambda_n}.
\]
We notice that the series on the right-hand side of this inequality can be written as
\[
\sum_{j=n+2}^{+\infty}\frac{1}{\lambda_j-\lambda_n}=
\sum_{k \textrm{ \ such\ that\ } 2(k-1)\geqslant n+2}^{+\infty}\frac{1}{\lambda_{2(k-1)}-\lambda_n}+
\sum_{k \textrm{ \ such\ that\ } 2(k-1)+1\geqslant n+2}^{+\infty}\frac{1}{\lambda_{2(k-1)+1}-\lambda_n}, \quad \forall n
>0.
\]
First, we want to prove that $\displaystyle G_{2(n-1)}\underset{n\to+\infty}{\longrightarrow}0$. We have
\begin{equation}
\abs{G_{2(n-1)}}\leqslant\sum_{k=n+1}^{+\infty}\frac{1}{\lambda_{2(k-1)}-\lambda_{2(n-1)}}+
\sum_{k=n+1}^{+\infty}\frac{1}{\lambda_{2(k-1)+1}-\lambda_{2(n-1)}}
,\quad\forall n>1.
\label{proofIndexG}
\end{equation}
If $\nu\in\left(0,\frac12\right)$, Corollary \ref{corollary1(0,1/2)} gives
\[\begin{split}
    \lambda_{2(n-1)}\underset{n\to+\infty}{=}\pi^2\left(n-\frac{\nu}{2}-\frac14\right)^2+\varepsilon_{2\nu-1}^{even}(n), \quad 
    \lambda_{2(n-1)+1}
    \underset{n\to+\infty}{=}\pi^2\left(n-\frac{\nu}{2}-\frac14\right)^2+\varepsilon_{2\nu-1}^{odd}(n),
\end{split}
\]
with $\displaystyle\varepsilon_{2\nu-1}^{even}(n),\varepsilon_{2\nu-1}^{odd}(n)=\mathcal{O}\left(\frac{1}{n^{2\nu-1}}\right)$. From there, it is not difficult to infer that

\begin{equation}
\begin{split}
    \exists N\in\N, \forall n>N, \forall k>n, \quad&\lambda_{2(k-1)}-\lambda_{2(n-1)}\geqslant (\pi^2-1)\left(k-\frac{\nu}{2}-\frac14\right)^2 - (\pi^2-1)\left(n-\frac{\nu}{2}-\frac14\right)^2,\\
    \textrm{and} \ &\lambda_{2(k-1)+1}-\lambda_{2(n-1)}\geqslant (\pi^2-1)\left(k-\frac{\nu}{2}-\frac14\right)^2 - (\pi^2-1)\left(n-\frac{\nu}{2}-\frac14\right)^2.
\end{split}
\label{proofIndex}
\end{equation}

\noindent If $\nu\in[\frac12,1)$, Corollary \ref{corollary1(1/2,1)} gives
\[\begin{split}
    \lambda_{2(n-1)}\underset{n\to+\infty}{=}\pi^2\left(n-\frac{\nu}{2}-\frac14\right)^2-2\dfrac{4\nu^2-1}{8}+\varepsilon_{2\nu-1}^{even}(n),\\ 
    \lambda_{2(n-1)+1}
    \underset{n\to+\infty}{=}\pi^2\left(n-\frac{\nu}{2}-\frac14\right)^2-2\dfrac{4\nu^2-1}{8}+\varepsilon_{2\nu-1}^{odd}(n),
\end{split}
\]
with $\displaystyle\varepsilon_{2\nu-1}^{even}(n),\varepsilon_{2\nu-1}^{odd}(n)=\mathcal{O}\left(\frac{1}{n^{2\nu-1}}\right)$. Therefore, it is very easy to obtain the same statement (\ref{proofIndex}).

\noindent Let $\nu\in(0,1)$. From (\ref{proofIndexG}) and (\ref{proofIndex}), we obtain that for any $n>N$,
\[\begin{split}
  \abs{G_{2(n-1)}} &\leqslant
  \frac{2}{\pi^2-1}\sum_{k=n+1}^{+\infty}\frac{1}{\left(k-\frac{\nu}{2}-\frac14\right)^2 - \left(n-\frac{\nu}{2}-\frac14\right)^2}
  =\frac{2}{\pi^2-1}\sum_{k=1}^{+\infty}\frac{1}{\left(k-\frac{\nu}{2}-\frac14\right)^2 +2n\left(k-\frac{\nu}{2}-\frac14\right)}\\
  &\leqslant \frac{2}{\pi^2-1}\left(\frac{1}{\left(1-\frac{\nu}{2}-\frac14\right)^2 +2n\left(1-\frac{\nu}{2}-\frac14\right)}+\sum_{k=2}^{+\infty}\int_{k-1}^k\frac{1}{\left(x-\frac{\nu}{2}-\frac14\right)^2 +2n\left(x-\frac{\nu}{2}-\frac14\right)}\textrm{d}x\right)\\
  &=\frac{2}{\pi^2-1}\left(\frac{1}{\left(1-\frac{\nu}{2}-\frac14\right)^2 +2n\left(1-\frac{\nu}{2}-\frac14\right)}+\left[\frac{1}{2n}\ln\left(\frac{x-\frac{\nu}{2}-\frac14}{x-\frac{\nu}{2}-\frac14+2n}\right)\right]^{+\infty}_1\right)
  \underset{n\to+\infty}{\longrightarrow}0.
\end{split}
\]
We can reiterate the proof by following the same steps to obtain that we have as well $\displaystyle G_{2(n-1)+1}\underset{n\to+\infty}{\longrightarrow}0$.
As a consequence, we have that the condensation index $c(\Lambda)$ is 0, which concludes the proof.
\qed

\vspace{0.3cm}

\noindent From the previous Proposition, we easily deduce the existence of biorthogonal families at any time $T>0$.
\begin{proposition}
Let $T>0$. We consider $(t\mapsto e^{-\lambda_nt})_{n\in\N}\in L^2(0,T)^\N$, with $\{\lambda_n,n\in\N\}$ the eigenvalues of $A_\nu$ defined in Proposition \ref{thDistribution} and \ref{PropKernel}. Then, there exists a biorthogonal family $(q_n)_{n\in\N}\in L^2(0,T)^\N$ to $(t\mapsto e^{-\lambda_nt})_{n\in\N}$ such that
\begin{equation}
\forall \varepsilon>0, \ \exists C_\varepsilon(T)>0, \ \forall n\in\N, \quad \norm{q_n}_{L^2(0,T)}\leqslant C_\varepsilon(T) e^{\varepsilon \lambda_n}.
    \label{inequalityBiorthogonal}
\end{equation}
\label{propBiorthogonal}
\end{proposition}

\noindent \textbf{\textit{Proof of Proposition \ref{propBiorthogonal}}}. Let $T>0$.  From (\ref{order}), Corollary \ref{corollary1(0,1/2)}, and Corollary \ref{corollary1(1/2,1)}, we obtain that $(\lambda_n)_{n\in\N}$ is a non-negative increasing sequence and  $\sum_{n\in\N\setminus{\{0\}}}\frac{1}{\lambda_n}<+\infty$. Therefore, we can use \cite[Theorem 4.1]{AmmarKhodja2013} (which can be easily adapted to include the case where $0$ is a simple eigenvalue by an appropriate translation argument in time), which gives the existence of the biorthogonal family $(q_n)_{n\in\N}$. Proposition \ref{propNullCondensIndex} gives that the condensation index is equal to 0. Therefore, \cite[Remark 4.3]{AmmarKhodja2013} gives (\ref{inequalityBiorthogonal}), which concludes the proof.
\qed
\vspace{0.3cm}

\noindent Now, we want to prove the following Proposition.

\begin{proposition} Let $\omega$ be a measurable subset of $(-1,1)$ of positive Lebesgue measure. Then,
    \[
       \inf_{n\in\N}\int_\omega \phi_n^2 > 0,
        \]
    with $(\phi_n)_{n\in\N}$ the Hilbert basis of eigenfunctions defined in Theorem \ref{propEigenBasis}.
    \label{corIntBorel}
\end{proposition}

\begin{remark}
\label{dif:NL}
Let us emphasize here that obtaining this proposition is not straightforward, and relies on the explicit form of our eigenfunctions. This comes from the fact that we are working here with singular potentials and nonstandard self-adjoint extensions, so that many tools from perturbation theory are not available. Indeed, for the usual eigenfunctions of the Dirichlet-Laplace operator on $(0,1)$ given by $\sqrt{2}\sin(n\pi x)$, we have, for any measurable set $\omega\subset (0,1)$ of positive measure, by the Riemann-Lebesgue Lemma,

    \[
     \lim_{n\rightarrow +\infty} \int_\omega 2\sin(n\pi x)^2\mathrm{d}x=  \lim_{n\rightarrow +\infty}\int_\omega (1-\cos(2n\pi x)\mathrm{d}x)=|\omega|,
        \]
                where $|\cdot|$ denotes the Lebesgue measure.
       So, we easily infer that 
    \[
     \inf_{n\in\mathbb N\setminus{\{0\}}} \int_\omega 2\sin(n\pi x)^2\mathrm{d}x>0.
        \]
    For the Dirichlet-Laplace operator with bounded potential, it is easy to prove by a perturbation argument that the same property holds (see \textit{e.g.} \cite[Appendix A]{LLP}). However, this reasoning is ineffective in our case, because our eigenfunctions  cannot be written as adequate perturbations of $\sin$ or $\cos$ functions.
\end{remark}

\noindent In order to prove Proposition \ref{corIntBorel}, we need several intermediate steps. First, we state several asymptotic results that will be useful later on.

\begin{lemma} Let $\beta_0\in(0,1)$.
    \begin{subequations}
    \setlength{\jot}{3pt}
    \renewcommand{\theequation}{\theparentequation.\alph{equation}}
    \begin{align}
    &J_{-\nu}\left(\sqrt{\lambda_{2(n-1)}}\beta\right)\underset{n\to+\infty}{=} \frac{1}{\sqrt{\pi}}\left(\frac{2}{\pi \beta n}\right)^{1/2} \left(\sin\left(\sqrt{\lambda_{2(n-1)}}\beta + \nu\frac{\pi}{2}+\frac{\pi}{4}\right)+\mathcal{O}_{\beta_0}\left(\frac{1}{n}\right)\right),  
    \forall \beta\in[\beta_0,1),\label{asympt J-nu(beta) (0,1)}\\
        &J_{-\nu}\left(\sqrt{\lambda_{2(n-1)}}\right)\underset{n\to+\infty}{=} \mathcal{O}\left(\frac{1}{n^{1/2+\min(1,2\nu)}}\right), \label{asympt J-nu (0,1)}\\
        &{J_{-\nu}}'\left(\sqrt{\lambda_{2(n-1)}}\beta\right)\underset{n\to+\infty}{=} \frac{1}{\sqrt{\pi}}\left(\frac{2}{\pi \beta n}\right)^{1/2} \left(\cos\left(\sqrt{\lambda_{2(n-1)}}\beta + \nu\frac{\pi}{2}+\frac{\pi}{4}\right)+\mathcal{O}_{\beta_0}\left(\frac{1}{n}\right)\right), \forall 1>\beta \geqslant\beta_0, \label{asympt J-nu'(beta) (0,1)}\\
        &{J_{-\nu}}'\left(\sqrt{\lambda_{2(n-1)}}\right)\underset{n\to+\infty}{=} \frac{(-1)^n}{\sqrt{\pi}}\left(\frac{2}{\pi n}\right)^{1/2}\left(1+\mathcal{O}\left(\frac{1}{n^{\min(1,4\nu)}}\right)\right), \label{asympt J-nu' (0,1)}\\
        &a_{2(n-1)}^2\underset{n\to+\infty}{=}\frac{1}{\pi}\frac{2}{\pi n}\left(1+\mathcal{O}\left(\frac{1}{n^{\min(1,2\nu)}}\right)\right), \label{asympt a2n-1 (0,1)}        
    \end{align}
    \end{subequations}
where $a_{2(n-1)}$ is defined in \eqref{a2n}, and $\mathcal{O}_{\beta_0}\left(\ \cdot \ \right)$ means: there exists $C(\beta_0)>0$ and $n_0\in\mathbb N\setminus{\{0\}}$, which might depend  on $\beta_0$, such that for any $n\geqslant n_0$ and any $\beta\geqslant \beta_0$,
$$\abs{\mathcal{O}_{\beta_0}\left(\ \cdot\ \right)}\leqslant C(\beta_0)\abs{ \ \cdot \ }. $$
    \label{lemmaAsymptotics -nu(0,1)}
\end{lemma}

\noindent \textbf{\textit{Proof of Lemma \ref{lemmaAsymptotics -nu(0,1)}}}. Let $\beta_0\in(0,1)$.\\
\noindent \textbf{Proof of (\ref{asympt J-nu(beta) (0,1)}) and (\ref{asympt J-nu (0,1)})}. From Lemma \ref{lemmaExpansion}, we obtain
\[
J_{-\nu}(x)\underset{x\to+\infty}{=}\sqrt{\frac{2}{\pi x}}\left(\cos\left(x+\nu\frac{\pi}{2}-\frac\pi4\right)+\mathcal{O}\left(\frac1x\right)\right).
\]
Let $\beta\geqslant \beta_0$. We know that  $\sqrt{\lambda_{2(n-1)}}\underset{n\to+\infty}{\longrightarrow}+\infty$. Therefore,
\[
\begin{split}
J_{-\nu}\left(\sqrt{\lambda_{2(n-1)}}\beta\right)&{=}\sqrt{\frac{2}{\pi \sqrt{\lambda_{2(n-1)}}\beta}}\left(\cos\left(\sqrt{\lambda_{2(n-1)}}\beta+\nu\frac{\pi}{2}-\frac\pi4\right)+\mathcal{O}_{\beta_0}\left(\frac{1}{\sqrt{\lambda_{2(n-1)}}}\right)\right)\\
&{=}\frac{1}{\sqrt{\pi}}\left(\frac{2}{\pi\beta n}\right)^{1/2}\left(1+\mathcal{O}_{\beta_0}\left(\frac1n\right)\right)^{-1/2}\left(\sin\left(\sqrt{\lambda_{2(n-1)}}\beta+\nu\frac{\pi}{2}+\frac\pi4\right)+\mathcal{O}_{\beta_0}\left(\frac1n\right)\right),
\end{split}
\]
where we used Corollary \ref{corrAsymptSqrtEven(0,1)} for the asymptotic expansion of $\sqrt{\lambda_{2(n-1)}}$. This concludes the proof of (\ref{asympt J-nu(beta) (0,1)}). Now, let us take $\beta=1$. Using the previous computations, we obtain
\[
\begin{split}
J_{-\nu}\left(\sqrt{\lambda_{2(n-1)}}\right)
&\underset{n\to+\infty}{=}\mathcal{O}\left(\frac{1}{n^{1/2}}\right)\left[\sin\left(\pi\left(n-\frac\nu2-\frac14\right)+\mathcal{O}\left(\frac{1}{n^{\min(1,2\nu)}}\right)+\nu\frac{\pi}{2}+\frac\pi4\right)+\mathcal{O}\left(\frac{1}{n}\right)\right]\\
&\underset{n\to+\infty}{=}\mathcal{O}\left(\frac{1}{n^{1/2}}\right)\left((-1)^n\mathcal{O}\left(\frac{1}{n^{\min(1,2\nu)}}\right)+\mathcal{O}\left(\frac{1}{n}\right)\right)\underset{n\to+\infty}{=}\mathcal{O}\left(\frac{1}{n^{1/2+\min(1,2\nu)}}\right),
\end{split}
\]
which concludes the proof of (\ref{asympt J-nu (0,1)}).
\vspace{0.2cm}

\noindent \textbf{Proof of (\ref{asympt J-nu'(beta) (0,1)}) and  (\ref{asympt J-nu' (0,1)})}. Lemma \ref{lemmaDerivateBessel} states that ${J_{-\nu}}'=\frac{1}{2}\left(J_{-\nu-1}-J_{1-\nu}\right)$. Therefore, by Lemma \ref{lemmaExpansion} on $J_{-\nu-1}$ and $J_{1-\nu}$, we obtain that
\[
\begin{split}
{J_{-\nu}}'(x)
&\underset{x\to+\infty}{=}\frac12\sqrt{\frac{2}{\pi x}}\left(\cos\left(x-(-\nu-1)\frac{\pi}{2}-\frac\pi4\right)-\cos\left(x-(1-\nu)\frac{\pi}{2}-\frac\pi4\right)+\mathcal{O}\left(\frac1x\right)\right)\\
&\underset{x\to+\infty}{=}\frac12\sqrt{\frac{2}{\pi x}}\left(\cos\left(x+\nu\frac\pi2+\frac\pi4\right)-\cos\left(x+\nu\frac\pi2-3\frac\pi4\right)+\mathcal{O}\left(\frac1x\right)\right)\\
&\underset{x\to+\infty}{=}\sqrt{\frac{2}{\pi x}}\left(\cos\left(x+\nu\frac\pi2+\frac\pi4\right)+\mathcal{O}\left(\frac1x\right)\right).
\end{split}
\]
Let $\beta\geqslant \beta_0$. Since $\sqrt{\lambda_{2(n-1)}}\underset{n\to+\infty}{\longrightarrow}+\infty$, we have
\[
\begin{split}
{J_{-\nu}}'\left(\sqrt{\lambda_{2(n-1)}}\beta\right)&\underset{n\to+\infty}{=}\sqrt{\frac{2}{\pi \sqrt{\lambda_{2(n-1)}}\beta}}\left(\cos\left(\sqrt{\lambda_{2(n-1)}}\beta+\nu\frac{\pi}{2}+\frac\pi4\right)+\mathcal{O}\left(\frac{1}{\sqrt{\lambda_{2(n-1)}}\beta}\right)\right)\\
&\underset{n\to+\infty}{=}\frac{1}{\sqrt{\pi}}\left(\frac{2}{\pi\beta n}\right)^{1/2}\left(1+\mathcal{O}_{\beta_0}\left(\frac1n\right)\right)^{-1/2}\left(\cos\left(\sqrt{\lambda_{2(n-1)}}\beta+\nu\frac{\pi}{2}+\frac\pi4\right)+\mathcal{O}_{\beta_0}\left(\frac1n\right)\right)
\end{split}
\]
which concludes the proof of (\ref{asympt J-nu'(beta) (0,1)}). Now, we fix $\beta=1$. Substituting $\sqrt{\lambda_{2(n-1)}}$ in (\ref{asympt J-nu'(beta) (0,1)}) with Corollary \ref{corrAsymptSqrtEven(0,1)} easily gives (\ref{asympt J-nu' (0,1)}).
\vspace{0.2cm}

\noindent \textbf{Proof of (\ref{asympt a2n-1 (0,1)})}. Let $n>1$. We recall that the expression of $a_{2(n-1)}$ is given in (\ref{a2n}).\\
Since $\sqrt{\lambda_{2(n-1)}}\underset{n\to+\infty}{=}\mathcal{O}(n)$, we obtain from Corollary \ref{corrAsymptBesselAndDerivatives} that 
\[
\begin{split}
    a_{2(n-1)}^2&=\left(1+\mathcal{O}\left(\frac{1}{n^{2}}\right)\right)J_{-\nu}\left(\sqrt{\lambda_{2(n-1)}}\right)^2+{J_{-\nu}}'\left(\sqrt{\lambda_{2(n-1)}}\right)^2\\
        &+\mathcal{O}\left(\frac{1}{n^{2\nu}}\right)\left[\left(1+\mathcal{O}\left(\frac{1}{n^{2}}\right)\right)\mathcal{O}\left(\frac{1}{n^{1/2}}\right)\mathcal{O}\left(\frac{1}{n^{1/2}}\right) + \mathcal{O}\left(\frac{1}{n^{1/2}}\right)\mathcal{O}\left(\frac{1}{n^{1/2}}\right)
        +\mathcal{O}\left(\frac{1}{n^{2}}\right)\right]\\
        &+\mathcal{O}\left(\frac{1}{n^{4\nu}}\right)\left[\left(1+\mathcal{O}\left(\frac{1}{n^{2}}\right)\right)\mathcal{O}\left(\frac{1}{n}\right)+\mathcal{O}\left(\frac{1}{n}\right)\right]\\
        &=\left(1+\mathcal{O}\left(\frac{1}{n^{2}}\right)\right)J_{-\nu}\left(\sqrt{\lambda_{2(n-1)}}\right)^2+{J_{-\nu}}'\left(\sqrt{\lambda_{2(n-1)}}\right)^2+\mathcal{O}\left(\frac{1}{n^{1+2\nu}}\right).
        \end{split}
\]
Using (\ref{asympt J-nu (0,1)}) and (\ref{asympt J-nu' (0,1)}), we obtain
\[
\begin{split}
    a_{2(n-1)}^2
    &\underset{n\to+\infty}{=}\left(1+\mathcal{O}\left(\frac{1}{n^{2}}\right)\right)\mathcal{O}\left(\frac{1}{n^{1+\min(2,4\nu)}}\right)+
   \frac{2}{\pi^2 n}\left(1+\mathcal{O}\left(\frac{1}{n^{\min(1,4\nu)}}\right)\right)+\mathcal{O}\left(\frac{1}{n^{1+2\nu}}\right)\\
    &\underset{n\to+\infty}{=}
  \frac{2}{\pi^2 n}\left(1+\mathcal{O}\left(\frac{1}{n^{\min(1,4\nu)}}\right)+\mathcal{O}\left(\frac{1}{n^{\min(2,4\nu)}}\right)+\mathcal{O}\left(\frac{1}{n^{2\nu}}\right)\right), \quad \mathrm{which\ gives\ (\ref{asympt a2n-1 (0,1)})}.
        \end{split}
\]
\qed

\begin{remark}
The reader can convince themselves
that all the asymptotic behaviors that we obtained so far in this section are still valid when the index $ 2(n-1)$ is replaced by $ 2(n-1)+1$ and $ \dfrac{\Gamma(\nu+1)}{\Gamma(1-\nu)}$ is replaced by  $\dfrac{1-2\nu}{1+2\nu}\dfrac{\Gamma(\nu+1)}{\Gamma(1-\nu)}$. In particular, Lemma \ref{lemmaAsymptotics -nu(0,1)} also holds for even indexes $ 2(n-1)+1$.
\label{remarkSameOdd}
\end{remark}

\noindent The previous asymptotic result is useful to prove the following non-concentration property of our eigenfunctions on any non-trivial interval that is far from $0$.
\begin{proposition}
Let $\alpha_0\in(0,1)$. There exists $C>0$ and $n_0\in\mathbb N\setminus{\{0\}}$ (which might depend on $\alpha_0$ and $\nu$) such that 
    \[
       \int_\alpha^\beta \phi_n^2\geqslant C (\beta-\alpha), \quad \alpha_0\leqslant \alpha<\beta<1,\quad \forall n\geqslant n_0.
        \]
Let $\beta_0\in(0,1)$. There exists $C>0$ and $n_0\in\mathbb N\setminus{\{0\}}$ (which might depend on $\beta_0$ and $\nu$) such that 

        \[
       \int_\alpha^\beta \phi_n^2\geqslant C (\beta-\alpha), \quad -1< \alpha<\beta\leqslant -\beta_0,\quad\forall n\geqslant n_0,
        \]    
with $(\phi_n)_{n\in\N}$ the Hilbert basis of eigenfunctions  in $L^2(-1,1)$ defined in Theorem \ref{propEigenBasis}.
    \label{propInt_a^b}
\end{proposition}

\noindent \textbf{\textit{Proof of Proposition \ref{propInt_a^b}}}.
Since $\phi_{2(n-1)}^2$ and $\phi_{2(n-1)+1}^2$ are even, there is no loss of generality by assuming that we are in the case $\alpha_0\leqslant \alpha<\beta<1$.
First, let us start with the even indexes $2(n-1)$. We have, for all $n>1$,
\begin{equation}
\label{ncab}
\begin{split}
    \int_\alpha^\beta \phi_{2(n-1)}^2=&\int_\alpha^\beta xJ_{-\nu}\left(\sqrt{\lambda_{2(n-1)}}x\right)^2\mathrm{d}x\\
    &-2\frac{\Gamma(\nu+1)}{\Gamma(1-\nu)}\left(\frac{\sqrt{\lambda_{2(n-1)}}}{2}\right)^{-2\nu}\int_\alpha^\beta xJ_{\nu}\left(\sqrt{\lambda_{2(n-1)}}x\right)J_{-\nu}\left(\sqrt{\lambda_{2(n-1)}}x\right)\mathrm{d}x \\
    &+\left(\frac{\Gamma(\nu+1)}{\Gamma(1-\nu)}\left(\frac{\sqrt{\lambda_{2(n-1)}}}{2}\right)^{-2\nu}\right)^2\int_\alpha^\beta xJ_\nu\left(\sqrt{\lambda_{2(n-1)}}x\right)^2\mathrm{d}x.
    \end{split}
\end{equation}
From the mean value theorem for integrals, there exists $r_n\in (\alpha,\beta)$ such that 
$$\int_\alpha^\beta xJ_{\nu}\left(\sqrt{\lambda_{2(n-1)}}x\right)J_{-\nu}\left(\sqrt{\lambda_{2(n-1)}}x\right)\mathrm{d}x=(\beta-\alpha) r_nJ_{\nu}\left(\sqrt{\lambda_{2(n-1)}}r_n\right)J_{-\nu}\left(\sqrt{\lambda_{2(n-1)}}r_n\right).$$
 From Corollary \ref{corrAsymptBesselAndDerivatives}, we deduce that there exists $C>0$ such that for any $n$ large enough, 
 \begin{equation*}\left |r_nJ_{\nu}\left(\sqrt{\lambda_{2(n-1)}}r_n\right)J_{-\nu}\left(\sqrt{\lambda_{2(n-1)}}r_n\right)\right|\leqslant \frac{C}{\lambda_{2(n-1)}}.\end{equation*}
Using Corollary \ref{corrAsymptSqrtEven(0,1)}, we deduce that for some new $C>0$, 
 \begin{equation}\label{ncab1}\left |-2\frac{\Gamma(\nu+1)}{\Gamma(1-\nu)}\left(\frac{\sqrt{\lambda_{2(n-1)}}}{2}\right)^{-2\nu}\int_\alpha^\beta xJ_{\nu}\left(\sqrt{\lambda_{2(n-1)}}x\right)J_{-\nu}\left(\sqrt{\lambda_{2(n-1)}}x\right)\mathrm{d}x\right |\leqslant \frac{(\beta-\alpha)C}{n^{1+2\nu}}.\end{equation}
 A similar reasoning also ensures that for $n$ large enough,
  \begin{equation}\label{ncab2}\left |\left(\frac{\Gamma(\nu+1)}{\Gamma(1-\nu)}\left(\frac{\sqrt{\lambda_{2(n-1)}}}{2}\right)^{-2\nu}\right)^2\int_\alpha^\beta xJ_\nu\left(\sqrt{\lambda_{2(n-1)}}x\right)^2\mathrm{d}x\right |\leqslant \frac{(\beta-\alpha)C}{n^{1+4\nu}}.\end{equation}
 Now, let us go back to the first term of \eqref{ncab}. Using (\ref{intSame-NuSameA - alpha-beta}), with (\ref{asympt J-nu(beta) (0,1)}) and (\ref{asympt J-nu'(beta) (0,1)}), we obtain that
 \[
\begin{split}
    \int_\alpha^\beta xJ_{-\nu}\left(\sqrt{\lambda_{2(n-1)}}x\right)^2
        &=\frac12 \beta^2\left(\left(1-\frac{\nu^2}{\sqrt{\lambda_{2(n-1)}}^2\beta^2}\right)J_{-\nu}(\sqrt{\lambda_{2(n-1)}}\beta)^2+{J_{-\nu}}'(\sqrt{\lambda_{2(n-1)}}\beta)^2\right) \\
        &-\frac12 \alpha^2\left(\left(1-\frac{\nu^2}{\sqrt{\lambda_{2(n-1)}}^2\alpha^2}\right)J_{-\nu}(\sqrt{\lambda_{2(n-1)}}\alpha)^2+{J_{-\nu}}'(\sqrt{\lambda_{2(n-1)}}\alpha)^2\right)\\
        &=\frac12 \beta^2\left[\left(1+{\mathcal{O}_{\alpha_0}\left(\frac{1}{ n^2}\right)}\right)\frac{2}{\pi^2 \beta n} \left(\sin^2\left(\beta\sqrt{\lambda_{2(n-1)}} + \nu\frac{\pi}{2}+\frac{\pi}{4}\right)+\mathcal{O}_{\alpha_0}\left(\frac{1}{n}\right)\right)\right.\\
        &\left.+\frac{2}{\pi^2 \beta n} \left(\cos^2\left(\beta\sqrt{\lambda_{2(n-1)}} + \nu\frac{\pi}{2}+\frac{\pi}{4}\right)+\mathcal{O}_{\alpha_0}\left(\frac{1}{n}\right)\right)\right] \\
        &-\frac12 \alpha^2\left[\left(1+{\mathcal{O}_{\alpha_0}\left(\frac{1}{n^2}\right)}\right)\frac{2}{\pi^2 \alpha n} \left(\sin^2\left(\alpha\sqrt{\lambda_{2(n-1)}} + \nu\frac{\pi}{2}+\frac{\pi}{4}\right)+\mathcal{O}_{\alpha_0}\left(\frac{1}{n}\right)\right)\right.\\
        &\left.+\frac{2}{\pi ^2\alpha n} \left(\cos^2\left(\alpha\sqrt{\lambda_{2(n-1)}} + \nu\frac{\pi}{2}+\frac{\pi}{4}\right)+\mathcal{O}_{\alpha_0}\left(\frac{1}{n}\right)\right)\right].
        \end{split}
\]
Since $\cos^2+\sin^2=1$, we have 
\[
\begin{split}
    \int_\alpha^\beta xJ_{-\nu}\left(\sqrt{\lambda_{2(n-1)}}x\right)^2\mathrm{d}x
        =\frac{1}{\pi^2 n}(\beta-\alpha)\left(1+\mathcal{O}_{\alpha_0}\left(\frac{1}{n}\right)\right).
        \end{split}
\]
Hence, we deduce that there exists $C>0$, which might depend on $\nu$ and $\alpha_0$, and that might from now on change from line to line, and there exists $\tilde n_0\in\mathbb N\setminus{\{0\}}$ (depending also on $\nu$ and $\alpha_0$) such that for any $n\geqslant \tilde n_0$,
   $$ \int_\alpha^\beta xJ_{-\nu}\left(\sqrt{\lambda_{2(n-1)}}x\right)^2\mathrm{d}x\geqslant \frac{C}{n}(\beta-\alpha).$$
   Coming back to \eqref{ncab} and using \eqref{ncab1} and \eqref{ncab2}, we deduce that there exists some $n_0\geqslant \tilde{n}_0$ such that for $n\geqslant n_0$, we have
\[
\begin{split}
    \int_\alpha^\beta \psi_{2(n-1)}^2
       \geqslant \frac{C}{n}(\beta-\alpha).
        \end{split}
\]

Therefore, using (\ref{asympt a2n-1 (0,1)}), we obtain that for some (possibly larger) $n_0$, we have 
\[
\begin{split}
    \int_\alpha^\beta \phi_{2(n-1)}^2
     \geqslant C(\beta-\alpha),\quad \forall n\geqslant n_0.
        \end{split}
\]
This concludes the proof for even indexes of eigenfunctions and for $\alpha_0\leqslant\alpha,\beta<1$.

\vspace*{0.2cm}

\noindent Then, we observe that for $n>0$, the expression of $\phi_{2(n-1)+1}^2$ is very similar to the one of $\phi_{2(n-1)}^2$ on $(0,1)$, except that the index $ 2(n-1)$ is replaced by $ 2(n-1)+1$ and that $ \dfrac{\Gamma(\nu+1)}{\Gamma(1-\nu)}$ is replaced by $\dfrac{1-2\nu}{1+2\nu}\dfrac{\Gamma(\nu+1)}{\Gamma(1-\nu)}$. Therefore, using Remark \ref{remarkSameOdd}, we also have (for some possibly larger $n_0$)
\[
\begin{split}
    \int_\alpha^\beta \phi_{2(n-1)+1}^2\geqslant 
     C(\beta-\alpha),\quad \forall n\geqslant n_0.
        \end{split}
\]

\qed
\vspace{0.3cm}

\noindent By using the structure of the open sets of $\mathbb R$, we can deduce the following result.

\begin{corollary} Let $\alpha_0\in (0,1)$. 
There exists $C>0$ and $n_0\in\mathbb N\setminus{\{0\}}$ (which might depend on $\alpha_0$ and $\nu$) such that for any open subset $\mathcal O$ that is included in $(\alpha_0,1)$ or $(-1,-\alpha_0)$, we have

    \[
       \int_{\mathcal O} \phi_n^2\geqslant C | \mathcal O|,\quad,\forall n\geqslant n_0.
        \]

    \label{corouvert}
\end{corollary}

\noindent\textbf{\textit{Proof of Corollary \ref{corouvert}.}}
As already explained, it suffices to treat the case where $\mathcal O \subset (\alpha_0,1)$. Then, $\mathcal O$ is a disjoint union of a countable collection of open intervals $\left\{(a_k, b_k)\right\}_{k\in\mathcal{I}}$ with $\mathcal{I}$ finite or countable, and for any $k\in I$, $\alpha_0<a_k<b_k<1$. With these notations, we have
\[
\int_\mathcal{O} \phi_n^2= \sum_{k\in\mathcal{I}} \int_{a_k}^{b_k} \phi_n^2.
\]
Using Lemma \ref{propInt_a^b}, this gives, for some $C>0$ and $n_0\in\mathbb N\setminus{\{0\}}$ which might depend on $\alpha_0$ and $\nu$, 
\[
\int_\mathcal{O} \phi_n^2\geqslant C  \sum_{k\in\mathcal{I}} (b_k-a_k)=C|\mathcal O|,\quad\forall n\geqslant n_0,
\]
whence the result.

\qed

\vspace*{0.3cm}

\noindent The last step is to pass from an open set to a measurable set of positive measure.

\vspace*{0.3cm}

\noindent \textbf{\textit{Proof of Proposition \ref{corIntBorel}}} Now, we consider $\omega$ a measurable subset of $(-1,1)$ assumed to have positive Lebesgue measure. Let $\tilde\omega$ be another measurable subset of $(-1,1)$. We will use repeatedly the following property: 
$$\tilde \omega\subset\omega \Rightarrow     
       \inf_{n\in\N}\int_\omega \phi_n^2 \geqslant     
       \inf_{n\in\N}\int_{\tilde \omega} \phi_n^2
        .$$

\color{black}
        \noindent Hence, to prove the desired result, it is enough to prove it for a measurable subset of $\omega$.
        
        \vspace{0.3cm}

\noindent Assume without loss of generality that $|\omega \cap[0,1)|>0$ (otherwise, we use a symmetry argument as before). Therefore, we can change $\omega$  into $\omega\cap [0,1)$, according to the discussion above, and we still call it $\omega$. By the Lebesgue density theorem, for almost every $x\in \omega$ and for $r_x>0$ small enough, we have $|\omega\cap B_x(r_x)|>0$. Since $\{0\}$ is of null measure, and $\omega\subset [0,1)$, 
there exists some $x_0 \in \omega$ with $x_0>0$ and some $r>0$ such that $|\omega\cap B_{x_0}(r)|>0$. Reducing $r$ if necessary, we can assume that $r< x_0/2$.
We set $2\alpha_0=x_0/2$. Then $\omega\cap B_{x_0}(r)\subset (2\alpha_0,1)$. Hence, we can change $\omega$ into $\omega\cap B_{x_0}(r)$, so we reduce to the case  where $\omega\subset (2\alpha_0,1)$ with $\alpha_0>0$.

\vspace{0.3cm}

 \noindent By exterior regularity of the Lebesgue measure \cite[Theorem 2.20]{Rudin1987}, 
 for all $k\in \mathbb N\setminus{\{0\}}$, there exists an open set $\mathcal{O}_k$ of $(-1,1)$ with $\omega\subseteq \mathcal{O}_k$ such that $\abs{\mathcal{O}_k \setminus \omega}<1/k$.
 Moreover, reducing $\mathcal {O}_k$ if necessary, we can assume that $\mathcal{O}_k \subset (\alpha_0,1)$ for any $k\in\mathbb N\setminus{\{0\}}$. By the ``reverse'' of the dominated convergence theorem, we know that for some extraction $\varphi$, 
  $$\mathbb {1}_{\mathcal{O}_{\varphi(k)}\setminus \omega }\underset{k\to+\infty}{\longrightarrow} 0, \quad \mbox{ a.e. on } (-1,1).$$

\noindent Since $\phi_n^2\in L^1(-1,1)$ for any $n\in\mathbb N\setminus{\{0\}}$, and 
  $$|\mathbb {1}_{\mathcal{O}_{\varphi(k)}} \phi_n^2|\leqslant |\phi_n^2|,$$
  the dominated convergence Theorem ensures that 
  $$\int_{-1}^1 \mathbb {1}_{\mathcal{O}_{\varphi(k)}} \phi_n^2 \underset{k\rightarrow +\infty}\longrightarrow \int_\omega \phi_n^2.$$
  In other words,
\begin{equation}\label{om}\int_{\mathcal O_{\varphi(k)}}  \phi_n^2\underset{k\rightarrow +\infty}\longrightarrow  \int_{\omega}  \phi_n^2.\end{equation}
Moreover, since $\mathcal O_{\varphi(k)}\subset (\alpha_0,1)$ for any $k\in\mathbb N\setminus{\{0\}}$,  using Proposition \ref{corouvert}, there exists $C>0$ and $n_0\in\mathbb N\setminus{\{0\}}$, which might depend on $\alpha_0$ and $\nu$,  such that for any $k\in\mathbb N\setminus{\{0\}}$, we have 
    \[
       \int_{\mathcal O_{\varphi(k)}} \phi_n^2\geqslant C | \mathcal O_{\varphi(k)}|,\quad\forall n\geqslant n_0.
        \]  
        Letting $k\to +\infty$ in this inequality, we deduce by \eqref{om} that 
 
 \[
       \int_{\omega} \phi_n^2\geqslant C |\omega|,\quad\forall n\geqslant n_0.
        \]  

\noindent For  $n<n_0$, as a consequence of the analyticity of the Bessel functions on $\mathbb C\setminus \{0\}$ together with the fact that every normalized eigenfunction  $\phi_n$ for $n\in \N\setminus{\{0\}}$ is of the form 
\[
\phi_n(x)=\alpha_n\sqrt{x}J_\nu(\delta x)+\beta_n\sqrt{x}J_{-\nu}(\delta x),\,x>0,
\]
for some $\alpha_n,\beta_n\in\mathbb R\setminus{\{0\}}$, for some $\delta\neq 0$ and some $\nu\in(0,1)$, we clearly have 
\[
\int_\omega \phi_n^2 > 0,\quad\forall n\leqslant n_0,n>0,
\]
and from Proposition \ref{PropKernel}, for any normalized eigenfunction  $\phi_0$ associated to the eigenvalue $0$, we also have 
\[
\int_\omega \phi_0^2 > 0.
\]
Hence, the desired result is proved.
\qed

\vspace{0.3cm}

\noindent We are now ready for the proof of our main Theorem.

\vspace{0.3cm}

\noindent \textit{\textbf{Proof of Theorem \ref{theoremNullControlIntern}}}. Let $T>0$, let $\omega$ be a measurable subset of $(-1,1)$ assumed to have positive Lebesgue measure. Let $f^0\in L^2(-1,1)$. We know that $(\phi_n)_{n\in\N}$ is a Hilbert basis of $L^2(-1,1)$, as defined in Theorem \ref{propEigenBasis}. Therefore,
\[
f^0=\sum_{j\in\N}\left\langle f^0, \phi_j\right\rangle_{L^2(-1,1)} \phi_j.
\]

\noindent  We denote the unique solution of \eqref{1} in $C^0\left([0,T],L^2(-1,1)\right)$ by $
f$. We recall that, by Duhamel's formula, we have
\begin{equation}
f(t)=e^{-A_\nu t}f^0+\int_0^t e^{-A_\nu (t-s)}\mathbb{1}_\omega u(s)\mathrm{d}s, \quad \forall t\in[0,T].
\label{weakSolProofNullcontrol}
\end{equation}

Choosing $t=T$ in (\ref{weakSolProofNullcontrol}), and taking the scalar product with $\phi_n$ gives
\begin{equation}
    \begin{split}
\forall n\in\N, \quad \left\langle f(T), \ \phi_n\right\rangle&= \left\langle e^{-A_\nu T}f^0, \ \phi_n \right\rangle_{L^2(-1,1)} +\left\langle \int_0^T e^{-A_\nu(T-s)}\mathbb{1}_\omega u(s)\mathrm{d}s, \ \phi_n \right\rangle_{L^2(-1,1)}\\
&=e^{-\lambda_n T}\left\langle f^0, \ \phi_n \right\rangle_{L^2(-1,1)} + \int_0^T e^{-\lambda_n(T-s)}\left\langle   u(s), \ \mathbb{1}_\omega\phi_n \right\rangle_{L^2(-1,1)}\mathrm{d}s.
\end{split}
\label{proofNullContro}
\end{equation}
Following \cite[Section 5.3.3]{Trelat2024} (for instance), we choose the following control
\begin{equation}
u(t,x)=-\sum_{n\in\N} \left\langle f^0, \ \phi_n\right\rangle_{L^2(-1,1)} e^{-\lambda_nT} {q_n(T-t)}
\dfrac{\mathbb{1}_\omega \phi_n(x)}{\norm{\mathbb{1}_\omega \phi_n}_{L^2(-1,1)}^2},
    \label{rightControl}
\end{equation}
after having checked that it makes sense. \\
Proposition \ref{corIntBorel} ensures that for any $n\in\N$, $\norm{\mathbb{1}_\omega \phi_n}_{L^2(-1,1)}>0$. We know need to check that $u\in L^2\left((0,T)\times(-1,1)\right)$. First, we have for all $t\in[0,T]$
\[\begin{split}
\norm{u(t)}^2_{L^2(-1,1)}&=
\sum_{n\in\N}\sum_{k\in\N} \left\langle f^0, \phi_n\right\rangle \left\langle f^0, \phi_k\right\rangle e^{-\lambda_nT} e^{-\lambda_kT} \dfrac{q_n(T-t)}{\norm{\mathbb{1}_\omega \phi_n}}
\dfrac{q_k(T-t)}{\norm{\mathbb{1}_\omega \phi_k}}
\left\langle\dfrac{\mathbb{1}_\omega \phi_n}{\norm{\mathbb{1}_\omega \phi_n}}, \ \dfrac{\mathbb{1}_\omega \phi_k}{\norm{\mathbb{1}_\omega \phi_k}}\right\rangle.
\end{split}
\]
Moreover,
\[\begin{split}
\forall n,k\in\N, \ \left |\int_0^T \left\langle f^0, \ \phi_n\right\rangle \left\langle f^0, \ \phi_k\right\rangle \right .&\left .e^{-\lambda_nT} e^{-\lambda_kT} \dfrac{q_n(T-t)}{\norm{\mathbb{1}_\omega \phi_n}}
\dfrac{q_k(T-t)}{\norm{\mathbb{1}_\omega \phi_k}}
\left\langle\dfrac{\mathbb{1}_\omega \phi_n}{\norm{\mathbb{1}_\omega \phi_n}}, \ \dfrac{\mathbb{1}_\omega \phi_k}{\norm{\mathbb{1}_\omega \phi_k}}\right\rangle \mathrm{d}t \right |\\
&\leqslant
\frac{\abs{\left\langle f^0, \ \phi_n\right\rangle \left\langle f^0, \ \phi_k\right\rangle}}{\underset{j\in\N}{\inf}\norm{\mathbb{1}_\omega \phi_j}^2} e^{-\lambda_nT} e^{-\lambda_kT}
\norm{q_n}_{L^2(0,T)}
\norm{q_k}_{L^2(0,T)}\\
&\leqslant \frac{\abs{\left\langle f^0, \ \phi_n\right\rangle \left\langle f^0, \ \phi_k\right\rangle}}{\underset{j\in\N}{\inf}\norm{\mathbb{1}_\omega \phi_j}^2} e^{-\lambda_nT} e^{-\lambda_kT}C_{T/2}^2 e^{\lambda_n \frac{T}2}e^{\lambda_k \frac{T}2},
\end{split}
\]
where we used Proposition \ref{propBiorthogonal} for $\varepsilon=\frac{T}{2}$ to obtain an upper bound to $\norm{q_n}^2_{L^2(0,T)}$, and Proposition \ref{corIntBorel}, which gives $\underset{n\in\N}{\inf}\norm{\mathbb{1}_\omega \phi_n}^2_{L^2(-1,1)}>0$.
Therefore, we can use Fubini's theorem to get
\[\begin{split}
\int_0^T\norm{u(t)}^2\mathrm{d}t&=
\sum_{n\in\N}\sum_{k\in\N}\int_0^T \left\langle f^0, \ \phi_n\right\rangle \left\langle f^0, \ \phi_k\right\rangle e^{-\lambda_nT} e^{-\lambda_kT} \dfrac{q_n(T-t)}{\norm{\mathbb{1}_\omega \phi_n}}
\dfrac{q_k(T-t)}{\norm{\mathbb{1}_\omega \phi_k}}
\left\langle\dfrac{\mathbb{1}_\omega \phi_n}{\norm{\mathbb{1}_\omega \phi_n}}, \ \dfrac{\mathbb{1}_\omega \phi_k}{\norm{\mathbb{1}_\omega \phi_k}}\right\rangle \mathrm{d}t\\
&\leqslant \frac{C_{T/2}^2}{\underset{j\in\N}{\inf}\norm{\mathbb{1}_\omega \phi_j}^2}
\left(\sum_{n\in\N}\abs{\left\langle f^0, \ \phi_n\right\rangle}e^{-\lambda_n\frac{T}{2}}\right)^2
< +\infty \mathrm{\ because \ } f^0\in L^2(-1,1).
\end{split}
\]
 This shows that $u\in L^2\left((0,T)\times (-1,1)\right)$. Now, we examine what (\ref{proofNullContro}) gives when we choose this control. Using the biorthogonal property of $(q_n)_{n\in\N}$, we have
\[
   \begin{split}
\forall n\in\N, \quad
\int_0^T e^{-\lambda_n(T-s)}\left\langle   u(s), \ \mathbb{1}_\omega\phi_n \right\rangle_{L^2(-1,1)}\mathrm{d}s&=
-\left\langle f^0, \ \phi_n\right\rangle e^{-\lambda_nT}
\left\langle\dfrac{\mathbb{1}_\omega \phi_n}{\norm{\mathbb{1}_\omega \phi_n}_{L^2(-1,1)}^2},\ \mathbb{1}_\omega \phi_n\right\rangle_{L^2(-1,1)}\\
&=-\left\langle f^0, \ \phi_n\right\rangle e^{-\lambda_nT}.
\end{split}
\]
Therefore, from (\ref{proofNullContro}), we have
\[
\left\langle f(T), \ \phi_n\right\rangle_{L^2(-1,1)}=0,\quad\forall n\in\N,
\]
which means that $f(T)=0$. This concludes the proof of Theorem \ref{theoremNullControlIntern}.

\qed

\subsection{Boundary Control}
Let us now explain briefly how one can deduce an appropriate boundary control result for the system

\begin{equation} \left \{
\begin{aligned}\partial_tf(t,x)-\partial_{xx}f(t,x)+\dfrac{c}{x^2}f(t,x)&=0,& (t,x)\in(0,T)\times(-1, 1),\\f(t,-1)=0,\,f(t, 1)&=u(t),&t\in(0,T),\\
f(0,x)&=f^0(x),&x\in(-1, 1),
\end{aligned} \right .
\label{1b}
\end{equation}
where $u\in L^2(0,T)$ (the case where the control acts at $x=-1$ can be treated similarly). The first important point is to understand the well-posedness of \eqref{1b}, which is less obvious than for \eqref{1.3}, since the control operator is now unbounded. The procedure is classical and we only give the main ingredients.
Without loss of generality, we study instead 

\begin{equation} \left \{
\begin{aligned}\partial_tf(t,x)-\partial_{xx}f(t,x)+\dfrac{c}{x^2}f(t,x) +f(t,x)&=0,& (t,x)\in(0,T)\times(-1, 1),\\f(t,-1)=0,\,f(t, 1)&=u(t),&t\in(0,T),\\
f(0,x)&=f^0(x),&x\in(-1, 1).
\end{aligned} \right .
\label{1b2}
\end{equation}
Indeed, one can pass from one equation to another by multiplying the solutions by $e^{-t}$ (or $e^{t}$), without changing the well-posedness or controllability results. The main interest of system \eqref{1b2} is that the underlying elliptic operator $A_\nu+\mathrm{Id}$ is now (self-adjoint and) positive. The eigenvectors are unchanged and the eigenvalues are shifted by $1$.

\begin{proposition} For any $n\in\mathbb N\setminus\{0\}$, we have $ \phi_{n}'(1)\neq0$. Moreover, 

    \begin{equation}
        \phi_{2(n-1)}'(1)\underset{n\to+\infty}{=}(-1)^n\pi n\left(1+\mathcal{O}\left(\frac{1}{n^{\min(1,2\nu)}}\right)\right),
        \label{DerivateEven(1)}
    \end{equation}
    \begin{equation}
        \phi_{2(n-1)+1}'(1)\underset{n\to+\infty}{=}(-1)^n\pi n\left(1+\mathcal{O}\left(\frac{1}{n^{\min(1,2\nu)}}\right)\right),
        \label{DerivateOdd(1)}
    \end{equation}
    with $(\phi_n)_{n\in\N}$ the Hilbert basis of eigenfunctions defined in Theorem \ref{propEigenBasis}.
    \label{propLimDerivative(1)}
\end{proposition}

\noindent \textbf{\textit{Proof of Proposition \ref{propLimDerivative(1)}}}. 
The first assertion comes from the fact that $\phi_n(1)=0$. Since $\phi_n$ is an eigenfunction of $A_\nu$, it satisfies a second-order ODE that is not singular in a neighborhood of  $x=1$, so we cannot have simultaneously $\phi_n(1)=0$ and $\phi_n'(1)=0$.

\noindent Let $n\in\N\setminus\{0\}$, let $x\in(0,1]$, by looking at the expression of $\psi_{2(n-1)}$ given in Theorem \ref{propEigenBasis}, we get
\[\begin{split}
    \psi_{2(n-1)}'(x)&=\frac{1}{2\sqrt{x}}J_{-\nu}(\sqrt{\lambda_{2(n-1)}}x)+\sqrt{x}\sqrt{\lambda_{2(n-1)}}J_{-\nu}'(\sqrt{\lambda_{2(n-1)}}x)\\
    &-\frac{\Gamma(\nu+1)}{\Gamma(1-\nu)}\left(\frac{\sqrt{\lambda_{2(n-1)}}}{2}\right)^{-2\nu}\left(
    \frac{1}{2\sqrt{x}}J_{\nu}(\sqrt{\lambda_{2(n-1)}}x)+\sqrt{x}\sqrt{\lambda_{2(n-1)}}J_{\nu}'(\sqrt{\lambda_{2(n-1)}}x)\right).
\end{split}
\]
Therefore,
\[
\begin{split}
    \psi_{2(n-1)}'(1)=\frac{1}{2}&\psi_{2(n-1)}(1)\\
    +
    &\sqrt{\lambda_{2(n-1)}}\left(J_{-\nu}'(\sqrt{\lambda_{2(n-1)}})-\frac{\Gamma(\nu+1)}{\Gamma(1-\nu)}\left(\frac{\sqrt{\lambda_{2(n-1)}}}{2}\right)^{-2\nu}J_{\nu}'(\sqrt{\lambda_{2(n-1)}})\right)\\
    \underset{n\to+\infty}{=}& \ 0+\pi n\left(1+\mathcal{O}\left(\frac1n\right)\right)\left(J_{-\nu}'(\sqrt{\lambda_{2(n-1)}})+\mathcal{O}\left(\frac{1}{n^{2\nu}}\right)\mathcal{O}\left(\frac{1}{n^{1/2}}\right)\right),
\end{split}
\]
where we used Corollary \ref{corrAsymptSqrtEven(0,1)} and Lemma \ref{lemmaDerivateBessel}. Now, from Lemma \ref{lemmaDerivateBessel} and Lemma \ref{lemmaExpansion}, we get
\[\begin{split}
J_{-\nu}'(x)&\underset{x\to+\infty}{=}\frac12 \left(J_{-\nu-1}(x)-J_{1-\nu}(x)\right)\\
&\underset{x\to+\infty}{=}\frac12\left(\frac{2}{\pi x}\right)^{1/2}\left(\cos\left(x+\nu\frac{\pi}{2}+\frac{\pi}{2}-\frac{\pi}{4}\right)-\cos\left(x+\nu\frac{\pi}{2}-\frac{\pi}{2}-\frac{\pi}{4}\right)+\mathcal{O}\left(\frac{1}{x}\right)\right)\\
&\underset{x\to+\infty}{=}\left(\frac{2}{\pi x}\right)^{1/2}\left(\cos\left(x+\nu\frac{\pi}{2}+\frac{\pi}{4}\right)+\mathcal{O}\left(\frac{1}{x}\right)\right).
\end{split}
\]
We then replace $x$ by $\sqrt{\lambda_{2(n-1)}}$ in the expression above. Using Corollary \ref{corrAsymptSqrtEven(0,1)}, we get
\[\begin{split}
J_{-\nu}'(\sqrt{\lambda_{2(n-1)}})
&\underset{n\to+\infty}{=}
\frac{1}{\pi}\left(\frac{2}{n}\right)^{1/2}\left[1+\mathcal{O}\left(\frac1n\right)\right]\left[(-1)^n\cos\left(\mathcal{O}\left(\frac{1}{n^{\min(1,2\nu)}}\right)\right)+\mathcal{O}\left(\frac1n\right)\right]\\
&\underset{n\to+\infty}{=}
\frac{(-1)^n}{\pi}\left(\frac{2}{n}\right)^{1/2}\left[1+\mathcal{O}\left(\frac{1}{n^{\min(1,4\nu)}}\right)\right].
\end{split}
\]
Thus,
\[
\begin{split}
    \psi_{2(n-1)}'(1)\underset{n\to+\infty}{=}&\pi n\left(1+\mathcal{O}\left(\frac1n\right)\right)\frac{(-1)^n}{\pi}\left(\frac{2}{n}\right)^{1/2}\left(1+\mathcal{O}\left(\frac{1}{n^{\min(1,2\nu)}}\right)\right).
\end{split}
\]
Using the asymptotic expression of $a_{2(n-1)}$ at (\ref{asympt a2n-1 (0,1)}), we directly obtain (\ref{DerivateEven(1)}). We can follow the same computations to prove (\ref{DerivateOdd(1)}).
\qed

\vspace{0.2cm}

From Proposition \ref{propLimDerivative(1)} and Corollary \ref{Corasvp}, we see that there exists some constant $c_1,c_2>0$ such that for any $n\in\mathbb N\setminus{\{0\}}$, we have 
\begin{equation}\label{der1}c_1 \sqrt{\lambda_n+1}\leqslant |\phi_n'(1)|\leqslant c_2\sqrt{\lambda_n+1},\quad \forall n\in\mathbb N. \end{equation}

From this estimate, one can obtain the following well-posedness result. Since $A_\nu+\mathrm{Id}$ is positive, one can define $(A_\nu+\mathrm{Id})^{1/2}$. Moreover, $A_\nu+\mathrm{Id}$ extends as an unbounded operator with domain $D((A_\nu)^{1/2})$ and state space $D((A_\nu)^{-1/2})=D((A_\nu)^{1/2})'$ (see \cite[Corollary 3.4.6]{TW}), where for any $s\in \mathbb R$, $D((A_\nu)^{s})$ is equipped with the Hilbert norm
$$||f||_s^2= \sum_{k=1}^{+\infty}(\lambda_k+1)^{2s} |\langle f,\phi_k \rangle_{L^2(-1,1)}|^2,$$
  which makes it a Hilbert space. In particular, we obtain a Hilbert basis  $(\tilde \phi_n)_{ n\in\mathbb N}$ of eigenfunctions in $D(A_\nu^{-1/2})$ as 
$$\tilde \phi_n=\sqrt{\lambda_n+1}\phi_n,$$
and estimate \eqref{der1} becomes 
 
\begin{equation}\label{der2} (\lambda_n+1)c_1\leqslant |\tilde \phi_n'(1)|\leqslant c_2 (\lambda_{n}+1),\quad\forall n\in\mathbb N. \end{equation}

\noindent For the sake of simplicity, from now on, we introduce 
$$H_{\pm1/2}=D((A_\nu)^{\pm 1/2}).$$

Standard computations (see for instance \cite[Section 10.7]{TW}, where the usual Laplacian operator is replaced by $A_\nu$) show that in the state space $H_{1/2}$, the first line of \eqref{1b2} can be written in an abstract way as 
$$y'=(A_\nu+\mathrm{Id}) y+bu,$$
where $b\in H_{-1/2}$  is a scalar control operator 
that can be written as 
$b=(A_\nu+\mathrm{Id}) \mathcal D$, where $\mathcal D: \mathbb R \rightarrow L^2(-1,1)$ is the Dirichlet map that is given by duality by  (see \cite[Proposition 10.6.1]{TW}) 
$$\mathcal D^*g=-((A_\nu+\mathrm{Id})^{-1}g)'(1),\,\forall g\in L^1(-1,1).$$
In particular, for any $k\in\mathbb N$, 
$$\begin{aligned}b_k&:=\langle b, \tilde \phi_k\rangle_{H_{-1/2}}\\&= \frac{1}{\lambda_k+1}\langle b,\tilde \phi_k\rangle_{L^2(-1,1)}\\&= \frac{1}{\lambda_k+1}\langle (A_\nu+\mathrm{Id}) \mathcal D,\tilde\phi_k\rangle_{L^2(-1,1)}\\&= \frac{1}{\lambda_k+1} \langle 1,\mathcal D^*(A_\nu+\mathrm{Id})\tilde \phi_k\rangle_{\mathbb R}\\&= -\frac{1}{1+\lambda_k} \tilde \phi_k'(1).\end{aligned}$$
Hence, \eqref{der2} becomes 

\begin{equation}\label{der3}c_1\leqslant |b_k|\leqslant c_2,\, k\in \mathbb N. \end{equation}

The upper bound in \eqref{der3} together with the fact that there exists some constant $C>0$ such that for any $k\in\mathbb N\setminus{\{0\}}$, we have $1\leqslant 1+\lambda_k\leqslant Ck^2$ (by Corollary \ref{Corasvp}) easily imply that the sequence
$(b_k)_{k\in\mathbb N}$ verifies the Carleson measure criterion for admissibility given in \cite[Definition 5.3.1]{TW}.

From \cite[Theorem 5.3.2]{TW}, we deduce that $b$ is an admissible control operator, and we have the following well-posedness result: for any $f^0 \in D((A_\nu)^{1/2})'$, there exists a unique weak solution $f$ to \eqref{1b2} verifying $f\in C^0([0,T],D((A_\nu)^{1/2})')$.

We are now ready to give our null boundary controllability result in arbitrary small time.
\begin{theorem} 
      For any $f^0 \in D((A_\nu)^{1/2})'$ and any $T>0$,  there exists a control  $u \in L^2(0,T)$ such that  the corresponding solution $f$ to \eqref{1b2} verifies $f(T) = 0$ in $D((A_\nu)^{1/2})'$.
    \label{mainb}
\end{theorem}

\noindent\textbf{Proof of Theorem \ref{mainb}.}
Since $b$ is an admissible control operator, since the eigenvalues $\lambda_k+1$ are positive for $k\in\mathbb N$, and 
$\sum_{k=0}^{+\infty} \frac{1}{\lambda_k+1}<+\infty$ by Corollary \ref{Corasvp}, one can apply directly \cite[Theorem 2.5]{AmmarKhodja2013} and obtain that 
\eqref{1b2} is null controllable for any $T>T_0$, where 
$$T_0=\limsup_{k\rightarrow +\infty} \left (\frac{\log\frac{1}{|b_k|}}{\lambda_k+1}+ c(\Lambda+1)\right).$$
Here, $c(\Lambda+1)$ is the condensation index of the sequence $\{\lambda_{n+1}+1\}_{n\in\mathbb N\setminus{\{0\}}}$, as defined in Definition \ref{defCondensationIndex}.
From the lower bound given in \eqref{der3}, we have that 
$$\frac{\log\frac{1}{|b_k|}}{\lambda_k+1}\underset{k\rightarrow +\infty}\rightarrow 0.$$
Moreover, it is easy to infer that $c(\Lambda+1)=c(\Lambda)$, where $c(\Lambda)$ is defined in Proposition \ref{propNullCondensIndex} and verifies $c(\Lambda)=0$. 
Hence, 
$T_0=0$, and our result follows.
\qed

\section{Further comments and open problems}
\label{sec:c}
\paragraph{The case $c_\nu=-\frac{1}{4}$.} We have cautiously excluded this case in our study (as in \cite{Morancey2015}), whereas this critical case does not lead necessarily to ill-posedness (the other critical case $c=\frac{3}{4}$ is impossible, on the contrary). However, the functional setting is not clear.  Indeed, this is equivalent to imposing $\nu=0$, and we lose the coercivity estimate \eqref{1.6}, which is in particular crucial for the definition of the domain of $A_\nu$. In the case $\nu=0$, we need to
change the definition of $D(A_\nu)$ according to \cite[Proposition 3.1]{AA}, and to find a suitable replacement for \eqref{1.6}.  In fact, \eqref{1.6} is proved by using the following Hardy inequality: 
$$\int_{-1}^1\frac{z(x)^2}{x^2}\leqslant 4\int_{-1}^1z_x(x)^2,\,\forall z\in H^{-1}(-1,1)\mbox{ such that }z(0)=0,$$
that might be replaced in the case $\nu=0$ by an improved Hardy-Poincar\'e type inequality like in \cite[Theorem 2.2]{VZ2}. This case would require an appropriate and extensive treatment that is likely to be different from the situation we presently studied, and is outside the scope of this article.

\paragraph{Backstepping.} Since we proved null controllability in any time $T>0$, we know that rapid  stabilization (\textit{i.e.} exponential stabilization at any arbitrarily large rate) holds (see \cite[Proposition 21 and Theorem 25]{TWZ}). However, it would be interesting to construct an explicit feedback using the backstepping method. This would require to solve a second order PDE of wave type,  with two singular potentials, one depending on space and one on  time. We were not able to understand in which functional setting such an equation could be well-posed.

\paragraph{Other self-adjoint extensions.} As highlighted before, our result highly depends on the self-adjoint extension that we use for the Laplace operator with core \(C_0^\infty((-1,0)\cup(0,1))\). The ``natural'' one (coming from \cite{V}) leads to a lack of null controllability from one side, whereas choosing the one of the present article leads to a positive null controllability result in arbitrary small time. However, amongst the extensions that lead to Dirichlet boundary conditions at the endpoints $\pm 1$, there are a lot of other extensions that are likely to lead to positive controllability results (the ones that are called the ``coupling extensions'' in definition \ref{def:coupled-ext}). Amongst these extensions, we need to understand which ones are non-negative, or at least semi-bounded from below, so that they generate a semigroup by \cite[Corollary 2.4.8]{CazenaveHaraux1990}. However, from expression \eqref{Affp}, it is not totally clear how to determine these extensions, which is a first technical difficulty. We also need to ensure that our extensions have compact resolvent. As soon as these difficulties are solved, we expect the rest of our computations to be quite similar to the present article, especially the study of the eigenfunctions (that are likely to be only slight modifications of \eqref{eigenSpace+} and \eqref{eigenSpace-}, with the same Bessel function appearing, but with different linear combinations) and the eigenvalues, for which we expect to have a similar behaviour, in the sense there are groups of two eigenvalues that become arbitrarily close, but  not too rapidly, so that the condensation index of the sequence of eigenvalues should still be $0$. Hence, we expect the same kind of null controllability result  in arbitrary small time to hold. There are still some technical points to understand, and the  computations are likely to be quite heavy, so we delay this general study to a potential future work. For the case of more general boundary conditions at the endpoints $\pm 1$, the question is totally open and likely to be more difficult.

\paragraph{Mixed singular/degenerate equations.} In \cite{V}, the author proved some controllability results for  parabolic equations posed on $(0,1)$, with degenerate diffusion and singular potential at the boundary $x=0$. It would be interesting to understand if these results can be extended in our case of an interior degeneracy/singularity and a control region that is only located in one side of the boundary. It would require to develop a similar functional setting, with adequate transmission conditions, but we expect the computations to be much heavier.

\paragraph{Grushin equation with a singularity.} The original motivation of \cite{Morancey2015} was to study the controllability properties of the Grushin equation with inverse square potential

$$\partial_tf(t,x,y)-\partial_{xx}f(t,x,y)-x^2\partial_{yy}f(t,x,y)+\dfrac{c}{x^2}f(t,x,y)=u(t, x,y)\mathbb{1}_\omega(x,y),\, (t,x,y)\in(0,T)\times(-1, 1)\times (0,1),$$
with Dirichlet boundary condition on the boundary of $(-1,1)\times(0,1)$, and $\omega$ is an open subset on $(-1,1)\times (0,1)$. This kind of model naturally arises when looking at the Laplace operator on almost-Riemannian manifolds. In \cite{Morancey2015}, approximate controllability is proved. It would be interesting to understand if null controllability can hold at least in some particular geometrical settings and in large enough time, when controlling for instance in a vertical strip that is located at one side of the singularity $x=0$, as in \cite{BCG} (where there is no singular potential), by using a Fourier decomposition in the $y$ variable and reducing the problem to a uniform observability problem. In view of the present study, we expect that this question may be rather difficult.

\section*{Acknowledgements}
This work was funded by the French Agence Nationale de la Recherche (Grant ANR-22-CPJ2-0138-01).

\vspace*{0.2cm}

\noindent The first author would like to thank Eric Cances, David Gontier and Dario Prandi for having suggested the strategy of Proposition \ref{p:app}.

\vspace*{0.2cm}

\noindent The authors would like to warmly thank the anonymous referees for their criticisms and remarks, that greatly helped us to improve the presentation and clarity of the paper.

\appendix

\section{A collection of well-known results on the Bessel functions}
\label{AB}

\begin{lemma}{\cite[Subsection 3.12]{Watson1944}}\\
   Assume that ${\tilde{\nu}}\in\R\setminus \mathbb Z$. Consider the Wronskian of $J_{\tilde{\nu}}$ and $J_{-{\tilde{\nu}}}$, defined as
    $$W\left(J_{\tilde{\nu}},J_{-{\tilde{\nu}}}\right):= J_{\tilde{\nu}} J_{-{\tilde{\nu}}}'-J_{\tilde{\nu}}' J_{-{\tilde{\nu}}}.
    $$ Then,
\[W\left(J_{\tilde{\nu}},J_{-{\tilde{\nu}}}\right)(x)=-\frac{2\sin({\tilde{\nu}}\pi)}{\pi x}, \quad \forall x>0.
\]
In particular, for ${\tilde{\nu}}\in\R\setminus \mathbb Z$,  $J_{{\tilde{\nu}}}$ and $J_{-{\tilde{\nu}}}$ form a basis of solutions for (\ref{besseleq}).

    \label{lemmaWronskien}
\end{lemma}

\begin{lemma}{\cite[(2), p. 82]{Watson1944}}
    \[
        \forall {\tilde{\nu}} \in \R, \quad {J_{{\tilde{\nu}}}}'=\frac12\left(J_{{\tilde{\nu}}-1}-J_{{\tilde{\nu}}+1}\right).
    \]
    \label{lemmaDerivateBessel}
    \end{lemma}
\begin{lemma}{\cite[(1), p.199]{Watson1944}}
    \[
        J_{\tilde{\nu}}(x) \underset{x\to+\infty}{=} \sqrt{ \frac{2}{\pi x}}  \left( 
    \cos \omega \sum_{k=0}^{n} (-1)^k \frac{a_{2k}({\tilde{\nu}})}{x^{2k}} 
    - \sin \omega \sum_{k=0}^{n} (-1)^k \frac{a_{2k+1}({\tilde{\nu}})}{x^{2k+1}} +\mathcal{O}\left(\frac{1}{x^{2n+2}}\right)
    \right), 
    \quad \forall {\tilde{\nu}}\in\R, \ \forall n\in\N,
    \]
    where \( \omega = x -{\tilde{\nu}}\dfrac{\pi}{2} - \dfrac{\pi}{4} \), and $\ a_0=1$, $a_k=\dfrac{(4{\tilde{\nu}}^2-1^2)(4{\tilde{\nu}}^2-3^2)...(4{\tilde{\nu}}^2-(2k-1)^2)}{k!8^k}$ for $k\in\N\setminus\{0\}$.
    
        \label{lemmaExpansion}
    \end{lemma}

    \begin{lemma}{\cite[(6.03), Section 6.5]{O4}} Assume that ${\tilde{\nu}} >-1$. We have 
    
\[            j_{{\tilde{\nu}},n} \underset{n\to+\infty}{=}\pi\left(n+\frac{{\tilde{\nu}}}{2}-\frac14\right)+\mathcal{O}\left(\frac{1}{n}\right).
\]
 \label{lemmaAsymptoticZero}
\end{lemma}

\section{Characterization and non-negativity of Self-adjoint extensions}
\label{appex}
In this whole section, we fix $\nu\in(0,1)$. Until Corollary \ref{corextdirichlet}, we follow the work of Morancey  \cite[Section 2.4]{Morancey2015}.

\subsection{Minimal and maximal domains, useful results                           }

According to \cite[Proposition 3.1]{AA}, the minimal and maximal domains associated to the
differential expression $A_\nu f=-\partial_{xx}f+\frac{c_{\nu}}{x^{2}}f
$ in $L^2(0,1)$ are respectively equal to
\[
H^2_0([0,1]) := \{\, y \in H^2([0,1]) \; ; \; y(0)=y(1)=y'(0)=y'(1)=0 \,\}
\]
and
\[
\{\, y \in H^2([0,1]) \,;\; y(0)=y'(0)=0 \,\} 
\oplus \mathrm{Span}\bigl\{\, x^{\nu+\frac12},\, x^{-\nu+\frac12} \,\bigr\}.
\]

\medskip

Then, \cite[Lemma~13.3.1]{Zettl} implies that the minimal and maximal domains associated to $A_\nu$ on the interval $(-1,1)$ are given by
\[
D_{\min} := 
\Bigl\{\, f \in \widetilde{H}^2_0(-1,1) \,;\;
f(-1)=f(1)=f'(-1)=f'(1)=0 \,\Bigr\},
\]
and
\[
D_{\max} = \widetilde{H}^2_0(-1,1) \oplus \mathcal{F}_s^\nu,
\]
with $\widetilde{H}^2_0(-1,1)$ and $\mathcal{F}_s^\nu$ respectively defined in \eqref{H20} and \eqref{Fs}.

\medskip

\noindent We define as well
\[
H^2_0([-1,0]) := \{\, y \in H^2([0,1]) \; ; \; y(-1)=y(0)=y'(-1)=y'(0)=0 \,\},
\]
which will be useful for Proposition \ref{propextautoadj}.
\medskip

We state below a useful Lemma from \cite{Zettl} on the Lagrange bracket, 
which will be instrumental in the characterization 
of the self-adjoint extensions of our operator.

\begin{lemma}{\cite[Lemma 9.2.3]{Zettl}}\\
\label{lem:green}For $f,g \in D_{max}=\widetilde{H}^2_0(-1,1) \oplus \mathcal{F}_s^\nu$, define
\[
[f,g](x) := (f g' - f' g)(x), \qquad \forall x \neq 0.
\]
Then,
\[
\begin{aligned}
\int_{-1}^{1} \left(-\partial_{xx}^2 f + \frac{c_\nu}{x^2} f\right)(x) g(x)\,\mathrm{d}x
&= \int_{-1}^{1} f(x) \left(-\partial_{xx}^2 g + \frac{c_\nu}{x^2} g\right)(x)\,\mathrm{d}x \\
& +  [f,g](1) - [f,g](0^+) + [f,g](0^-) - [f,g](-1).
\end{aligned}
\]
\end{lemma}
We will also need the following result on the behavior at $x=0$ of functions in $\widetilde{H}^2_0(-1,1)$.
\begin{lemma}{\cite[Lemma 2.2]{Morancey2015}}\\
\label{lem:behavior-H20}For any $f\in\widetilde{H}^2_0(-1,1)$,
\[
\lim_{x\to \pm0} \frac{f(x)}{x^{3/2}} = 0
\qquad\text{and}\qquad
\lim_{x\to \pm0} \frac{f'(x)}{x^{1/2}} = 0.
\]
\end{lemma}

\subsection{All self-adjoint extensions}

The question of finding the self-adjoint extensions of a given closed symmetric operator is classical. The particular case of Sturm-Liouville operators has been widely studied. Most of these results
are contained in \cite{Zettl}. In our case, we are concerned with the Sturm-Liouville operator $A_\nu$ on the interval $(-1,1)$. This fits in the setting of
\cite[Chapter 13]{Zettl}. The number of transmission conditions to impose is given by the deficiency index. Following \cite[Proposition 3.1]{AA}, it comes that our operator on the interval $(0, 1)$ has deficiency index $2$. Then, \cite[Lemma 13.3.1]{Zettl} implies that the deficiency index for the interval $(-1, 1)$ is $4$. We thus get the following proposition, which is simply a rewriting of \cite[Theorem 13.3.1, Case 5]{Zettl}.

\medskip

\begin{proposition}{\cite[Theorem 13.3.1]{Zettl}}.
\label{propextautoadj}
Let $u$ and $v$ in $D_{\max}$ be such that their restriction on $(0,1)$ 
(resp.\ $(-1,0)$) are linearly independent modulo $H^2_0([0,1])$ 
(resp.\ $H^2_0([-1,0])$), and
\[
 \left[u, v\right](-1) = \left[u, v\right](0^-) = \left[u, v\right](0^+) = \left[u, v\right](1) = 1.
\]
Let $N_1,N_2,N_3,N_4$ be $4\times 2$ matrices. 
Then every self-adjoint extension of the minimal operator is given by the restriction 
of $D_{\max}$ to the functions $f$ satisfying the transmission conditions
\begin{equation}
\label{condN0}    
N_1
\begin{pmatrix}
[f,u](-1)\\[0.3em]
[f,v](-1)
\end{pmatrix}
+
N_2
\begin{pmatrix}
[f,u](0^-)\\[0.3em]
[f,v](0^-)
\end{pmatrix}
+
N_3
\begin{pmatrix}
[f,u](0^+)\\[0.3em]
[f,v](0^+)
\end{pmatrix}
+
N_4
\begin{pmatrix}
[f,u](1)\\[0.3em]
[f,v](1)
\end{pmatrix}= 0,
\end{equation}
where the matrices satisfy that $(N_1\ N_2\ N_3\ N_4)$ has full rank and
\begin{equation}
N_1 E N_1^T - N_2 E N_2^T + N_3 E N_3^T - N_4 E N_4^T = 0,
\quad 
\textrm{with} \quad  E := 
\begin{pmatrix}
0 & -1\\
1 & 0
\end{pmatrix}.
\label{condN}
\end{equation}
Conversely, every choice of such matrices defines a self-adjoint extension.
\end{proposition}

\subsection{Reduction under homogeneous Dirichlet boundary conditions}
Now, we focus on the case where homogeneous Dirichlet conditions are imposed at the endpoints
so that we are in the framework that we studied.

\medskip

For $f\in D_{max}$, we now define the quantities $\alpha^+(f)$, $\alpha^-(f)$, $\beta^+(f)$, $\beta^-(f)$, which will be useful to compactly define the self-adjoint extensions in our case.
\medskip

\begin{definition}
Let $f \in D_{\max}= \widetilde{H}^2_0(-1,1) \oplus \mathcal{F}_s^\nu$.  
As explained in Remark \ref{rqUniqueDecomposition}, $f$ admits a unique decomposition
\[
f = f_r + f_s, \quad \textrm{such \ that \ } 
f_r \in \widetilde{H}^2_0(-1,1) \mathrm{\ and \ } 
f_s \in \mathcal{F}_s,
\]
where
\[
f_s(x) =
\begin{cases}
c_1^-\, |x|^{\nu+\frac12} + c_2^-\, |x|^{-\nu+\frac12}, & \quad \forall x\in(-1,0),\\[0.3em]
c_1^+\, |x|^{\nu+\frac12} + c_2^+\, |x|^{-\nu+\frac12}, &\quad \forall x\in(0,1).
\end{cases}
\]
We now define the boundary coefficients
\[
\alpha^+(f) := c_1^+ + c_2^+, 
\qquad 
\alpha^-(f) := c_1^- + c_2^-,
\]
and
\[
\beta^+(f) := \bigl(\nu+\tfrac12\bigr)c_1^+ + \bigl(-\nu+\tfrac12\bigr)c_2^+, 
\qquad 
\beta^-(f) := \bigl(\nu+\tfrac12\bigr)c_1^- + \bigl(-\nu+\tfrac12\bigr)c_2^-.
\]

\medskip

\noindent When there is no ambiguity, we simply write $
\alpha^\pm \text{ and } \beta^\pm
$ instead of $\alpha^\pm(f)$ and $\beta^\pm(f)$.
\label{def:alphasbetas}
\end{definition}

\begin{proposition}
\label{corextdirichlet}
Consider any self-adjoint extension  $(\tilde A,D(\tilde A))$ of the minimal operator that imposes homogeneous Dirichlet conditions at the endpoints:
\begin{equation}\label{condd}
f\in D(\tilde A) \Rightarrow f(-1)=f(1)=0.
\end{equation}

\medskip

\noindent Then, such a self-adjoint extension is given by the restriction of $D_{\max}$ to the functions $f$ satisfying the homogeneous Dirichlet conditions $f(-1)=f(1)=0$ and the boundary condition
\begin{equation}
M_2
\begin{pmatrix}
\alpha^-(f)\\[0.3em]
-\beta^-(f)
\end{pmatrix}
+
M_3
\begin{pmatrix}
\alpha^+(f)\\[0.3em]
\beta^+(f)
\end{pmatrix}= \begin{pmatrix}0\\0\end{pmatrix},
\label{condbord}
\end{equation}
where $M_2,M_3$ are $2\times2$ matrices satisfying
\begin{equation}
\operatorname{rank}(M_2\ M_3)=2,
\qquad
\det(M_2)=\det(M_3).
\label{condrankdetM}
\end{equation}

Conversely, every choice of matrices $M_2$ and $M_3$ satisfying the above conditions defines a self-adjoint extension under homogeneous Dirichlet boundary conditions.
\end{proposition}
\noindent\textbf{Proof of Proposition \ref{corextdirichlet}}.\\
We apply Proposition \ref{propextautoadj} with the following choices of $u$ and $v$. We define on $(0,1]$ the functions $u$ and $v$ as the solutions of
\[
- f''(x) + \frac{c_\nu}{x^2} f(x) = 0
\]
with 
\[
u(1)=0, \ u'(1)=1, \quad v(1)=-1,\ v'(1)=0,
\]
i.e.
\[
u(x) = \frac{1}{2\nu}x^{\nu+\frac12} - \frac{1}{2\nu}x^{-\nu+\frac12}, 
\quad
v(x) = -\frac{\nu-\frac12}{2\nu}x^{\nu+\frac12} - \frac{\nu+\frac12}{2\nu}x^{-\nu+\frac12},\quad \forall x\in(0,1].
\]
Thus, by Lemma \ref{lem:green}, one has $[u,v](x) =1$, for all $x\in (0,1]$.

\medskip
\noindent We define $u$ and $v$ similarly on $[-1,0)$, i.e.
\[
u(x) = -\frac{1}{2\nu}|x|^{\nu+\frac12} + \frac{1}{2\nu}|x|^{-\nu+\frac12}, 
\quad
v(x) = -\frac{\nu-\frac12}{2\nu}|x|^{\nu+\frac12} - \frac{\nu+\frac12}{2\nu}|x|^{-\nu+\frac12}, \quad \forall x\in[-1,0)
\]
Using Lemma \ref{lem:green}, one has $
[u,v](x) =1$, for all $x\in [-1,0)$.
\medskip

Therefore, we can apply Proposition \ref{propextautoadj} with this choice of $u$ and $v$.

\medskip
Using once again Lemma \ref{lem:green}, we compute
\[
[f,u](1) =f(1), \quad  [f,v](1)=f'(1), \quad 
[f,u](-1) =f(-1), \quad  [f,v](-1)=f'(-1), \quad \forall f\in D_{max},
\]
and
\[
[f,u](0^+) =\alpha^+, \quad  [f,v](0^+)=\beta^+, \quad 
[f,u](0^-) =\alpha^-, \quad  [f,v](0^-)=-\beta^-, \quad \forall f\in D_{max}.
\]
Therefore, \eqref{condN0} becomes 
\begin{equation}
    \label{condN1}
N_1
\begin{pmatrix}
f(1)\\[0.3em]
f'(-1)
\end{pmatrix}
+
N_2
\begin{pmatrix}
\alpha^-\\[0.3em]
\beta^-
\end{pmatrix}
+
N_3
\begin{pmatrix}
\alpha^+\\[0.3em]
-\beta^+
\end{pmatrix}
+
N_4
\begin{pmatrix}
f(1)\\[0.3em]
f'(1)
\end{pmatrix}= 0.\end{equation}

Let us first  assume  that we have a self-adjoint extension $(\tilde A,D(\tilde A))$ that imposes  homogeneous Dirichlet conditions at the endpoints, \textit{i.e.} \eqref{condd}.
According to  Proposition \ref{propextautoadj}.
$D(\tilde A)$ is described as those $f\in D_{max}$ which verify  \eqref{condN1} and \eqref{condN}, where $N=(N_1\ N_2\ N_3\ N_4)$ has full rank.

Moreover, it is easy to prove that for any real numbers $a_1,\ldots a_8$, there exists a nonzero $f\in D_{max}$ such that 

\[
f(1)=a_1, \quad  f'(1)=a_2, \quad 
f(-1)=a_3, \quad  f'(-1)=a_4,
\]
and
\[\alpha^+=a_5, \quad  \beta^+=a_6, \quad \alpha^-=a_7, \quad  -\beta^-=a_8.
\]
Hence, in view of \eqref{condN1}, \eqref{condd}
in particular implies that 
$$\mathrm{Ker}(N)\subset \mathrm{Ker}(P),$$
where 
$$P=\begin{pmatrix} 
1&0&0&0&0&0&0&0\\0&0&0&0&0&0&1&0
\end{pmatrix}.$$
By usual factorization theorem, there exists some matrix $Q$ of size $2\times 4$ such that 
$P=QN$. Remark that $Q$ is necessarily of maximum rank $2$.
Write $Q$ as
$$Q=\begin{pmatrix} 
q_1&q_2&q_3&q_4\\q_5&q_6&q_7&q_8
\end{pmatrix}.$$
We introduce 
$$\tilde Q=\begin{pmatrix} 
q_1&q_2&q_3&q_4\\\cdot&\cdot&\cdot&\cdot\\\cdot&\cdot&\cdot&\cdot\\q_5&q_6&q_7&q_8
\end{pmatrix},$$
where the $\cdot$ are chosen in such a way that $\tilde Q$ is invertible, which is always possible since $Q$ is of rank $2$. Then, we observe that 
$$\tilde Q N=\begin{pmatrix} 
1&0&0&0&0&0&0&0\\\multicolumn{2}{c}{ M_1}&\multicolumn{2}{c}{ M_2}&\multicolumn{2}{c}{ M_3}&\multicolumn{2}{c}{ M_4}\\0&0&0&0&0&0&1&0
\end{pmatrix},$$
for some matrices $M_1,\ldots M_4$ of size $2\times 2$, and some invertible $\tilde Q$. Therefore, 
$$\mathrm{Ker}(N)=\mathrm{Ker}(\tilde QN).$$
Moreover, clearly, if for some column vector $X=(x_1,\ldots, x_8)^t$, we have 
$\tilde QN X=0$, then, necessarily, $x_1=x_7=0$, so that without loss of generality, we may assume that $M_1$ and $M_4$ are of the form
$$ M_i = \begin{pmatrix} 0&m_{i,1}\\0&m_{i,2} \end{pmatrix},\, i=1,4,$$
since it will not change the kernel of $\tilde Q N$. 
Now, explicit computations give that Condition (\ref{condN}) is equivalent to 

$$\begin{pmatrix} 
0&-m_{11}&-m_{12}&0\\m_{11}&0&\det(M_3)-\det(M_2)& -m_{41}\\m_{12}&\det(M_2)-\det(M_3)&0&-m_{42}\\0&0&m_{41}&m_{42}
\end{pmatrix}=\begin{pmatrix} 
0&0&0&0\\0&0&0&0\\0&0&0&0\\0&0&0&0
\end{pmatrix}.$$
We deduce that $M_1=M_4=0$, and $det(M_2)=det(M_3)$. So, finally, we can always assume that 

$$\tilde Q N=\begin{pmatrix} 
1&0&0&0&0&0&0&0\\\multicolumn{2}{c}{ 0}&\multicolumn{2}{c}{ M_2}&\multicolumn{2}{c}{ M_3}&\multicolumn{2}{c}{0}\\0&0&0&0&0&0&1&0
\end{pmatrix},$$
so that we deduce that condition \eqref{condd} implies that for some $2\times2$ matrices $M_2,M_3$, \eqref{condbord} and  (\ref{condrankdetM}) hold. 

The converse part is obvious. This concludes the proof of Corollary \ref{corextdirichlet}.
\qed 
\begin{remark}
\label{rq:types_extensions}
Corollary~\ref{corextdirichlet} provides a complete characterization of all self-adjoint extensions 
of the operator $A_\nu$ under homogeneous Dirichlet boundary conditions. 
It shows that any such extension can be encoded by a pair of matrices $(M_2,M_3)$ 
satisfying the rank and determinant constraints (\ref{condrankdetM}).

\medskip

\noindent
Of course, different pairs of matrices may define the same self-adjoint extension 
(for instance, $(M_2,M_3)$ and $(2M_2,2M_3)$ yield identical domains). 

\medskip

\noindent

\end{remark}

Among these extensions, two distinct situations can be identified, depending of whether $\det(M_2)=\det(M_3)$ is zero or not. This is the purpose of the two following definitions.

\begin{definition}[Coupled self-adjoint extensions]
\label{def:coupled-ext}
Let $A_\nu$ be endowed with homogeneous Dirichlet conditions at $\pm 1$.
According to Corollary~\ref{corextdirichlet}, any self-adjoint extension is then determined by a pair of matrices
$(M_2,M_3)$. We say that such an extension is a coupled extension whenever
\[
\det(M_2)=\det(M_3)\neq 0.
\]
We denote by
$
\mathfrak{E}_{\mathrm{cpl}}$ 
the set of all coupled self-adjoint extensions of $A_\nu$ under homogeneous
Dirichlet boundary conditions.
\end{definition}
Corollary~\ref{corextdirichlet} together with easy algebraic manipulations  directly imply the following property (that justifies the terminology ``coupled'', since we have nontrivial relations between $f(0^-), f'(0^-),f(0^+)$ and $f'(0^+)$.
\begin{proposition}\label{propcoup}
Assume that $A_\nu$ is a coupled extension. The extension defined by the invertible matrices $M_2$ and $M_3$ is equivalent to an extension
parameterized by a unique matrix $M$ with $\det(M)=1$,
through the transmission condition
\begin{equation}\label{eq:transmission-coupled}
\begin{pmatrix}
\alpha^-\\[0.3em]
-\beta^-
\end{pmatrix}
+
M
\begin{pmatrix}
\alpha^+\\[0.3em]
\beta^+
\end{pmatrix}
= 0.
\end{equation}
Its domain is
\[
D(A_\nu^M)
=
\Bigl\{
f\in D_{\max} ; \  f(\pm 1)=0 \text{ and \ (\ref{eq:transmission-coupled}) \ is  \ satisfied}
\Bigr\}.
\]
Hence,
\[
\mathfrak{E}_{\mathrm{cpl}}
:=
\Bigl\{
\,(A_\nu,\, D(A_\nu^M)) \ ;\ M\in \mathrm{SL}_2(\R)\,\Bigr\}.
\]
\end{proposition}
Let us now introduce the other family of extensions.

\begin{definition}[Decoupled self-adjoint extensions]
Let $A_\nu$ be endowed with homogeneous Dirichlet conditions at $\pm 1$.
According to Corollary~\ref{corextdirichlet}, any self-adjoint extension is then determined by a pair of matrices
$(M_2,M_3)$. We say that such an extension is a decoupled extension whenever
\[
\det(M_2)=\det(M_3)=0,
\qquad
\operatorname{rank}(M_2\ M_3)=2.
\]
\end{definition}

Then, we have the following proposition, which justifies the terminology of a ``decoupled'' extension.
\begin{proposition}
\label{propDecoupledExt}
Assume that $A_\nu$ is a decoupled extension.
Then the transmission condition \eqref{condbord} reduces to two independent
linear relations
\begin{equation}\label{eq:decoupledm}
\ell^-_1\,\alpha^- + \ell^-_2\,\beta^- = 0,
\end{equation}
\begin{equation}\label{eq:decoupledp}
\ell^+_1\,\alpha^+ + \ell^+_2\,\beta^+ = 0,
\end{equation}
for some nonzero vectors $(\ell^-_1,\ell^-_2)$ and $(\ell^+_1,\ell^+_2)$.
In particular, no transmission condition couples the left and right coefficients
$(\alpha^-,\beta^-)$ and $(\alpha^+,\beta^+)$.
\end{proposition}

\noindent\textbf{Proof of Proposition \ref{propDecoupledExt}}.\\
We consider a self-adjoint extension defined as in Corollary \ref{corextdirichlet} by the choice of two matrices $M_2$ and $M_3$ such that $\det(M_2)=\det(M_3)=0$ and $\operatorname{rank}(M_2 \ M_3)=2$.
Then each matrix $M_2$ and $M_3$ is nonzero of rank $1$. Hence, one of its two rows is nonzero and the other is proportional to it. Writing
\[
 M_2=\begin{pmatrix}a_1&b_1\\ k_1 a_1&k_1 b_1\end{pmatrix}
\quad\text{with }(a_1,b_1)\neq(0,0)\quad\text{or}\quad
 M_2=\begin{pmatrix}k_1 c_1&k_1 d_1\\ c_1&d_1\end{pmatrix}
\quad\text{with }(c_1,d_1)\neq(0,0),
\]
and similarly
\[
M_3=\begin{pmatrix}a_2&b_2\\ k_2 a_2&k_2 b_2\end{pmatrix}
\quad\text{with }(a_2,b_2)\neq(0,0)\quad\text{or}\quad
 M_3=\begin{pmatrix}k_2 c_2&k_2 d_2\\ c_2&d_2\end{pmatrix}
\quad\text{with }(c_2,d_2)\neq(0,0),
\]
we obtain four canonical placements ($k_j$ on the first or on the second row).
The rank–$2$ condition on the concatenation imposes the independence of the two
resulting nonzero rows, which is equivalent to either $k_1\neq k_2$ when $k_1,k_2$ sit on the same row, or $k_1k_2\neq 1$ when they sit on different rows.
\medskip

Since the four cases are very similar, we compute the transmission conditions only for the case $k_1$ and $k_2$ on the second row. We have
\[
 M_2=\begin{pmatrix}a_1&-\,b_1\\ k_1 a_1&-\,k_1 b_1\end{pmatrix},
\qquad
 M_3=\begin{pmatrix}a_2&b_2\\ k_2 a_2&k_2 b_2\end{pmatrix},
\qquad (a_1,b_1)\neq(0,0),\ (a_2,b_2)\neq(0,0),\ k_1\neq k_2.
\]
With these matrices, the transmission conditions (\ref{condbord}) are
\[
\begin{cases}
a_1\alpha^- - b_1\beta^- + a_2\alpha^+ + b_2\beta^+ = 0,\\
k_1 a_1\alpha^- - k_1 b_1\beta^- + k_2 a_2\alpha^+ + k_2 b_2\beta^+ = 0.
\end{cases}
\]
which gives, since $k_1\neq k_2$,
\[
a_2\alpha^- = b_2\beta^-
\qquad\text{and}\qquad
a_3\alpha^+ = -b_3\beta^+.
\]
If $k_1$ was sitting on the first row of $ M_2$ with $(c_1,d_1)\neq(0,0)$, we would get the relation $c_1\alpha^- = d_1\beta^-$, and it works the same for $M_3$. Therefore, the general conclusion is that the transmission conditions are of the form \eqref{eq:decoupledm} and \eqref{eq:decoupledp}.

\medskip
To conclude, observe that the relations \eqref{eq:decoupledm} and \eqref{eq:decoupledp} impose one
boundary condition at $0^-$ and one at $0^+$, with no coupling between the two
sides.  On each half-interval, $A_\nu$ is a Sturm–Liouville operator with a
regular endpoint at $\pm1$ and a single singular endpoint at $0^\pm$. We are exactly in the framework
of \cite[Proposition~10.4.2]{Zettl}, which classifies self-adjoint extensions
for Sturm–Liouville problems with one singular endpoint.
\qed

\begin{remark}
\label{rem:decoupled-not-interesting}
Proposition~\ref{propDecoupledExt} shows that when the matrices $M_2$ and $M_3$
have rank $1$, the resulting self-adjoint extension  $A_\nu$ splits into two independent
self-adjoint problems on $(-1,0)$ and $(0,1)$. 
More precisely, 
\eqref{eq:decoupledm} and \eqref{eq:decoupledp} define two self-adjoint extensions of $A_\nu$, denoted
$A_\nu^-$ and $A_\nu^+$.  Their domains are
\[
D(A_\nu^-)=\bigl\{ f\in D_{\max}((-1,0)) : f(-1)=0,\ 
\ell^-_1\alpha^-+\ell^-_2\beta^-=0 \bigr\},
\]
\[
D(A_\nu^+)=\bigl\{ f\in D_{\max}((0,1)) : f(1)=0,\ 
\ell^+_1\alpha^+ + \ell^+_2\beta^+=0 \bigr\}.
\]
Both extensions correspond to the Laplace operator with Robin boundary conditions, so that both generate a semigroup. As an immediate consequence, for any initial condition $f^0 \in L^2(-1,1)$,
 $$e^{-t A_\nu}f^0=(e^{-t A_\nu^-}f_0 \mathbb{1}_{(-1,0)})^r+(e^{-t A_\nu^+}f_0 \mathbb{1}_{(0,1)})^l,$$
 where $^r$ denotes the extension operator by $0$ on $(0,1)$ and  $^l$ denotes the extension operator by $0$ on $(-1,0)$.

As a consequence, such extensions are of no interest for our controllability
problem. Indeed, if the control acts only on one side of the singularity (for
instance inside $(0,1)$, to fix ideas), then, by the Duhamel formula, for any $f^0 \in L^2(-1,1)$, the restriction of the corresponding solution to \eqref{1} to the interval $(-1,0)$ is given by the restriction of $(-1,0)$ to $e^{-t A_\nu^-}f_0 \mathbb{1}_{(-1,0)}$, which follows a free evolution equation and then cannot be controlled.
This is notably the case for the uncoupled extension with domain 
\[
\{\, f \in L^{2}(-1,1) ;
f_{|[0,1]} \in H^2_{\mathrm{loc}}((0,1]) \cap H^1_0(0,1),\;
f_{|[-1,0]} \in H^2_{\mathrm{loc}}([-1,0)) \cap H^1_0(-1,0),\ 
-\partial_{xx}^2 f + \frac{c_\nu}{x^2} f \in L^2(-1,1)\},
\]
that is inherited from \cite{V} and was presented in the introduction of our paper. It turns out that this extension correspond, after easy computations, to the case 
$$M_2=\left(
\begin{array}{cc}
 2 v+1 & 2 \\
 2 v+1 & 2 \\
\end{array}
\right),\,\, M_3=\left(
\begin{array}{cc}
 -2 v-1 & 2 \\
 0 & 0 \\
\end{array}
\right).$$

For this reason, in order to have controllability from one side of the singularity, one has to consider {coupled extensions} as defined in Definition \ref{def:coupled-ext}.

\end{remark}

\subsection{Non-negativity of coupled self-adjoint extensions}

In the context of our controllability problem, the non-negativity of the underlying
self-adjoint extension plays a fundamental role. Indeed, it ensures the
well-posedness of the parabolic dynamics generated by $A_\nu$ thanks to \cite[Corollary 2.4.8]{CazenaveHaraux1990}. Moreover, it seems to be some reasonable assumption for a diffusion operator. For
this reason, we wish to restrict our attention to coupled self-adjoint
extensions that are positive.
To study the positivity of our extensions, let us compute $\langle Af, f\rangle$ under homogeneous Dirichlet conditions in the following Lemma. 
\begin{lemma}
For all $f\in D_{max}$ that verifies the homogeneous Dirichlet conditions $f(\pm1)=0$, one has
\begin{equation}\label{Affp}
\langle Af,  f\rangle =  \int_{-1}^{1} \left(\partial_x f_r\right)^2(x) + \frac{c_\nu}{x^2}f_r^2(x) \,\mathrm{d}x
-\alpha^+ \beta^+ - \alpha^-\beta^-.
\end{equation}
\label{lem:easyAff}
\end{lemma}

\noindent\textbf{Proof of Lemma \ref{lem:easyAff}}.\\
Let all $f\in D_{max}$ that verifies the homogeneous Dirichlet conditions $f(\pm1)=0$. Using Lemma \ref{lem:green} and integration by parts it comes that
\begin{align*}
\langle A_n f,f\rangle
&= \int_{-1}^{1} \left(-\partial_{xx}^2 f_r + \frac{c_\nu}{x^{2}} f_r\right)(x) f(x)\,\mathrm{d}x\\
&= \int_{-1}^{1} (\partial_x f_r)^{2}(x) + \frac{c_\nu}{x^{2}} f_r^{2}(x)\,\mathrm{d}x + (-\partial_x f_r)(1) f_r(1)
       + \partial_x f_r(-1) f_r(-1)\\
       &+ [f_r,f_s](1)  - [f_r,f_s](0^+) + [f_r,f_s](0^-) - [f_r,f_s](-1).
\end{align*}
Using Lemma \ref{lem:behavior-H20}, it comes that $[f_r,f_s](0^+) = [f_r,f_s](0^-) = 0$.
Gathering the boundary terms and using $f(1)=f(-1)=0$, we get
\begin{equation}
\label{eq:2.10}
\begin{aligned}
\langle A_n f,f\rangle
&= \int_{-1}^{1} (\partial_x f_r)^{2}(x)
   + \frac{c_\nu}{x^{2}} f_r^{2}(x)\,\mathrm{d}x\\
&\quad + f_r(1)\partial_x f_s(1) - f_r(-1)\partial_x f_s(-1).
\end{aligned}
\tag{2.10}
\end{equation}
As $f(1) = f(-1) = 0$, we obtain
\[
f_r(1)\partial_x f_s(1)
= -\alpha^+\beta^+, \quad f_r(-1)\partial_x f_s(-1)
= \alpha^-\beta^-,
\]
which concludes the proof.
\qed
\begin{remark}\label{rem:pos}
From \eqref{Affp}, we see that it is not obvious to decide if an extension is non-negative or not, because of the punctual terms $\alpha^{\pm}$ and $\beta^{\pm}$. 
A good way to ensure that we have a non-negative extension is to choose our extension in such a way that 
\begin{equation}\label{surp}\alpha^+\beta^++\alpha^-\beta^-=0.\end{equation}
As already remarked in \cite[Section 2.4]{Morancey2015}, there is one natural choice in order to ensure \eqref{surp}, which corresponds to the specific extension studied in this article, whose
domain is defined at \eqref{1.5}. It corresponds exactly to the coupled
extension associated with the identity matrix $M=I_2$ in Proposition \ref{propcoup}. 
In this case, as explained \cite[Section 2.4]{Morancey2015}, the transmission condition \eqref{eq:transmission-coupled}
reduces to
\[
\alpha^- = -\alpha^+,
\qquad
\beta^- = \beta^+,
\]
which verifies \eqref{surp}. It expresses a natural symmetric matching across the singularity
$x=0$.  
This choice is particularly relevant for controllability: it ensures that no
artificial reflection or decoupling occurs at the singular point.
In other words, $M=I_2$ imposes a transmission mechanism that preserves the
natural structure of the underlying Sturm–Liouville problem and leads to the
extension with a balanced coupling across the
singularity.
Notice that another choice would be to choose $M=-I_2$ in Proposition \ref{propcoup}, which would lead to the relation
\[
\alpha^- = \alpha^+,
\qquad
\beta^- = -\beta^+,
\]
for which \eqref{surp} is also verified. We expect this extension to be very close to the one studied in this article, and to lead to the same kind of controllability results.

\end{remark}

\begin{remark}\label{posrem}
Let us go back to the case of the usual Laplace operator $-\Delta_D$ on $(-1,1)$, with Dirichlet boundary conditions, which is a positive extension of $A_{1/2}$ on $D_{\max}(-1,1)$. Since the domain in this case is $$D(\Delta_D)=H^2(-1,1)\cap H^1_0(-1,1),$$ the ``singular'' and ``regular'' parts of a function $f \in D(\Delta_D)$ is given by 
$$f(x)=f_r(x)+f_s(x)= (f(x)-f(0)+f'(0)x) + (f(0)+f'(0)x).$$ One easily verifies thanks to the definitions of $\alpha^{\pm}$ and $\beta^{\pm}$ that this corresponds to choosing $M$ in Proposition \ref{propcoup} as 
$$M = \begin{pmatrix}-1&2\\0&-1 \end{pmatrix}.$$
In particular, $-\Delta_D$ is in the case $\nu=1/2$ a coupled extension, in the sense of Definition \ref{def:coupled-ext} (for $f\in D_{\max}$, it corresponds to imposing continuity conditions $f(0^+)=f(0^-)$ and $f'(0^+)=f'(0^-)$). Hence, it would be very tempting to consider the same matrix $M$ for any $\nu \in (0,1)$. However, it is not clear from \eqref{Affp} if this extension is non-negative, and we leave this question for future research.
\end{remark}

\section{Ill-posedness of \eqref{1} for $c<-\frac{1}{4}$}
\label{appen}
The goal of this appendix is to prove the following Proposition.
\begin{proposition}\label{p:app}
Assume we work in the state space $L^2(-1,1)$. Assume  that $c<-\frac{1}{4}$. Then, there does not exist any self-adjoint extension  of $\partial_{xx}-\dfrac{c}{x^2}\mathrm{Id}$, posed on $C^\infty_0((-1,1)\setminus \{0\})$, for which this operator generates a $C^0$-semigroup.
\end{proposition}
\noindent \textbf{Proof of Proposition \ref{p:app}.}
    We reason by contradiction. 
Consider $(A,D(A))$ any self-adjoint extension of $(\partial_{xx}-\dfrac{c}{x^2}\mathrm{Id},C^\infty_0((-1,1)\setminus \{0\}))$, such that $A$ generates a $C^0$-semigroup. Then, it is well-known (see \cite[Proposition, p. 91]{EN}) that $A$ is semibounded in the following sense: there exists $C>0$ such that for any $f\in D(A)$, we have 
\begin{equation}\label{bb}
\langle f,Af\rangle_{L^2(-1,1)} \leqslant C||f||_{L^2(-1,1)}^2.
    \end{equation}

\noindent In particular, for $f\in C^\infty_0((0,1))$ extended by $0$ on $(-1,0)$, \eqref{bb} becomes (after one integration by parts)

\begin{equation}\label{bb2}
-\int_{0}^{1}f'(x)^2\mathrm{d}x-c\int_{0}^1 \frac{f(x)^2}{x^2}\mathrm{d}x\leqslant C\int_{0}^1 {f(x)^2}\mathrm{d}x.
    \end{equation}

\noindent By an easy  density argument, \eqref{bb2} also holds for any $f\in H^1_0((0,1))$.  Let $\varepsilon \in (0,1)$. Let 
$f(x)=x^{1/2+\varepsilon}(1-x)$. We have $f\in H^1_0(0,1)$.
Explicit computations give 

$$\int_{0}^1 f(x)^2\mathrm{d}x=\frac{1}{4 \varepsilon ^3+18 \varepsilon ^2+26 \varepsilon +12},$$
$$\int_{0}^1 \frac{f(x)^2}{x^2}\mathrm{d}x=\frac{1}{2\varepsilon + 6\varepsilon^2 + 4\varepsilon^3},$$
$$\int_{0}^1 f'(x)^2\mathrm{d}x=\frac{2 \epsilon +1}{8 \epsilon ^2+8 \epsilon }.$$
Hence,  as $\varepsilon\rightarrow 0$, the right-hand side of \eqref{bb2} tends to $C/12$, whereas the left-hand side is equivalent to
$$\frac{1}{\varepsilon} \left (-\frac{1}{8}-\frac{c}{2} \right).$$

\noindent Since $c<-1/4$, this quantity goes to $+\infty$ as $\varepsilon\rightarrow 0$. This is in contradiction with \eqref{bb2}, which concludes the proof.
\qed

\end{document}